\documentclass{article}

\usepackage[utf8]{inputenc}
\usepackage[T1]{fontenc}
\usepackage[british]{babel}
\usepackage[style=british]{csquotes}
\usepackage[en-GB]{datetime2}
\DTMlangsetup[en-GB]{ord=omit}

\usepackage{amssymb}
\usepackage{amsmath}
\usepackage[amsmath,hyperref,thmmarks]{ntheorem}
\usepackage{mathtools}
\usepackage{hyperref}
\hypersetup{bookmarksnumbered}

\usepackage{tikz}
\usepackage{tikz-cd}
\usetikzlibrary{calc}

\usepackage[tt=false,semibold]{libertine}
\usepackage[libertine,smallerops,vvarbb]{newtxmath}
\usepackage[lining,scaled=.8]{FiraMono}
\usepackage[cal=euler,scr=boondox,frak=euler]{mathalpha}
\tikzcdset{
    arrow style=tikz,
    arrows={/tikz/line width=.5pt},
    diagrams={>={Straight Barb[scale=0.8]}}
}

\usepackage{bm}
\usepackage{cases}

\usepackage{setspace}
\newcommand{\initlengths}{%
    \setlength{\abovedisplayshortskip}{3pt plus 9pt minus 3pt}%
    \setlength{\belowdisplayshortskip}{9pt plus 9pt minus 9pt}%
    \setlength{\abovedisplayskip}{9pt plus 9pt minus 9pt}%
    \setlength{\belowdisplayskip}{9pt plus 9pt minus 9pt}%
}

\usepackage{calc}
\newlength{\dwidth}
\usepackage{titling}

\usepackage[immediate]{silence}
\usepackage{sectsty}
\DeactivateWarningFilters[temp]
\allsectionsfont{\bfseries\boldmath}

\makeatletter
\def\theorem@optheaderfont{\normalfont}
\renewtheoremstyle{plain}%
    {\item[\theorem@headerfont\hskip\labelsep ##1\ ##2\theorem@separator]}%
    {\item[\theorem@headerfont\hskip\labelsep ##1\ ##2\ {\theorem@optheaderfont(##3)}\theorem@separator]}
\newtheoremstyle{prooflike}%
    {\item[\theorem@headerfont\hskip\labelsep ##1\theorem@separator]}%
    {\item[\theorem@headerfont\hskip\labelsep ##3\theorem@separator]}
\makeatother

\theoremseparator{.}
\newtheorem{theorem}{Theorem}[section]
\newtheorem{lemma}[theorem]{Lemma}
\newtheorem{corollary}[theorem]{Corollary}
\newtheorem{proposition}[theorem]{Proposition}
\theoremsymbol{\ensuremath{\Box}}

\newtheorem{corollaryq}[theorem]{Corollary}

\theorembodyfont{\normalfont}
\theoremsymbol{\ensuremath{\lhd}}
\newtheorem{definition}[theorem]{Definition}
\newtheorem{remark}[theorem]{Remark}
\newtheorem{construction}[theorem]{Construction}

\newtheorem{notation}[theorem]{Notation}
\theoremstyle{prooflike}
\theoremheaderfont{\normalfont\itshape}
\theoremsymbol{\ensuremath{\Box}}
\newtheorem{proof}{Proof}
\theoremheaderfont{\bfseries}
\newtheorem{bproof}{Proof}

\numberwithin{equation}{section}
\newcounter{subequation}
\numberwithin{subequation}{equation}

\usepackage[labelsep=period, labelfont={bf}]{caption}
\numberwithin{figure}{section}

\usepackage{enumitem}
\setlist[enumerate]{label=\textnormal{(\roman*)}}

\usepackage[
    backend=biber,
    style=oxnum,
    giveninits,
    scnames,
    maxbibnames=4,
    sorting=nyt
]{biblatex}

\DeclareFontFamily{OMX}{MnSymbolE}{}
\DeclareSymbolFont{MnLargeSymbols}{OMX}{MnSymbolE}{m}{n}
\SetSymbolFont{MnLargeSymbols}{bold}{OMX}{MnSymbolE}{b}{n}
\DeclareFontShape{OMX}{MnSymbolE}{m}{n}{
    <-6>  MnSymbolE5
   <6-7>  MnSymbolE6
   <7-8>  MnSymbolE7
   <8-9>  MnSymbolE8
   <9-10> MnSymbolE9
  <10-12> MnSymbolE10
  <12->   MnSymbolE12
}{}
\DeclareFontShape{OMX}{MnSymbolE}{b}{n}{
    <-6>  MnSymbolE-Bold5
   <6-7>  MnSymbolE-Bold6
   <7-8>  MnSymbolE-Bold7
   <8-9>  MnSymbolE-Bold8
   <9-10> MnSymbolE-Bold9
  <10-12> MnSymbolE-Bold10
  <12->   MnSymbolE-Bold12
}{}
\DeclareMathDelimiter{[}    {\mathopen} {MnLargeSymbols}{'000}{MnLargeSymbols}{'000}
\DeclareMathDelimiter{]}    {\mathclose}{MnLargeSymbols}{'005}{MnLargeSymbols}{'005}
\DeclareMathDelimiter{\llbr}{\mathopen} {MnLargeSymbols}{'102}{MnLargeSymbols}{'102}
\DeclareMathDelimiter{\rrbr}{\mathclose}{MnLargeSymbols}{'107}{MnLargeSymbols}{'107}
                     
\makeatletter
\def\big#1{{\hbox{$\left#1\vbox to8.5\p@{}\right.\n@space$}}}
\makeatother

\newcommand{\bbC}{\mathbb{C}}
\newcommand{\bbP}{\mathbb{P}}
\newcommand{\bbQ}{\mathbb{Q}}
\newcommand{\bbR}{\mathbb{R}}
\newcommand{\bbT}{\mathbb{T}}
\newcommand{\bbZ}{\mathbb{Z}}
\newcommand{\bia}{\bm{a}}
\newcommand{\bid}{\bm{d}}
\newcommand{\bie}{\bm{e}}
\newcommand{\bif}{\bm{f}}

\newcommand{\bit}{\bm{t}}

\newcommand{\bix}{\bm{x}}
\newcommand{\biy}{\bm{y}}
\newcommand{\biz}{\bm{z}}
\newcommand{\calE}{\mathcal{E}}
\newcommand{\calF}{\mathcal{F}}

\newcommand{\calI}{\mathcal{I}}

\newcommand{\calL}{\mathcal{L}}
\newcommand{\calM}{\mathcal{M}}
\newcommand{\calO}{\mathcal{O}}
\newcommand{\calP}{\mathcal{P}}
\newcommand{\calU}{\mathcal{U}}
\newcommand{\calV}{\mathcal{V}}
\newcommand{\calX}{\mathcal{X}}
\newcommand{\ch}{\operatorname{ch}}
\newcommand{\Coh}{\textsf{Coh}}
\newcommand{\Cohb}{\acute{\textsf{C}}\textsf{oh}}
\newcommand{\Db}{\textsf{D}^{\mathrm{b}}}
\newcommand{\Dbb}{\acute{\textsf{D}}{}^{\mathrm{b}}}
\renewcommand{\dbinom}[2]{\Bigl({\setlength{\arraycolsep}{1pt}\begin{array}{cc}{#1}\\{#2}\end{array}}\Bigr)}
\renewcommand{\det}{\operatorname{det}}

\newcommand{\Eb}{\acute{\mathcal{E}}}
\newcommand{\Extc}{\mathcal{E}\mspace{-2mu}\mathit{xt}}
\newcommand{\Extcb}{\acute{\mathcal{E}}\mspace{-2mu}\mathit{xt}}
\newcommand{\frM}{\mathfrak{M}}
\newcommand{\fund}[1]{\smash{[#1]}_{\mathrm{fund}}}

\newcommand{\Gm}{\mathbb{G}_{\mathrm{m}}}

\newcommand{\id}{\mathrm{id}}
\newcommand{\im}{\operatorname{im}}
\newcommand{\inv}[1]{\smash{[#1]}_{\mathrm{inv}}}

\newcommand{\Mb}{\acute{\mathcal{M}}}
\newcommand{\Mt}{\tilde{\mathcal{M}}}
\newcommand{\Mbnpl}[1]{\acute{\mathcal{M}}^{\mspace{1mu}#1,\mspace{1mu}\smash\pl}}
\newcommand{\Mbpl}[1][]{\acute{\mathcal{M}}^{\mspace{2mu}\smash\pl#1}}

\newcommand{\Mbss}[1][]{\acute{\mathcal{M}}^{\mspace{1mu}\mathrm{ss}#1}}
\newcommand{\Mbnss}[1]{\acute{\mathcal{M}}^{\mspace{1mu}#1,\mspace{1mu}\mathrm{ss}}}
\newcommand{\Mbnst}[1]{\acute{\mathcal{M}}^{\mspace{1mu}#1,\mspace{1mu}\mathrm{st}}}
\newcommand{\Mpl}{\mathcal{M}^{\mspace{2mu}\smash\pl}}

\newcommand{\Mtpl}{\tilde{\mathcal{M}}{}^{\mspace{2mu}\smash\pl}}
\newcommand{\Mss}{\mathcal{M}^{\mspace{1mu}\smash{\mathrm{ss}}}}
\newcommand{\Mssfd}{\mathcal{M}^{\mspace{1mu}\smash{\mathrm{ss,\mspace{2mu}fd}}}}
\newcommand{\numberthis}{\addtocounter{equation}{1}\tag{\theequation}}
\newcommand{\PD}{\operatorname{PD}}
\newcommand{\Perf}{\textsf{Perf}}
\newcommand{\PF}{\mathrm{PF}}
\newcommand{\PFtilde}{\smash{\widetilde{\mathrm{PF}}}}
\newcommand{\pl}{{\smash{\mathrm{pl}}}}
\newcommand{\pr}{\mathrm{pr}}
\newcommand{\pt}{{*}}
\newcommand{\res}{\operatorname{res}}

\newcommand{\second}{\mathchoice{^{\prime\mspace{-2mu}\prime}}{^{\prime\mspace{-2mu}\prime}}{^{\prime\mspace{-1mu}\prime}}{^{\prime\mspace{-1mu}\prime}}}
\newcommand{\sumbar}{\mathop{\mathrlap{\mathchoice{\mspace{4mu}\smash{\frac{\mspace{12mu}}{}}}{\mspace{1mu}\smash{\frac{\mspace{11mu}}{}}}{\smash{\frac{\mspace{8mu}}{}}}{\smash{\frac{\mspace{8mu}}{}}}}{\sum}}}
\newcommand{\supp}{\operatorname{supp}}
\newcommand{\sym}{\mathrm{sym}}
\newcommand{\td}{\operatorname{td}}
\newcommand{\third}{\mathchoice{^{\prime\mspace{-2mu}\prime\mspace{-2mu}\prime}}{^{\prime\mspace{-2mu}\prime\mspace{-2mu}\prime}}{^{\prime\mspace{-1mu}\prime\mspace{-1mu}\prime}}{^{\prime\mspace{-1mu}\prime\mspace{-1mu}\prime}}}
\newcommand{\upe}{\mathrm{e}}
\newcommand{\upi}{\mathrm{i}}
\newcommand{\upL}{\mathrm{L}}
\newcommand{\virt}[1]{\smash{[#1]}_{\mathrm{virt}}}

\newcommand{\leftsubstack}[2][6em]{\substack{\makebox[#1][l]{\scriptsize$\begin{aligned}#2\end{aligned}$}}}

\DeclareFontFamily{OMX}{MnSymbolE}{}
\DeclareSymbolFont{MnLargeSymbols}{OMX}{MnSymbolE}{m}{n}
\DeclareFontShape{OMX}{MnSymbolE}{m}{n}{
    <-6>  MnSymbolE5
   <6-7>  MnSymbolE6
   <7-8>  MnSymbolE7
   <8-9>  MnSymbolE8
   <9-10> MnSymbolE9
  <10-12> MnSymbolE10
  <12->   MnSymbolE12
}{}
\DeclareMathDelimiter{\llbr}{\mathopen} {MnLargeSymbols}{'102}{MnLargeSymbols}{'102}
\DeclareMathDelimiter{\rrbr}{\mathclose}{MnLargeSymbols}{'107}{MnLargeSymbols}{'107}

\makeatletter
\def\blfootnote{\xdef\@thefnmark{}\@footnotetext}
\makeatother

\title{Counting sheaves on curves}
\author{Chenjing Bu}
\date{}

\begin{document}

\initlengths

\maketitle

\begin{abstract}
    We compute Joyce's~\cite{Joyce2021} enumerative invariants
    $\inv{\Mss_{\smash{(r,d)}}}$
    for semistable rank~$r$ degree~$d$ coherent sheaves
    on a complex projective curve.
    These invariants are a generalization of the fundamental class
    of the moduli of semistable sheaves.
    We express the invariants as a regularized sum,
    which is a way to assign finite values to divergent series,
    and we obtain explicit expressions for the invariants.
    
    From these invariants,
    one can extract cohomology pairings on the moduli
    of semistable sheaves.
    When $r$ and $d$ are coprime,
    formulae for such pairings were found by Witten~\cite{Witten1992}
    and proved by Jeffrey and Kirwan~\cite{jeffrey-kirwan-2}.
    Our results provide a new point of view on this classical problem,
    and can be seen as a generalization of this
    to the case when $r$ and $d$ are not coprime.
\end{abstract}

\tableofcontents

\section{Introduction}
\label{sect-intro}

Let $X$ be a smooth, projective algebraic curve over $\bbC$.
Let $\Mss_{(r, d)}$ be the moduli of semistable rank~$r$, degree~$d$
coherent sheaves on~$X$,
and $\Mssfd_{(r, d)}$ the fixed determinant version,
consisting of coherent sheaves with a
chosen determinant line bundle $L \to X$ of degree~$d$.

It has been a classical problem to determine the structure
of the cohomology of these moduli spaces,
starting from the seminal work of
Atiyah--Bott~\cite{AtiyahBott1983},
who described a set of generators
for the cohomology ring of $\Mssfd_{(r,d)}$~%
\cite[Theorem~9.11]{AtiyahBott1983} when $r$ and $d$ are coprime.
It then became an important question to study
cohomology pairings on these moduli spaces.

When $r$ and~$d$ are coprime,
the spaces $\Mss_{(r, d)}$ and $\Mssfd_{(r,d)}$
are smooth, projective $\bbC$-varieties.
Witten~\cite[(5.21)]{Witten1992} obtained a formula for
intersection pairings on the moduli space $\Mssfd_{(r,d)}$,
employing methods from physics.
His formula was subsequently proved by
Jeffrey--Kirwan,
who expressed the intersection pairing
as a repeated residue of a meromorphic function~%
\cite[Theorem~9.12]{jeffrey-kirwan-2},
and established its equivalence to Witten's formula.

On the other hand, when $r$ and~$d$ are not coprime,
the moduli spaces $\Mss_{(r,d)}$ and $\Mssfd_{(r,d)}$ can be singular,
and intersection pairings are not directly defined.
In this case, Jeffrey \emph{et~al.}~\cite[Theorem~20]{JKKW2006}
computed cohomology pairings on
a partial desingularization of $\Mss_{(r, d)}$,
generalizing the aforementioned work of Jeffrey--Kirwan.

A new approach to this classical problem
was facilitated by the recent work of Joyce~\cite{Joyce2021},
in which a generalized notion of enumerative invariants was introduced.
For counting sheaves on curves,
Joyce's invariants are homology classes denoted by $\inv{\Mss_{(r,d)}}$\,,
which coincide with the fundamental classes $\fund{\Mss_{(r,d)}}$
when $r$ and~$d$ are coprime.
Joyce also defined his invariants for stable pairs,
and established wall-crossing formulae
for these invariants when varying the stability condition.
Applying these wall-crossing formulae
leads to a method to compute the invariants $\inv{\Mss_{(r,d)}}$\,,
from which one can extract cohomology pairings on
$\Mss_{(r,d)}$ and $\Mssfd_{(r,d)}$\,.

In this paper, we compute Joyce's invariants
$\inv{\Mss_{(r,d)}}$ for all possible $r$ and~$d$,
as well as the fixed determinant version $\inv{\Mssfd_{(r,d)}}$\,,
and we obtain explicit expressions for these invariants
in~\S\S\ref{sect-main-results}--\ref{sect-fixed-determinant}.
We show in~\S\ref{sect-comparison}
that when $r$ and~$d$ are coprime,
at least in the simple case of the symplectic volume,
our formula for intersection pairings extracted from the invariants
$\inv{\Mssfd_{(r,d)}}$
is equivalent to Jeffrey--Kirwan's repeated residue formula,
and hence, to Witten's formula.

We will express the invariants $\inv{\Mss_{(r,d)}}$
as a \emph{regularized sum} over a lattice.
The regularized sum, denoted by $\sumbar$,
is a way to assign finite values to divergent series.
Roughly, this is defined by imposing the relation
\begin{equation}
    \sumbar_{n = 0}^{\infty} a^n = \frac{1}{1 - a}
\end{equation}
for any expression $a \neq 1$;
the precise definition will be given in \S\ref{sect-regularized-sum}.
The way that the regularized sum appears in our situation
is that we often have a geometric series of meromorphic functions,
and we take the residue of the regularized sum of this series,
as if we were summing over the residues of the individual terms.
For example, if we consider the regularized sum
\begin{equation}
    \sumbar_{n = 0}^{\infty} \frac{\upe^{n z}}{z^2} = \frac{1}{z^2 (1 - \upe^z)},
\end{equation}
then the $z = 0$ residue of the right-hand side is $-1/12$,
but we think of this as summing over the
residues on the left-hand side, $1 + 2 + 3 + 4 + \cdots$.
This can be seen as a model of
how regularized sums work in our formalism.

Our main result,
Theorem~\ref{thm-inv-mss-main-as-reg-sum},
states that when $r \geq 1$,
the invariant $\inv{\Mss_{(r,d)}}$ can be expressed as a regularized sum
\begin{align*}
    \numberthis
    \label{eq-intro-main-1}
    & \inv{\Mss_{(r,d)}} = 
    \frac{1}{r} \cdot {}
    \\ & \hspace{1em}
    \sumbar_{ \leftsubstack[6em]{ 
        \\[-2ex]
        & d = d_1 + \cdots + d_{r} \\[-.6ex]
        & (d_1 + \cdots + d_{i})/i \leq d/r, \ i = 1, \dotsc, r-1 \\[-.6ex]
        & \text{with $m$ equalities} 
    } } \mspace{-9mu}
    \frac{1}{m+1} \cdot
    \bigl[ \bigl[ \dotsc \bigl[
        \fund{\Mss_{(1,d_1)}} \, , \ 
        \fund{\Mss_{(1,d_2)}} \bigr] , \ 
        \dotsc \bigr] , \ 
        \fund{\Mss_{(1,d_r)}}
    \bigr],
\end{align*}
where the Lie brackets come from
Joyce's \cites[\S4]{Joyce2021}{joyce-hall}
vertex algebra structure on the homology of moduli stacks,
which we will discuss in detail in~\S\ref{sect-background-sheaves} below.
Note that these Lie brackets are defined as
residues of meromorphic functions obtained as vertex operations,
and the regularized sum here is to be understood as
the residue of the regularized sum of those meromorphic functions.

In particular, when $r$ and~$d$ are coprime,
\eqref{eq-intro-main-1} simplifies to
\begin{align*}
    \numberthis
    \label{eq-intro-main-2}
    & \fund{\Mss_{(r,d)}} = 
    \frac{1}{r} \cdot {}
    \\[1ex] & \hspace{1em}
    \sumbar_{ \leftsubstack[6em]{ 
        \\[-2ex]
        & d = d_1 + \cdots + d_{r} \\[-.6ex]
        & (d_1 + \cdots + d_{i})/i < d/r, \ i = 1, \dotsc, r-1
    } } \mspace{-9mu}
    \bigl[ \bigl[ \dotsc \bigl[
        \fund{\Mss_{(1,d_1)}} \, , \ 
        \fund{\Mss_{(1,d_2)}} \bigr] , \ 
        \dotsc \bigr] , \ 
        \fund{\Mss_{(1,d_r)}}
    \bigr].
\end{align*}

From~\eqref{eq-intro-main-1},
one can extract cohomology pairings on the moduli spaces
$\Mss_{(r,d)}$ as well as the fixed determinant moduli spaces $\Mssfd_{(r,d)}$.
Explicit formulae for these pairings will be given in Theorems~\ref{thm-inv-mss-main}
and~\ref{thm-inv-mssfd-pairing}, respectively.

In particular, one can explicitly compute
the symplectic volume of $\Mssfd_{(r,d)}$,
in the sense of Definition~\ref{def-symp-vol}.
For example, when $r = 2$, the volume of $\Mssfd_{(2,d)}$
is given by the formula
\begin{equation}
    \mathrm{vol} (\Mssfd_{(2,d)}) =
    \frac{(-1)^{g+d} \, 2^{g-1}}{(2g-2)!} \, B_{2g-2} \Bigl( \Bigl\{ \frac{d}{2} \Bigr\} \Bigr),
\end{equation}
which is a special case of~\eqref{eq-pairing-mssfd-2-d},
where $g$ is the genus of the curve,
$B_n (-)$ denotes the $n$-th Bernoulli polynomial,
and $\{ d/2 \}$ denotes the fractional part of $d/2$.
When $d$ is odd, this agrees with the formulae of
Witten~\cites[(3.27)]{Witten1991} and
Jeffrey--Kirwan~\cite[(2.9)]{jeffrey-kirwan-2}.
When $d$ is even, this agrees with Witten's formula~%
\cite[(3.11)]{Witten1991}.



A question that remains open is that,
in the case when $r$ and $d$ are not coprime,
whether our formula for the invariants~$\inv{\Mss_{(r,d)}}$
agrees with the cohomology pairings on the
partial desingularization of $\Mss_{(r,d)}$
given by Jeffrey \emph{et~al.}~\cites[Theorem~20]{JKKW2006}.

Another interesting question is
whether one could relate our formula with the Verlinde formula,
as in Jeffrey--Kirwan~\cite[Theorem~11.2]{jeffrey-kirwan-2},
which is a formula for the dimension of
$H^0 (\smash{\Mss_{(r,d)}}, \calL^{\otimes k})$,
where $\calL \to \smash{\Mss_{(r,d)}}$ is a certain line bundle, and $k > 0$.
In the coprime case,
Jeffrey--Kirwan proved the Verlinde formula
using intersection pairings on the moduli space
and the Riemann--Roch formula,
and one might expect to obtain a non-coprime version
following a similar approach,
using our formula for the invariants~$\inv{\Mss_{(r,d)}}$.

This paper is organized as follows.
In \S\ref{sect-background},
we provide background material on moduli stacks of sheaves and pairs,
and their homology.
We also discuss Joyce's vertex algebra construction
on the homology of these moduli stacks.
In \S\ref{sect-rank-0} and \S\ref{sect-rank-1},
we compute the invariants $\inv{\Mss_{(0,d)}}$ and $\inv{\Mss_{(1,d)}}$,
respectively.
In \S\ref{sect-higher-rank},
we state our main results on the invariants
$\inv{\Mss_{(r,d)}}$ and $\inv{\Mssfd_{(r,d)}}$ for $r > 1$,
and compare them to Jeffrey--Kirwan's formulae.
In \S\ref{sect-elliptic-curves},
we compute the invariants for elliptic curves
using the Fourier--Mukai transform,
and we check that the result is consistent with the results of \S\ref{sect-higher-rank}.
The appendix \S\ref{sect-regularized-sum}
is devoted to defining the notion of the regularized sum.
Finally, in \S\ref{sect-proof-of-main-theorem},
we present the proof of our main results in \S\ref{sect-higher-rank}.

\subsection*{Acknowledgements}

The author is grateful to his supervisor, Dominic Joyce,
for introducing the subject of enumerative geometry to the author,
and providing many insightful comments and suggestions on this project. 
The author would also like to thank
Christopher Beem, Zhengping Gui, Frances Kirwan, Si Li, Balázs Szendrői, and Kai Xu
for helpful discussions.
Finally, the author would like to thank
the anonymous referee for giving numerous valuable suggestions.

This work was done during the author's PhD programme funded by the
Mathematical Institute, University of Oxford.

\section{Background material}
\label{sect-background}

\subsection{Moduli stacks of sheaves}
\label{sect-background-sheaves}

We review some background on the
moduli stack of coherent sheaves on an algebraic curve,
and we state Joyce's~\cites[\S4]{Joyce2021}{joyce-hall}
vertex algebra construction
on the homology of this moduli stack.
Most of this material can be found in Joyce~\cite[\S7]{Joyce2021}.

Let $X$ be a smooth, projective curve over $\bbC$, of genus $g$.
For $k = 0, 1, 2$,
let $b^k$ be the $k$-th Betti number of $X$,
so that $b^0 = b^2 = 1$ and $b^1 = 2g$.
Choose a basis $( \epsilon_{j,k} )_{j=1, \dotsc, b^k}$
for $H^k (X; \bbQ)$, so that
\begin{align*}
    \epsilon_{1,0} & = 1, \\
    \textstyle
    \int_X \epsilon_{1,2} & = 1, \\
    \textstyle
    \int_X \epsilon_{j,1} \cup \epsilon_{j',1} & =
    \begin{cases}
        1, & j' = j + g, \\
        -1, & j = j' + g, \\
        0, & \text{otherwise}.
    \end{cases}
\end{align*}
For $k = 0,1,2$, let $(e_{j,k})_{j=1,\dotsc,b^k}$
be the dual basis for $H_k (X; \bbQ)$, so that
\[
    \epsilon_{j,k} \cdot e_{j',k'} =
    \begin{cases}
        1, & (j,k) = (j',k'), \\
        0, & \text{otherwise}.
    \end{cases}
\]
Write
\begin{equation}
    J = \{ (j, k) \mid k = 0,1,2, \ j = 1,\dotsc,b^k \}.
\end{equation}

For $(j, k) \in J$, define
\begin{equation}
    M_{j,k}^{j',k'} =
    \int_X \epsilon_{j,k} \cup \epsilon_{j',k'} \cup \td (X),
\end{equation}
where $\td (X) = 1 + (1-g) \, \epsilon_{1,2}$
is the Todd class of $X$.
Explicitly, we have
\begin{align*}
    M_{1,0}^{1,0} & = 1 - g, \\
    M_{1,0}^{1,2} = M_{1,2}^{1,0} & = 1, \\
    M_{j,1}^{j+g,1} = -M_{j+g,1}^{j,1} & = 1 \quad (j = 1, \dotsc, g), \\
    M_{j,k}^{j',k'} & = 0 \quad
    \text{in all other cases}.
\end{align*}

Let $\Coh (X)$ be the abelian category
of coherent sheaves on $X$,
and let $\Db (X) = \Db \Coh (X)$ be its bounded derived category.

Let $K (X) = K^{\mathrm{num}} (\Coh (X))$
be the \emph{numerical Grothendieck group} of $\Coh (X)$,
as in~\cite[Example~1.4]{Joyce2021}. In this case, we have
\begin{equation}
    K (X) = \{ (r, d) \mid r, d \in \bbZ \},
\end{equation}
where $r$ and $d$ represent the rank and the degree, respectively.
Let
\begin{equation}
    \label{eq-def-cx}
    C (X) = \{ (r, d) \in K (X) \mid r > 0 \text{ or } (r = 0 \text{ and } d > 0) \}
\end{equation}
be the \emph{positive cone}.

By To\"en--Vaqui\'e \cite[Theorem~3.6]{toen-vaquie},
there is a derived moduli stack $\frM$
of objects of $\Db (X)$,
which is a derived stack over~$\bbC$.
Let $\calM$ denote its classical truncation,
which is a higher stack over~$\bbC$.
It comes with a universal perfect complex
\begin{equation}
    \calU \longrightarrow X \times \calM.
\end{equation}

Note that we choose to work with the full moduli stack
of complexes of coherent sheaves on $X$,
as opposed to the moduli of coherent sheaves,
since the homology of the former is much simpler than that of the latter,
as in~\eqref{eq-coh-free}--\eqref{eq-homol-free} below.

The moduli stack $\calM$ can be decomposed as
\begin{equation}
    \calM = \coprod_{\alpha \in K(X)} \calM_{\alpha}\,,
\end{equation}
where each $\calM_{\alpha}$ is a connected component of $\calM$,
and can be seen as the moduli stack of complexes of the class $\alpha$.
Let
\begin{equation}
    \calU_\alpha \longrightarrow X \times \calM_\alpha
\end{equation}
be the restriction of $\calU$ to
$X \times \calM_\alpha$\,.

As in Joyce~\cite[Definition~7.35]{Joyce2021},
we consider the $\Gm$-\emph{rigidification}
\begin{equation}
    \Mpl = \calM / [{*}/\Gm] \, ,
\end{equation}
where $[{*}/\Gm]$ acts on $\calM$
as scalar automorphisms of sheaves.
This removes the scalar automorphisms
from stabilizer groups of $\calM$.
We have a natural projection
\begin{equation}
    \Pi^\pl \colon
    \calM \longrightarrow \Mpl,
\end{equation}
which is a principal $[{*}/\Gm]$-bundle
away from the point $0 \in \Mpl$.
We also have the decomposition into connected components as
\begin{equation}
    \Mpl = \coprod_{\alpha \in K(X)} \Mpl_{\alpha}\,.
\end{equation}

Next, we wish to study the \emph{ordinary homology} and \emph{cohomology}
of stacks, as in Joyce~\cite[Definition~2.2]{Joyce2021}.
This is defined as follows.
For a higher stack~$\calX$ over $\bbC$,
as in Simpson~\cite[\S1]{Simpson1996} and Blanc~\cite[\S3]{Blanc2016},
one can define the \emph{topological realization} of~$\calX$,
which is a topological space $|\calX|$.
The topological realization is $\mathbb{A}^1$-homotopy invariant,
meaning that it sends $\mathbb{A}^1$-homotopy equivalences
to homotopy equivalences.

For example, the topological realization of a $\bbC$-variety
coincides with the underlying topological space of its analytification;
the topological realization of $[*/G]$,
where $G$ is an algebraic group over $\bbC$,
is the classifying space of the analytification of $G$.

As in Joyce~\cite[Example~2.35]{joyce-hall},
or as can be deduced from Blanc~\cite[Theorems~4.5 and~4.21]{Blanc2016},
there is a natural map from the algebraic $K$-theory of $\calX$
to the topological $K$-theory of $|\calX|$.
This implies that one can define Chern classes for
perfect complexes on $\calX$, as cohomology classes of $|\calX|$.

In the following,
whenever we mention the \emph{ordinary homology} or \emph{cohomology}
of~$\calX$, we always refer to the homology or cohomology
of its topological realization.

For $\alpha \in K(X)$, $(j, k) \in J$, and $l \geq k/2$, let
\begin{equation}
    \label{eq-def-sjkl}
    S_{\alpha; j,k,l} =
    \ch_l (\calU_{\alpha}) \setminus e_{j,k} \quad
    \in H^{2l-k} (\calM_\alpha; \bbQ),
\end{equation}
where
\begin{equation}
    {\setminus} \colon
    H^\bullet (X \times \calM; \bbQ) \otimes H_{-\bullet} (X; \bbQ)
    \longrightarrow H^\bullet (\calM; \bbQ)
\end{equation}
denotes the slant product.
When $k$ is even and $l = k/2$,
the variable $S_{\alpha;j,k,k/2}$ is a number,
and is denoted by $\alpha_{j,k}$\,,
so that if $\alpha = (r,d)$, then
\begin{equation}
    \alpha_{1,0} = r, \quad
    \alpha_{1,2} = d.
\end{equation}

By a result of Gross~\cite[Theorem~4.15]{gross},
there is an isomorphism
\begin{equation}
    \label{eq-coh-free}
    H^\bullet (\calM_\alpha; \bbQ) \simeq
    \bbQ [S_{\alpha; j,k,l} : (j,k) \in J, \ l > k/2],
\end{equation}
where the right-hand side is a graded commutative
(or supercommutative) polynomial ring,
where $\deg S_{\alpha; j,k,l} = 2l-k$,
with no extra relations on these variables.

One can define an isomorphism
\begin{equation}
    \label{eq-homol-free}
    H_\bullet (\calM_\alpha; \bbQ) \simeq
    \upe^{\alpha} \cdot
    \bbQ [s_{j,k,l} : (j,k) \in J, \ l > k/2],
\end{equation}
where $\upe^{\alpha}$ is a formal symbol,
with $\deg s_{j,k,l} = 2l-k$, so that we have the dual pairing
\begin{align*}
    \numberthis
    \label{eq-dual-pairing-sheaf}
    & \biggl(
        \prod_{ \substack{ (j, k) \in J \\ l > k/2 } }
        S_{\alpha; j,k,l}^{m_{j,k,l}}
    \biggr) \cdot
    \biggl(
        \upe^\alpha \cdot
        \prod_{ \substack{ (j, k) \in J \\ l > k/2 } }
        s_{j,k,l}^{m'_{j,k,l}}
    \biggr)
    \\ &
    = \begin{cases}
        \displaystyle
        \prod_{ \substack{
            (j, k), (j', k') \in J \\
            l > k/2, \ l' > k'/2, \\
            (j',k',l') \prec (j,k,l) \\
            k, k' \text{ odd}
        } } (-1)^{m_{j,k,l} \, m'_{j',k',l'}}
        \cdot
        \prod_{ \substack{ (j, k) \in J \\ l > k/2 } }
        m_{j,k,l}!, &
        \text{if } m_{j,k,l} = m'_{j,k,l} \text{ for all } j,k,l, \\
        0, &
        \text{otherwise},
    \end{cases}
\end{align*}
where the order of the products is determined by
the total order $\prec$
on the set of all $(j,k,l)$ with $(j,k) \in J$ and $l>k/2$,
given by $(j',k',l') \prec (j,k,l)$
if $l'<l$, or $l'=l$ and $k'<k$, or
$l'=l$, $k'=k$ and $j'<j$.

Define the \emph{Ext complex},
$\Extc \to \calM \times \calM$, by
\begin{equation}
    \Extc =
    (\pr_{23})_* \, \bigl(
        \pr_{12}^* (\calU^\vee) \otimes
        \pr_{13}^* (\calU)
    \bigr),
\end{equation}
where we are using derived versions of the
pushforward, pullback, dual, and tensor product functors,
and $\pr_{ij}$ denotes the projection from
$X \times \calM \times \calM$
to the product of the $i$-th and $j$-th factors.
Define a perfect complex
$\calE \to \calM \times \calM$ by
\begin{equation}
    \calE = \Extc^\vee.
\end{equation}

For $\alpha, \alpha' \in K (X)$,
let $\Extc_{\alpha,\alpha'}$ and $\calE_{\alpha,\alpha'}$
be the restriction of $\Extc$ and  $\calE$
to the component $\calM_\alpha \times \calM_{\alpha'}$.
Then, by the Riemann--Roch formula, we have
\begin{align*}
    & \ch (\Extc_{\alpha, \alpha'}) \\
    = {} &
    \int_X {} \ \biggl(
        \sum_{ \substack{ (j,k) \in J \\ l \geq k/2 } }
        (-1)^{l} \, \epsilon_{j,k} \boxtimes S_{\alpha;j,k,l} \boxtimes 1
    \biggr) \cdot \biggl(
        \sum_{ \substack{ (j',k') \in J \\ l' \geq k'/2 } }
        \epsilon_{j',k'} \boxtimes 1 \boxtimes S_{\alpha';j',k',l'}
    \biggr) \cdot \td (X)
    \\
    = {} &
    \sum_{ \substack{
        (j, k), \, (j', k') \in J \\
        l \geq k/2, \, l' \geq k'/2
    } }
    (-1)^l \, M_{j,k}^{j',k'} \,
    S_{\alpha; j,k,l} \boxtimes
    S_{\alpha'; j',k',l'} \, ,
    \numberthis
    \label{eq-ch-ext}
\end{align*}
so that
\begin{equation}
    \ch (\calE_{\alpha, \alpha'}) =
    \sum_{ \substack{
        (j, k), \, (j', k') \in J \\
        l \geq k/2, \, l' \geq k'/2
    } }
    (-1)^{l' - (k+k')/2} \, M_{j,k}^{j',k'} \,
    S_{\alpha; j,k,l} \boxtimes
    S_{\alpha'; j',k',l'} \, .
    \label{eq-ch-cale}
\end{equation}
The rank of $\Extc_{\alpha, \alpha'}$, or equivalently
$\calE_{\alpha, \alpha'}$,
is thus given by
\begin{align*}
    \chi (\alpha, \alpha')
    & = \sum_{ \substack{
        (j, k), (j', k') \in J \\
        k, k' \text{ even}
    } }
    (-1)^{k/2} \, M_{j,k}^{j',k'} \,
    \alpha_{j,k} \,
    \alpha'_{j',k'} \\
    & = (1 - g) \, r r' + r d' - d r',
    \numberthis
    \label{eq-def-chi}
\end{align*}
where we write $\alpha = (r, d)$ and $\alpha' = (r', d')$.

For $\alpha \in K (X)$, let
\begin{align}
    \hat{H}_\bullet (\calM_\alpha; \bbQ) & =
    H_{\bullet - 2 \chi (\alpha, \alpha)} (\calM_\alpha; \bbQ), \\
    \check{H}_\bullet (\Mpl_\alpha; \bbQ) & =
    H_{\bullet + 2 - 2 \chi (\alpha, \alpha)} (\Mpl_\alpha; \bbQ),
\end{align}
and write
\begin{align}
    \hat{H}_\bullet (\calM; \bbQ) & =
    \bigoplus_{\alpha \in K (X)}
    \hat{H}_\bullet (\calM_\alpha; \bbQ), \\
    \check{H}_\bullet (\Mpl_{>0}; \bbQ) & =
    \bigoplus_{\alpha \in C (X)}
    \check{H}_\bullet (\Mpl_\alpha; \bbQ).
\end{align}
In the following, 
$\hat{H}_\bullet (\calM; \bbQ)$ and 
$\check{H}_\bullet (\Mpl_{>0}; \bbQ)$ will acquire structures of
a graded vertex algebra and a graded Lie algebra,
respectively, and the shifts are necessary
to make the grading compatible with those algebraic structures.

Let $\alpha \in C (X)$,
and let $\mu$ denote slope stability of sheaves on $X$.
We have the open substack of $\mu$-semistable coherent sheaves
of class $\alpha$,
\begin{equation}
    \Mss_{\alpha} =
    \Mss_{\alpha} (\mu) \subset 
    \Mpl_\alpha\,.
\end{equation}
This is a smooth projective variety when
$\alpha = (r,d)$ with $r$ and $d$ coprime,
essentially by the classical result of
Narasimhan--Seshadri~\cite{narasimhan-seshadri}.

For $\alpha \in C (X)$,
Joyce \cite[Theorem~7.63]{Joyce2021} defined
\emph{enumerative invariants}
\begin{equation}
    \inv{\Mss_{\alpha}} \in \check{H}_0 (\Mpl_{\alpha}; \bbQ),
\end{equation}
which coincide with the fundamental classes $\fund{\Mss_{\alpha}}$
when $r$ and $d$ are coprime.
The actual homological degree of this class is the expected dimension of $\Mss_{\alpha}$\,, 
which is
\begin{equation}
    \operatorname{vdim}_{\bbR} \Mss_{\alpha} =
    2 - 2 \chi (\alpha, \alpha) =
    2(g-1)r^2 + 2.
\end{equation}

As in Joyce~\cites[\S4]{Joyce2021}{joyce-hall},
the homology $\hat{H}_\bullet (\calM; \bbQ)$
carries the structure of a \emph{graded vertex algebra},
in the sense of~\cite[Definition~4.1]{Joyce2021};
see also Frenkel--Ben-Zvi~\cite[Definition~1.3.1, Remark~1.3.2~(1)]{FrenkelBenZvi2004}.
We state Joyce's vertex algebra construction
in this situation as follows.

\begin{theorem}
    \label{thm-va-sheaves}
    Denote by
    \begin{align}
        \odot \colon [*/\Gm] \times \calM
        & \longrightarrow \calM, \\
        \oplus \colon \calM \times \calM
        & \longrightarrow \calM
    \end{align}
    the maps given by scalar multiplication and taking the direct sum,
    respectively.

    There is a graded vertex algebra structure
    on the graded vector space $\hat{H}_\bullet (\calM; \bbQ)$,
    defined as follows.

    \begin{enumerate}
        \item 
            The unit element $1 \in \hat{H}_0 (\calM; \bbQ)$
            is given by the class 
            \begin{equation}
                \label{eq-va-sheaf-unit}
                1 = [0] \in H_0 (\calM_{(0, 0)}; \bbQ),
            \end{equation}
            where $0 \in \calM_{(0, 0)}$ is the point
            corresponding to the zero sheaf.
        
        \item 
            The translation operator~$D$ is defined as follows.
            For $A \in H_{\bullet} (\calM_\alpha; \bbQ)$, we have
            \begin{equation}
                \label{eq-va-sheaf-transl}
                D (A) = \odot_* (t \boxtimes A) \in H_{\bullet+2} (\calM_\alpha; \bbQ),
            \end{equation}
            where $t \in H_2 ([*/\Gm]; \bbQ) \simeq H_2 (\mathrm{BU} (1); \bbQ)$
            is the dual of the first Chern class.
            
        \item 
            For $A \in \hat{H}_\bullet (\calM_\alpha; \bbQ)$ and
            $A' \in \hat{H}_\bullet (\calM_{\alpha'}; \bbQ)$, we have
            \begin{multline}
                \label{eq-va-sheaf-prod}
                Y (A, z) \, A' =
                (-1)^{\chi (\alpha, \alpha')} \cdot
                z^{\chi (\alpha, \alpha') + \chi (\alpha', \alpha)} \cdot {} \\
                (\oplus_{\alpha, \alpha'})_* \circ
                (\upe^{z D} \otimes \id) \Bigl(
                    (A \boxtimes A') \cap
                    c_{1/z} \bigl(
                        \calE_{\alpha, \alpha'} \oplus
                        \sigma_{\smash{\alpha, \alpha'}}^* (\calE_{\alpha', \alpha})
                    \bigr)
                \Bigr),
            \end{multline}
            where $\chi$ is as in~\textnormal{\eqref{eq-def-chi},}
            $\oplus_{\alpha, \alpha'} \colon \calM_\alpha \times \calM_{\alpha'}
            \to \calM_{\alpha + \alpha'}$ is the restriction of the map~$\oplus,$
            $\sigma_{\alpha, \alpha'} \colon \calM_\alpha \times \calM_{\alpha'}
            \to \calM_{\alpha'} \times \calM_\alpha$
            is the map swapping the two factors, and
            $c_{1/z} (-) = \sum_{i \geq 0} z^{-i} \, c_i (-)$.
    \end{enumerate}
\end{theorem}

By a direct computation (see also~\cite{joyce-su}),
using the identification~\eqref{eq-homol-free},
we can rewrite~\eqref{eq-va-sheaf-unit}--\eqref{eq-va-sheaf-prod}
as follows.
For $\alpha = (r, d) \in K (X)$, \ 
$\alpha' = (r', d') \in K (X)$, \ 
$p, p' \in \bbQ [s_{j,k,l} : (j, k) \in J, \ l > k/2]$, 
we have
\begin{align*}
    \numberthis
    & 1 = \upe^{(0, 0)} \cdot 1
    \in \hat{H}_0 (\calM_{(0,0)}; \bbQ),
    \\
    \numberthis
    & D (\upe^{\alpha} \cdot p) =
    \upe^{\alpha} \cdot \biggl(
        r \, s_{1,0,1} +
        d \, s_{1,2,2} +
        \sum_{ \substack{
            (j, k) \in J \\
            l > k/2
        } }
        s_{j,k,l+1}
        \frac{\partial}{\partial s_{j,k,l}}
    \biggr) \, p,
    \\
    \numberthis
    \label{eq-vertex-algebra-sheaves}
    & Y (\upe^{\alpha} \cdot p, \, z)
    (\upe^{\alpha'} \cdot p') =
    \\ & \hspace{2em}
    \upe^{\alpha + \alpha'} \cdot
    (-1)^{\chi (\alpha, \alpha')} \,
    z^{\chi (\alpha, \alpha')
        + \chi (\alpha', \alpha)} \cdot {}
    \\ & \hspace{2em}
    \mathrm{exp} \biggl[
        z \biggl(
            r \, s_{1,0,1} +
            d \, s_{1,2,2} +
            \sum_{ \substack{
                (j, k) \in J \\
                l > k/2
            } }
            s_{j,k,l+1}
            \frac{\partial}{\partial s_{j,k,l}}
        \biggr) 
    \biggr] \circ {}
    \\ & \hspace{2em}
    \mathrm{exp} \biggl[
        \sum_{ \substack{
            l, l' \geq 0 \\
            l + l' > 0
        } }
        (-1)^l \, (2g - 2) \, ( l + l' - 1 )! \,
        z^{-(l+l')} \,
        \partial_{1,0,l} \, \partial'_{1,0,l'}
    \biggr]
    \\[-3ex] && \hspace{-6em}
    \bigl( p (s_{j,k,l}) \cdot p' (s'_{j,k,l}) \bigr) \,
    \bigg|_{s_{j,k,l} = s'_{j,k,l}}, 
\end{align*}
where
\[
    \partial_{1,0,l} =
    \begin{cases}
        r, &
        l = 0, \\
        \dfrac{\partial}{\partial s_{1,0,l}}, &
        l > 0,
    \end{cases}
    \hspace{2em}
    \partial'_{1,0,l} =
    \begin{cases}
        r', &
        l = 0, \\
        \dfrac{\partial}{\partial s'_{1,0,l}}, &
        l > 0.
    \end{cases}
\]
For $\alpha = (r,d)$,
we will frequently write
\begin{equation}
    \label{eq-def-d-r-d}
    D_{(r,d)} =
    r \, s_{1,0,1} +
    d \, s_{1,2,2} +
    \sum_{ \substack{
        (j, k) \in J \\
        l > k/2
    } }
    s_{j,k,l+1}
    \frac{\partial}{\partial s_{j,k,l}}
\end{equation}
for the translation operator $D$
restricted to the class $(r, d)$.

As in Joyce~\cite[Theorem~4.8]{Joyce2021},
for $\alpha \in C (X)$,
the projection $\calM_\alpha \to \Mpl_\alpha$ induces an isomorphism
\begin{equation}
    \label{eq-h-mpl-is-quot-by-imd}
    \check{H}_\bullet (\Mpl_\alpha; \bbQ) \simeq
    \hat{H}_{\bullet+2} (\calM_\alpha; \bbQ) / \im D.
\end{equation}
The graded vertex algebra structure on
$\hat{H}_\bullet (\calM; \bbQ)$
induces a graded Lie algebra structure on
$\hat{H}_{\bullet+2} (\calM; \bbQ) / \im D$
by a standard construction, defined by
\begin{equation}
    \label{eq-def-lie-brack-sheaf}
    \bigl[ A + \im D, \ A' + \im D \bigr] =
    \res_z Y (A, \, z) \, A' + \im D
\end{equation}
in $\hat{H}_{\bullet+2} (\calM; \bbQ) / \im D$,
where $A, A' \in \hat{H}_{\bullet+2} (\calM; \bbQ)$.
The graded Lie algebra involved here is the \emph{super} version,
meaning that a sign rule is imposed on odd variables,
as in~\cite[Definition~4.2]{Joyce2021}.

In particular, this equips the graded vector space
\begin{equation}
    \check{H}_\bullet (\Mpl_{>0}; \bbQ) \simeq
    \bigoplus_{\alpha \in C (X)} 
    \hat{H}_{\bullet+2} (\calM_\alpha; \bbQ) / \im D
\end{equation}
with the structure of a graded Lie algebra.
Moreover, the subspace
\begin{equation}
    \check{H}_0 (\Mpl_{>0}; \bbQ) =
    \bigoplus_{\alpha \in C (X)} \check{H}_0 (\Mpl_\alpha; \bbQ)
\end{equation}
is an ordinary Lie algebra,
and this is where the invariants $\inv{\Mss_{\alpha}}$ live in.

\subsection{Moduli stacks of pairs}
\label{sect-background-pairs}

Next, we discuss the notion of \emph{stable pairs}
of Joyce~\cite[\S8]{Joyce2021} and Joyce--Song~\cite[\S5.4]{JoyceSong2012},
which originates from ideas of Bradlow~\cite{Bradlow1991}
and has been extensively studied in the literature.
Stable pairs are a key ingredient in Joyce's~\cite[\S5]{Joyce2021}
definition of his enumerative invariants.
Most of the material below
can be found in \cite[\S8]{Joyce2021} or \cite{joyce-su}.

As above, we fix a smooth, projective curve $X$ over $\bbC$.

\begin{definition}
    \label{def-pair}
    Fix a line bundle $L \to X$.
    As in~\cite[Definition~8.4]{Joyce2021},
    a \emph{pair} on $X$ is a triple $(E, V, \rho)$,
    where $E \in \Coh (X)$,
    $V$ is a finite-dimensional $\bbC$-vector space,
    and $\rho \colon V \otimes_{\bbC} L \to E$
    is a morphism in $\Coh (X)$.

    A \emph{morphism of pairs} $(E, V, \rho) \to (E', V', \rho')$
    consists of a morphism $\theta \colon E \to E'$ in $\Coh (X)$
    and a linear map $\phi \colon V \to V'$, such that the diagram
    \[
        \begin{tikzcd}
            V \otimes_{\bbC} L \ar[r, "\phi \otimes \id_L"] \ar[d, "\rho"'] &
            V' \otimes_{\bbC} L \ar[d, "\rho'"] \\
            E \ar[r, "\theta"] &
            E' \rlap{ .}
        \end{tikzcd}
    \]
    commutes in $\Coh (X)$.

    Let $\Cohb (X)^L$ be the abelian category of such pairs.
\end{definition}

As in~\cite[Definition~8.4]{Joyce2021},
there is a derived category of pairs $\Dbb (X)^L$,
which can be defined as a comma $\infty$-category
involving the bounded derived $\infty$-category of vector spaces
and the bounded derived $\infty$-category of coherent sheaves.
It is a $\bbC$-linear stable $\infty$-category,
and can also be viewed as a dg-category.
It is similar to but different from
the derived category $\Db \Cohb (X)^L$.

Let
\begin{equation}
    \acute{K} (X) =
    \{ (\alpha, e) \mid \alpha \in K (X), \ e \in \bbZ \},
\end{equation}
where $\alpha$ is the class of the coherent sheaf $E$ in a pair,
and $e$ is the rank of the vector space $V$ in a pair.
Define the \emph{positive cone}
\begin{equation}
    \acute{C} (X) =
    \{ (\alpha, e) \in \acute{K} (X) \mid
    \alpha \in C (X) \text{ or } (\alpha = 0 \text{ and } e > 0) \}.
\end{equation}

As in~\cite[Definition~8.4]{Joyce2021},
one can consider the To\"en--Vaqui\'e~\cite{toen-vaquie}
derived moduli stack $\acute{\frM}^L$
of objects of $\Dbb (X)^L$.
Let $\Mb^L$ denote its classical truncation,
which is a higher stack over $\bbC$.
It comes with universal complexes
\begin{align}
    \calU & \longrightarrow X \times \Mb^L , \\
    \calV & \longrightarrow \Mb^L,
\end{align}
and a universal pair
\begin{equation}
    \rho \colon L \boxtimes \calV \longrightarrow \calU
\end{equation}
on $X \times \Mb^L$.
The universal complexes give rise to classifying maps
\begin{align}
    \Pi^{\calV} & \colon \Mb \longrightarrow \Perf , \\
    \Pi^{\calU} & \colon \Mb \longrightarrow \calM ,
\end{align}
which classify the families $\calV$
and $\calU$\,, respectively,
where $\Perf$ denotes the moduli stack
of perfect complexes on $\operatorname{Spec} \mathbb{C}$,
as in To\"en--Vaqui\'e~\cite[Definition~3.28]{toen-vaquie}.

The moduli stack $\Mb^L$ can be decomposed as
\begin{equation}
    \Mb^L = \coprod_{(\alpha, e) \in \acute{K} (X)} \Mb^L_{\alpha,e}\,,
\end{equation}
where each $\Mb^L_{\alpha,e}$ is a connected component of $\Mb^L$,
consisting of pairs of class $(\alpha, e)$. Let
\begin{equation}
    \rho_{\alpha, e} \colon
    L \boxtimes \calV_{\alpha,e} \longrightarrow
    \calU_{\alpha,e}
\end{equation}
be the restriction of the universal pair to $X \times \Mb^L_{\alpha,e}$\,.

There is also a projective linear version of $\Mb^L$,
defined as the $\Gm$-rigidification
\begin{equation}
    \Mbnpl{L} =
    \Mb^L / [{*}/\Gm].
\end{equation}
There is a projection
\begin{equation}
    \acute{\Pi}^\pl \colon
    \Mb^L \to \Mbnpl{L},
\end{equation}
which is a $[{*}/\Gm]$-bundle away from the point $0 \in \Mbnpl{L}$.
We also have the decomposition
\begin{equation}
    \Mbnpl{L} =
    \coprod_{(\alpha, e) \in \acute{K} (X)} \Mbnpl{L}_{\alpha,e}\,.
\end{equation}

\begin{notation}
    \label{ntn-n}
    Typically, we will take $L = \calO_X (-N')$,
    where $\calO_X (1)$ is an ample line bundle
    and $N'$ is a large integer.
    We will write
    \begin{equation}
        c_1 (L) = -N \epsilon_{1,2} \, ,
    \end{equation}
    where $N$ is an integer.
    We sometimes also write $\nu = N-g+1$,
    so that the Hilbert polynomial
    of a rank $r$, degree $d$ coherent sheaf can be written as
    \begin{equation}
        P_{(r,d)}^{\calO_X(1)} (N') =
        r \nu + d.
    \end{equation}
    In this case, we may sometimes drop the superscript $L$
    from our notation, and write
    $\Cohb (X)$, $\Mb$, $\Mbpl$, etc., for
    $\Cohb (X)^L$, $\Mb^L$, $\Mbnpl{L}$, etc.
\end{notation}

The map
\begin{equation}
    ( \Pi^{\calU}, \ \Pi^{\calV} ) \colon
    \Mb^L \longrightarrow \calM \times \Perf
\end{equation}
is an $\mathbb{A}^1$-homotopy equivalence,
as its composition with the map $\calM \times \Perf \to \Mb^L$
sending $(E, V)$ to the pair $(E, V, 0)$ is homotopic to the identity map on $\Mb^L$,
via the homotopy $((E, V, \rho), t) \mapsto (E, V, t \rho)$,
where $t$ goes from $0$ to $1$.
As a result, by the $\mathbb{A}^1$-homotopy invariance
of the topological realization, we have
\begin{equation}
    \begin{aligned}[b]
        H^\bullet (\Mb^L_{\alpha, e}; \bbQ)
        & \simeq
        H^\bullet (\calM_\alpha; \bbQ) \otimes H^\bullet (\Perf_e; \bbQ) \\
        & \simeq
        \bbQ [S_{\alpha; j,k,l} : (j, k) \in J, \ l > k/2] \otimes
        \bbQ [R_{e; l} : l > 0],
    \end{aligned}
\end{equation}
where $S_{\alpha; j,k,l}$ is as in \eqref{eq-def-sjkl};
$S_{\alpha; j,k,l}$ is an even (or odd) variable if $k$ is even (or odd);
$\Perf_e \subset \Perf$ is the component consisting of complexes of rank~$e$.
The generator $R_{e; l}$ is of degree $2l$, defined by 
\begin{equation}
    R_{e; l} = \ch_l (\calV_e).
\end{equation}
For simplicity, we sometimes write
$S_{j,k,l}$ for $S_{\alpha; j,k,l}$\,,
and $R_l$ for $R_{e; l}$.

Write $\acute{J} = J \sqcup \{ ({+}, 0) \}$,
where ${+}$ is a symbol.
For $(j, k) \in \acute{J}$ and $l \geq k/2$, write
\begin{equation}
    \acute{S}_{\alpha, e; j,k,l} =
    \begin{cases}
        (\Pi^{\calU}_{\alpha, e})^* \, S_{\alpha; j,k,l}, &
        (j, k) \in J, \\
        (\Pi^{\calV}_{\alpha, e})^* \, R_{e; l}, &
        (j, k) = ({+}, 0)
    \end{cases}
    \quad \in H^{2l-k} (\Mb^L_{\alpha, e}),
\end{equation}
with the convention that
$S_{\alpha; j,k,k/2} = \alpha_{j,k}$ and
$R_{e; 0} = e$,
so that
\begin{equation}
    \label{eq-pair-cohom}
    H^\bullet (\Mb^L_{\alpha, e}) \simeq
    \bbQ [\acute{S}_{\alpha, e; j,k,l} : (j, k) \in \acute{J}, \ l > k/2] .
\end{equation}

We define an isomorphism
\begin{equation}
    \label{eq-pair-homol}
    H_\bullet (\Mb^L_{\alpha, e}) \simeq
    \upe^{(\alpha, e)} \cdot
    \bbQ [\acute{s}_{j,k,l} : (j, k) \in \acute{J}, \ l > k/2]
\end{equation}
where $\upe^{(\alpha, e)}$ is a formal symbol,
with $\deg s_{j,k,l} = 2l-k$,
by the dual pairing
\begin{align*}
    \numberthis
    & \biggl(
        \prod_{ \substack{ (j, k) \in \acute{J} \\ l > k/2 } }
        \acute{S}_{\alpha,e; j,k,l}^{m_{j,k,l}}
    \biggr) \cdot
    \biggl(
        \upe^\alpha \cdot
        \prod_{ \substack{ (j, k) \in \acute{J} \\ l > k/2 } }
        \acute{s}_{j,k,l}^{m'_{j,k,l}}
    \biggr) \\
    & \quad = \begin{cases}
        \displaystyle
        \prod_{ \substack{
            (j, k), (j', k') \in \acute{J} \\
            l > k/2, \ l' > k'/2, \\
            (j',k',l') \prec (j,k,l) \\
            k, k' \text{ odd}
        } } (-1)^{m_{j,k,l} \, m'_{j',k',l'}}
        \cdot
        \prod_{ \substack{ (j, k) \in \acute{J} \\ l > k/2 } }
        m_{j,k,l}!, &
        \text{if } m_{j,k,l} = m'_{j,k,l} \text{ for all } j,k,l, \\
        0, &
        \text{otherwise},
    \end{cases}
\end{align*}
where $\prec$ is as in the paragraph below \eqref{eq-dual-pairing-sheaf}.
In particular, for $(j,k) \in J$ and $l > k/2$, we have
\begin{equation}
    \Pi^{\calU}_* (\upe^{(\alpha,e)} \cdot \acute{s}_{j,k,l}) =
    \upe^{\alpha} \cdot s_{j,k,l}\,.
\end{equation}
When no confusion will arise, we sometime write
\begin{equation}
    \upe^{(\alpha,e)} \cdot s_{j,k,l} =
    \upe^{(\alpha,e)} \cdot \acute{s}_{j,k,l}
\end{equation}
for $(j, k) \in \acute{J}$ and $l > k/2$.

Define a matrix
$\bigl( \acute{M}_{j,k}^{j',k'} \bigr)
_{(j, k) \in \acute{J}}
^{(j', k') \in \acute{J}}$
by
\begin{equation}
    \acute{M}_{j,k}^{j',k'} =
    \begin{cases}
        M_{j,k}^{j',k'}, &
        (j, k) \neq ({+}, 0) \neq (j', k'), \\
        -\int_X \ch(L^\vee) \, \epsilon_{j',k'} \td(X), &
        (j, k) = ({+}, 0) \neq (j', k'), \\
        0, &
        (j, k) \neq ({+}, 0) = (j', k'), \\
        1, &
        (j, k) = ({+}, 0) = (j', k'),
    \end{cases}
\end{equation}
where $\td (X) = 1 + (1-g) \, \epsilon_{1,2}$
is the Todd class of $X$.
Explicitly, in the situation of Notation~\ref{ntn-n},
we have
\begin{equation}
    \acute{M}_{{+},0}^{1,0} = -\nu, \qquad
    \acute{M}_{{+},0}^{j,1} = 0, \qquad
    \acute{M}_{{+},0}^{1,2} = -1.
\end{equation}

As in~\cite[Definition~8.14]{Joyce2021},
we define the Ext complex
$\Extcb \to \Mb^L \times \Mb^L$
by completing the exact triangle
\begin{equation}
    \label{eq-def-ext-acute}
    \begin{tikzcd}
        \Extcb \ar[d] \\
        (\pr_{23})_* \, \bigl(
            \pr_{12}^* (\calU^\vee) \otimes 
            \pr_{13}^* (\calU)
        \bigr) \ \oplus \ 
        \bigl(
            \pr_{23 \to 2}^* (\calV^\vee) \otimes
            \pr_{23 \to 3}^* (\calV)
        \bigr) \ar[d, "{
            \rho_{\alpha, e}^* \oplus \, {-(\rho_{\alpha', e'})_*}
        }"] \\
        (\pr_{23})_* \, \bigl(
            \pr_{12}^* ((L \boxtimes \calV)^\vee) \otimes 
            \pr_{13}^* (\calU)
        \bigr) \ar[d] \\
        \Extcb [1] \rlap{\ ,}
    \end{tikzcd}
\end{equation}
where $\smash{\pr_{ij}}$ denotes the projection from
$X \times \Mb^L \times \Mb^L$
to the product of the $i$-th and $j$-th factors,
and $\pr_{23 \to 2}$ and $\pr_{23 \to 3}$
denote the projections from
$\Mb^L \times \Mb^L$
to its factors.
The map $\rho_{\alpha,e}^*$ is defined in the obvious way,
and the map $(\rho_{\alpha',e'})_*$ is defined using the isomorphism
\begin{equation}
    \pr_{23}^* \bigl(
        \pr_{23 \to 2}^* (\calV^\vee) \otimes
        \pr_{23 \to 3}^* (\calV)
    \bigr)
    \simeq
    \pr_{12}^* ((L \boxtimes \calV)^\vee) \otimes 
    \pr_{13}^* (L \boxtimes \calV).
\end{equation}

Define a perfect complex $\Eb \to \Mb^L \times \Mb^L$ by
\begin{equation}
    \Eb = (\Extcb)^\vee.
\end{equation}
For $(\alpha, e), (\alpha', e') \in \acute{K} (X)$,
write
\begin{equation}
    \Extcb_{(\alpha, e), (\alpha', e')}
    \quad \text{and} \quad
    \Eb_{(\alpha, e), (\alpha', e')}
\end{equation}
for the restriction of
$\Extcb$ and $\Eb$, respectively,
to the component
$\Mb^L_{\alpha, e} \times \Mb^L_{\alpha', e'}$.

For $(\alpha,e), (\alpha',e') \in \acute{K} (X)$,
let $\Eb_{(\alpha, e), (\alpha', e')}$
be the restriction of $\Eb$
to the component $\Mb^L_{\alpha,e} \times \Mb^L_{\alpha',e'}$.
Then, a direct computation similar to
\eqref{eq-ch-ext}--\eqref{eq-ch-cale}
shows that the Chern character of 
$\smash{\Eb_{(\alpha, e), (\alpha', e')}}$
is given by
\begin{equation}
    \ch (\Eb_{(\alpha, e), (\alpha', e')}) =
    \sum_{ \substack{
        (j, k), \, (j', k') \in \acute{J} \\
        l \geq k/2, \, l' \geq k'/2
    } }
    (-1)^{l' - (k+k')/2} \, \acute{M}_{j,k}^{j',k'} \,
    \acute{S}_{\alpha, e; j,k,l} \boxtimes
    \acute{S}_{\alpha', e'; j',k',l'} \, .
\end{equation}
As a consequence,
the rank of $\Eb_{(\alpha, e), (\alpha', e')}$\,,
or equivalently $\Extcb_{(\alpha, e), (\alpha', e')}$\,,
is given by
\begin{equation}
    \label{eq-def-chi-acute}
    \acute{\chi} \bigl( (\alpha, e), (\alpha', e') \bigr)
    = \sum_{ \substack{
        (j, k), (j', k') \in \acute{J} \\
        k, k' \textnormal{ even}
    } }
    (-1)^{k/2} \, \acute{M}_{j,k}^{j',k'} \,
    \acute{S}_{\alpha, e; j,k,k/2} \,
    \acute{S}_{\alpha', e'; j',k',k'/2} \, .
\end{equation}
Explicitly, if we write $\alpha = (r, d)$ and $\alpha' = (r', d')$,
and use the convention of Notation~\ref{ntn-n}, then
\begin{equation}
    \acute{\chi} \bigl( ((r, d), e), ((r', d'), e') \bigr)
    = (1 - g) r r' +
    (r - e) d' - (d + \nu e) r' + e e'.
\end{equation}

For $(\alpha, e) \in \acute{K} (X)$, let
\begin{align}
    \hat{H}_\bullet (\Mb^L_{\alpha,e}; \bbQ) & =
    H_{\bullet - 2 \acute{\chi} ((\alpha, e), (\alpha, e))} (\Mb^L_{\alpha,e}; \bbQ), \\
    \check{H}_\bullet (\Mbnpl{L}_{\alpha,e}; \bbQ) & =
    H_{\bullet + 2 - 2 \acute{\chi} ((\alpha, e), (\alpha, e))} (\Mbnpl{L}_{\alpha,e}; \bbQ),
\end{align}
and write
\begin{align}
    \hat{H}_\bullet (\Mb^L; \bbQ) & =
    \bigoplus_{(\alpha, e) \in \acute{K} (X)}
    \hat{H}_\bullet (\Mb^L_{\alpha, e}; \bbQ), \\
    \check{H}_\bullet (\Mbnpl{L}_{>0}; \bbQ) & =
    \bigoplus_{(\alpha, e) \in \acute{C} (X)}
    \check{H}_\bullet (\Mbnpl{L}_{\alpha, e}; \bbQ).
\end{align}

Let $\acute{\mu}$ be the stability condition on $\Cohb (X)^L$
denoted by $\bar{\mu}^0_1$ in \cite[Definition~8.6]{Joyce2021}.
This is defined as a map
$\acute{\mu} \colon \acute{C} (X) \to (\bbR^2 \cup \{ +\infty \}, \leq)$,
where the total order $\leq$ is defined by
$(a_1, b_1) \leq (a_2, b_2)$ if and only if either $a_1 < b_2$\,,
or $a_1 = a_2$ and $b_1 \leq b_2$;
the element $+\infty$ is the greatest element.
For a class $(\alpha, e) = ((r, d), e) \in \acute{C} (X)$, we have
\begin{equation}
    \label{eq-def-mu-acute}
    \acute{\mu} (\alpha, e) = \begin{cases}
        (d / r, \, e / r), & r > 0, \\
        +\infty, & r = 0.
    \end{cases}
\end{equation}

For a class $(\alpha, e) = ((r, d), e) \in \acute{C} (X)$,
we have the substack of $\acute{\mu}$-semistable pairs
of class~$(\alpha, e)$,
\begin{equation}
    \Mbnss{L}_{\alpha, e} =
    \Mbnss{L}_{\alpha, e} (\acute{\mu}) \subset 
    \Mbnpl{L}_{\alpha, e}\,.
\end{equation}
When $r > 0$ and $e = 0, 1$,
Joyce~\cite[Theorem~8.24]{Joyce2021} defined \emph{enumerative invariants}
\begin{equation}
    \inv{\Mbnss{L}_{\alpha, e}} \in \check{H}_0 (\Mbnpl{L}_{\alpha, e}; \bbQ),
\end{equation}
which coincide with the Behrend--Fantechi \cite{BehrendFantechi1997}
virtual fundamental classes $\virt{\Mbnss{L}_{\alpha, e}}$
when the stable locus $\Mbnst{L}_{\alpha, e} (\acute{\mu})$
is equal to the semistable locus $\Mbnss{L}_{\alpha, e}  (\acute{\mu})$.
For example, this is always the case when $e = 1$.

For $\alpha \in C (X)$, the inclusion
\begin{equation}
    \Mss_\alpha \to
    \Mbnss{L}_{\alpha, 0}
\end{equation}
sending a sheaf $E$ to the pair $(E, 0, 0)$
identifies the sheaf invariant $\smash{\inv{\Mss_\alpha}}$ with the pair invariant
$\smash{\inv{\Mbnss{L}_{\alpha, 0}}}$\,.

Joyce~\cites[\S8]{Joyce2021}{joyce-hall} constructed
a graded vertex algebra structure on
the homology $\hat{H}_\bullet (\Mb^L; \bbQ)$,
which we state as follows.

\begin{theorem}
    \label{thm-va-pairs}
    Denote by
    \begin{align}
        \odot \colon [*/\Gm] \times \Mb^L
        & \longrightarrow \Mb^L, \\
        \oplus \colon \Mb^L \times \Mb^L
        & \longrightarrow \Mb^L
    \end{align}
    the maps given by scalar multiplication and taking the direct sum,
    respectively.

    There is a graded vertex algebra structure
    on the graded vector space $\hat{H}_\bullet (\Mb^L; \bbQ)$,
    defined as follows.

    \begin{enumerate}
        \item 
            The unit element $1 \in \hat{H}_0 (\Mb^L; \bbQ)$
            is given by the class 
            \begin{equation}
                \label{eq-va-pair-unit}
                1 = [0] \in H_0 (\Mb^L_{0, 0}; \bbQ),
            \end{equation}
            where $0 \in \Mb^L_{0, 0}$ is the point
            corresponding to the zero pair.
        
        \item 
            The translation operator~$\acute{D}$ is defined as follows.
            For $A \in H_{\bullet} (\Mb^L_{\alpha, e}; \bbQ)$, we have
            \begin{equation}
                \label{eq-va-pair-transl}
                \acute{D} (A) = \odot_* (t \boxtimes A) \in
                H_{\bullet+2} (\Mb^L_{\alpha, e}; \bbQ),
            \end{equation}
            where $t \in H_2 ([*/\Gm]; \bbQ) \simeq H_2 (\mathrm{BU} (1); \bbQ)$
            is the dual of the first Chern class.
            
        \item 
            For $A \in \hat{H}_\bullet (\Mb^L_{\alpha, e}; \bbQ)$ and
            $A' \in \hat{H}_\bullet (\Mb^L_{\alpha', e'}; \bbQ)$, 
            the vertex operation $\acute{Y} (A, z) \, A'$
            is given by
            \begin{align*}
                \numberthis
                \label{eq-va-pair-prod}
                & \acute{Y} (A, z) \, A' =
                (-1)^{\acute{\chi} ((\alpha, e), (\alpha', e'))} \cdot
                z^{\acute{\chi} ((\alpha, e), (\alpha', e')) +
                \acute{\chi} ((\alpha', e'), (\alpha, e))} \cdot {} \\
                & \hspace{1em} (\oplus_{(\alpha, e), (\alpha', e')})_* \circ
                (\upe^{z \acute{D}} \otimes \id) \Bigl(
                    (A \boxtimes A') \cap
                    c_{1/z} \bigl(
                        \Eb_{(\alpha, e), (\alpha', e')} \oplus
                        \sigma^*
                        (\Eb_{(\alpha', e'), (\alpha, e)})
                    \bigr)
                \Bigr),
            \end{align*}
            where $\acute{\chi}$ is as in~\textnormal{\eqref{eq-def-chi-acute},}
            $\oplus_{(\alpha, e), (\alpha', e')}$
            is the restriction of the map~$\oplus$
            to $\Mb^L_{\alpha, e} \times \Mb^L_{\alpha', e'} \, ,$
            $\sigma$~is the map swapping the two factors, and
            $c_{1/z} (-) = \sum_{i \geq 0} z^{-i} \, c_i (-)$.
    \end{enumerate}
\end{theorem}

By a direct computation (see also \cite{joyce-su}),
using the identification~\eqref{eq-pair-homol},
the vertex algebra structure can be explicitly written as follows:
\begin{align*}
    \numberthis
    & 1 =
    \upe^{(0, 0)} \cdot 1
    \in \hat{H}_0 (\Mb_{0,0}), \\
    \numberthis
    & \acute{D} (\upe^{(\alpha, e)} \cdot p) =
    \upe^{(\alpha, e)} \cdot \biggl(
        r \, s_{1,0,1} +
        d \, s_{1,2,2} +
        e \, s_{{+},0,1} +
        \sum_{ \substack{
            (j, k) \in \acute{J} \\
            l > k/2
        } }
        s_{j,k,l+1}
        \frac{\partial}{\partial s_{j,k,l}}
    \biggr) \, p, \\
    \numberthis
    \label{eq-vertex-algebra-pairs}
    & \acute{Y} (\upe^{(\alpha, e)} \cdot p, \, z)
    (\upe^{(\alpha', e')} \cdot p') =
    \\ & \hspace{2em}
    \upe^{(\alpha + \alpha', e + e')} \cdot
    (-1)^{\acute{\chi} ((\alpha, e), (\alpha', e'))} \,
    z^{\acute{\chi} ((\alpha, e), (\alpha', e'))
        + \acute{\chi} ((\alpha', e'), (\alpha, e))} \cdot {}
    \\ & \hspace{2em}
    \mathrm{exp} \biggl[
        z \biggl(
            r \, s_{1,0,1} +
            d \, s_{1,2,2} +
            e \, s_{{+},0,1} +
            \sum_{ \substack{
                (j, k) \in \acute{J} \\
                l > k/2
            } }
            s_{j,k,l+1}
            \frac{\partial}{\partial s_{j,k,l}}
        \biggr) 
    \biggr] \circ {}
    \\ & \hspace{2em}
    \mathrm{exp} \biggl[
        - \sum_{ \substack{
            (j, k), (j', k') \in \acute{J} \\
            l \geq k/2, \ l' \geq k'/2 \\
            l + l' - (k+k')/2 > 0
        } }
        (-1)^l \, \Bigl( l + l' - \frac{k+k'}{2} - 1 \Bigr)! \ 
        z^{-( l + l' - (k + k')/2 )} \cdot {}
    \\ & \hspace{4em}
        \bigl(
            \acute{M}_{j,k}^{j',k'} +
            (-1)^{kk' + (k+k')/2} \, \acute{M}_{j',k'}^{j,k}
        \bigr) \ 
        \partial_{j,k,l} \, \partial'_{j',k',l'}
    \biggr] \ 
    \bigl( p (s_{j,k,l}) \cdot p' (s'_{j,k,l}) \bigr) \,
    \bigg|_{s_{j,k,l} = s'_{j,k,l}}, \nonumber
\end{align*}
where $\alpha = (r, d) \in K (X)$, \ 
$\alpha' = (r', d') \in K (X)$, \ 
$p, p' \in \bbQ [s_{j,k,l} : (j, k) \in \acute{J}, \ l > k/2]$, and
\[
    \partial_{j,k,l} =
    \begin{cases}
        r, &
        (j, k, l) = (1, 0, 0), \\
        d, &
        (j, k, l) = (1, 2, 1), \\
        e, &
        (j, k, l) = ({+}, 0, 0), \\
        \dfrac{\partial}{\partial s_{j,k,l}}, &
        l > k/2,
    \end{cases}
    \hspace{2em}
    \partial'_{j,k,l} =
    \begin{cases}
        r', &
        (j, k, l) = (1, 0, 0), \\
        d', &
        (j, k, l) = (1, 2, 1), \\
        e', &
        (j, k, l) = ({+}, 0, 0), \\
        \dfrac{\partial}{\partial s'_{j,k,l}}, &
        l > k/2.
    \end{cases}
\]

For $(\alpha, e) \neq (0, 0)$,
the projection $\Mb^L \to \Mbnpl{L}$ induces an isomorphism
\begin{equation}
    \label{eq-homol-mbpl-quot}
    \check{H}_\bullet (\Mbnpl{L}_{\alpha, e}; \bbQ) \simeq
    \hat{H}_{\bullet+2} (\Mb^L_{\alpha, e}; \bbQ) / \im \acute{D}.
\end{equation}
The vertex algebra structure on
$\hat{H}_\bullet (\Mb^L; \bbQ)$,
as described above,
induces a graded Lie algebra structure on
$\hat{H}_\bullet (\Mb^L; \bbQ) / \im \acute{D}$, defined by
\begin{equation}
    \bigl[ A + \im \acute{D}, \ A' + \im \acute{D} \bigr] =
    \res_z \acute{Y} (A, \, z) \, A' + \im \acute{D}
\end{equation}
in $\hat{H}_{\bullet+2} (\Mb^L; \bbQ) / \im \acute{D}$,
where $A, A' \in \hat{H}_{\bullet+2} (\Mb^L; \bbQ)$.

In particular, this equips the graded vector space
\begin{equation}
    \check{H}_\bullet (\Mbnpl{L}_{>0}; \bbQ) \simeq
    \bigoplus_{(\alpha, e) \in \acute{C} (X)} 
    \hat{H}_{\bullet+2} (\Mb^L_{\alpha, e}; \bbQ) / \im \acute{D}
\end{equation}
with the structure of a graded Lie algebra.
Moreover, the subspace
\begin{equation}
    \check{H}_0 (\Mbnpl{L}_{>0}; \bbQ) =
    \bigoplus_{(\alpha, e) \in \acute{C} (X)}
    \check{H}_0 (\Mbnpl{L}_{\alpha, e}; \bbQ)
\end{equation}
is an ordinary Lie algebra,
and this is where the invariants $\inv{\Mbnss{L}_{\alpha, e}}$ live in.

\subsection{Homology of the projective linear moduli stacks}
\label{sect-def-xi}

In this section, we study the relationship between the homology of
the moduli stacks $\calM$ and $\Mpl$ considered above.
The purpose of doing so is that,
although the isomorphism~\eqref{eq-h-mpl-is-quot-by-imd}
allows us to write elements of
$\check{H}_\bullet (\Mpl_{>0}; \bbQ)$
as equivalence classes of polynomials in the variables
$s_{j,k,l}$, and similarly in the case of pairs,
we wish to obtain a more direct description
by choosing a canonical representative in each of these
equivalence classes.
This will enable us to identify each element of 
$\check{H}_\bullet (\Mpl_{>0}; \bbQ)$ with a polynomial
in the variables~$s_{j,k,l}$.

In the following, we fix a class
$\alpha = (r, d) \in C (X)$ with $r > 0$.

By \cite[Proposition~3.24]{joyce-hall},
the principal $[*/\Gm]$-bundle
\begin{equation}
    \Pi_\alpha^\pl \colon
    \calM_\alpha \longrightarrow \Mpl_\alpha
\end{equation}
is \emph{rationally trivial} in the sense of
\cite[Definition~2.26]{joyce-hall}.
That is, there exists an integer~$n$
and a locally trivial fibration
\begin{equation}
    p_\alpha \colon
    \Mtpl_\alpha \longrightarrow \Mpl_\alpha\,,
\end{equation}
with fibre $[*/\bbZ_n]$, such that
the pullback principal $[*/\Gm]$-bundle
\begin{equation}
    p_\alpha^* \, \Pi_\alpha^\pl \colon
    \Mt_\alpha =
    \calM_\alpha \underset{\Mpl_\alpha}{\times} \Mtpl_\alpha
    \longrightarrow \Mtpl_\alpha
\end{equation}
is trivial.
Note that $p_\alpha$ induces isomorphisms
on rational homology and cohomology,
so that $\Pi_\alpha^\pl$
looks like a trivial bundle
from the viewpoint of rational homology.

Let $i'_\alpha \colon \Mtpl_\alpha \to \calM_\alpha$
be the composition of a section
$\Mtpl_\alpha \to \Mt_\alpha$
and the projection $\Mt_\alpha \to \calM_\alpha$,
and let $\calU'^{\,\pl}_\alpha \to X \times \Mtpl_\alpha$
be the pullback of the universal sheaf
$\calU_\alpha \to X \times \calM_\alpha$
along the map $\id_X \times i'_\alpha$.
For a chosen point $x_0 \in X$, let
\begin{equation}
    \label{eq-def-upl}
    \calU^{\pl} =
    \calU'^{\,\pl} \otimes
    \pr_2^* \det \Bigl(
        \calU'^{\,\pl}_\alpha |_{\{x_0\} \times \Mtpl_\alpha}
    \Bigr)^\vee
\end{equation}
be the universal sheaf over $X \times \Mtpl_{\alpha}$,
so that it satisfies the following properties:
If $i_\alpha \colon \Mtpl_\alpha \to \calM_\alpha$
denotes the map classifying the sheaf $\calU^{\pl}$,
then
\begin{equation}
    \Pi_\alpha^\pl \circ i_\alpha = p_\alpha\,.
\end{equation}
Moreover, we have 
\begin{equation}
    \label{eq-ch1-upl-equals-zero}
    \ch_1 \bigl( \calU^{\pl}|_{\{x_0\} \times \Mtpl_\alpha} \bigr) = 0.
\end{equation}

From the viewpoint of rational homology,
the above construction gives a map
\begin{equation}
    (i_\alpha)_* \colon
    H_\bullet (\Mtpl_\alpha; \bbQ) \simeq
    H_\bullet (\Mpl_\alpha; \bbQ) \longrightarrow
    H_\bullet (\calM_\alpha; \bbQ),
\end{equation}
which is a section of the projection
\begin{equation}
    (\Pi_\alpha^\pl)_* \colon
    H_\bullet (\calM_\alpha; \bbQ) \longrightarrow
    H_\bullet (\Mpl_\alpha; \bbQ) \simeq
    H_\bullet (\calM_\alpha; \bbQ) / {\im D},
\end{equation}
where the latter isomorphism is by \cite[Theorem~4.8]{Joyce2021}.
Moreover, since \eqref{eq-ch1-upl-equals-zero}
implies that $(i_\alpha)^* \, S_{\alpha;1,0,1} = 0$, we have
\begin{equation}
    \im {} (i_\alpha)_* \subset
    \ker \frac{\partial}{\partial s_{1,0,1}} =
    \upe^{\alpha} \cdot
    \bbQ[s_{j,k,l} : (j, k, l) \neq (1, 0, 1)].
\end{equation}

Similarly, in the case of pairs,
if we fix a class $(\alpha, e) \in \acute{C} (X)$
with $\alpha = (r, d)$ and $r > 0$,
we can also construct a map
\begin{equation}
    (\acute{\imath}_{\alpha,e})_* \colon
    H_\bullet (\Mbpl_{\alpha,e}; \bbQ) \longrightarrow
    H_\bullet (\Mb_{\alpha,e}; \bbQ),
\end{equation}
which is a section of the projection
\begin{equation}
    (\acute{\Pi}_{\alpha,e}^\pl)_* \colon
    H_\bullet (\Mb_{\alpha,e}; \bbQ) \longrightarrow
    H_\bullet (\Mbpl_{\alpha,e}; \bbQ) \simeq
    H_\bullet (\Mb_{\alpha,e}; \bbQ) / {\im \acute{D}},
\end{equation}
satisfying
\begin{equation}
    \im {} (\acute{\imath}_{\alpha,e})_* \subset
    \ker \frac{\partial}{\partial \acute{s}_{1,0,1}} =
    \upe^{(\alpha, e)} \cdot
    \bbQ[\acute{s}_{j,k,l} : (j, k, l) \neq (1, 0, 1)].
\end{equation}

Next, we express the maps $(i_\alpha)_*$
and $(\acute{\imath}_{\alpha,e})_*$
in terms of the variables $s_{j,k,l}$ and $\acute{s}_{j,k,l}$,
using the isomorphisms \eqref{eq-homol-free} and~\eqref{eq-pair-homol}.

Define linear isomorphisms
\begin{align}
    \xi \colon
    H_\bullet (\Mpl_\alpha) \simeq
    \upe^{\alpha} \cdot \bbQ[s_{j,k,l}] / \im D
    & \longrightarrow
    \upe^{\alpha} \cdot
    \bbQ[s_{j,k,l} : (j, k, l) \neq (1, 0, 1)], \\
    \acute{\xi} \colon
    H_\bullet (\Mbpl_{\alpha, e}) \simeq
    \upe^{(\alpha, e)} \cdot \bbQ[\acute{s}_{j,k,l}] / \im \acute{D}
    & \longrightarrow
    \upe^{(\alpha, e)} \cdot
    \bbQ[\acute{s}_{j,k,l} : (j, k, l) \neq (1, 0, 1)]
\end{align}
by
\begin{align}
    \xi (\upe^{\alpha} \cdot p) & =
    \upe^{\alpha} \cdot \biggl[
        \sum_{i \geq 0}
        \frac{1}{i! \, (-r)^i} \,
        D^i \, \Bigl( \frac{\partial}{\partial s_{1,0,1}} \Bigr)^i
    \biggr] \, p, \\
    \acute{\xi} (\upe^{(\alpha, e)} \cdot p) & =
    \upe^{(\alpha, e)} \cdot \biggl[
        \sum_{i \geq 0}
        \frac{1}{i! \, (-r)^i} \,
        \acute{D}^i \, \Bigl( \frac{\partial}{\partial \acute{s}_{1,0,1}} \Bigr)^i
    \biggr] \, p,
\end{align}
where $r = \operatorname{rank} \alpha$.
That these maps are well-defined
will follow from the next lemma.

\begin{lemma}
    \label{lem-xi-well-def}
    For any $p \in \bbQ[s_{j,k,l}]$,
    the element $\xi (\upe^{\alpha} \cdot p)$
    is the unique element of
    \[
        (\upe^{\alpha} \cdot p + \im D) \ \cap \ 
        \upe^{\alpha} \cdot \bbQ[s_{j,k,l} : (j, k, l) \neq (1, 0, 1)].
    \]
    Similarly, for any $p \in \bbQ[\acute{s}_{j,k,l}]$,
    the element $\acute{\xi} (\upe^{(\alpha, e)} \cdot p)$
    is the unique element of
    \[
        (\upe^{(\alpha, e)} \cdot p + \im \acute{D}) \ \cap \ 
        \upe^{(\alpha, e)} \cdot \bbQ[\acute{s}_{j,k,l} : (j, k, l) \neq (1, 0, 1)].
    \]
\end{lemma}

\begin{proof}
    For the first part, notice that
    $\xi (\upe^{\alpha} \cdot p) \in \upe^{\alpha} \cdot p + \im D$.
    To see that
    $\xi (\upe^{\alpha} \cdot p) \in
    \upe^{\alpha} \cdot \bbQ[s_{j,k,l} : (j, k, l) \neq (1, 0, 1)]$,
    we notice that
    \begin{equation}
        \biggl[ \frac{\partial}{\partial s_{1,0,1}}, \ D \biggr] =
        \biggl[ \frac{\partial}{\partial s_{1,0,1}}, \ r s_{1,0,1} \biggr] = r,
    \end{equation}
    so that $[\partial / \partial s_{1,0,1}, \ D^i] = i r D^{i-1}$
    for all $i \geq 0$.
    Therefore,
    \begin{align*}
        & \phantom{{} = {}}
        \frac{\partial}{\partial s_{1,0,1}} \, \biggl[
            \sum_{i \geq 0}
            \frac{1}{i! \, (-r)^i} \,
            D^i \, \Bigl( \frac{\partial}{\partial s_{1,0,1}} \Bigr)^i
        \biggr] \, p \\
        & = \biggl[
            \sum_{i \geq 0}
            \frac{1}{i! \, (-r)^i} \biggl(
                D^i \, \Bigl( \frac{\partial}{\partial s_{1,0,1}} \Bigr)^{i+1} +
                i r D^{i-1} \, \Bigl( \frac{\partial}{\partial s_{1,0,1}} \Bigr)^i
            \biggr)
        \biggr] \, p \\
        & = 0.
        \numberthis
    \end{align*}
    To see the uniqueness,
    we notice that
    $(\im D) \cap \upe^{\alpha} \cdot \bbQ[s_{j,k,l} : (j, k, l) \neq (1, 0, 1)] = 0$.
    
    Finally, the second part is analogous to the first part.
\end{proof}

The maps $\xi$ and $\acute{\xi}$
induce dual isomorphisms
\begin{align}
    \Xi \colon
    \bbQ[S_{\alpha;j,k,l} : (j, k, l) \neq (1, 0, 1)]
    & \longrightarrow
    H^\bullet (\Mpl_{\alpha}) \subset
    \bbQ[S_{\alpha;j,k,l}], \\
    \acute{\Xi} \colon
    \bbQ[\acute{S}_{\alpha,e;j,k,l} : (j, k, l) \neq (1, 0, 1)]
    & \longrightarrow
    H^\bullet (\Mbpl_{\alpha, e}) \subset
    \bbQ[\acute{S}_{\alpha,e;j,k,l}],
\end{align}
given by
\begin{align}
    \Xi (p) & =
    \biggl[
        \sum_{i \geq 0}
        \frac{1}{i! \, (-r)^i} \,
        S_{\alpha; 1,0,1}^i \, (D^\vee)^i
    \biggr] \, p, \\
    \acute{\Xi} (p) & =
    \biggl[
        \sum_{i \geq 0}
        \frac{1}{i! \, (-r)^i} \,
        \acute{S}_{\alpha, e; 1,0,1}^i \, (\acute{D}^\vee)^i
    \biggr] \, p,
\end{align}
where
\begin{align}
    D^\vee & =
    r \, \frac{\partial}{\partial S_{\alpha; 1,0,1}}
    + d \, \frac{\partial}{\partial S_{\alpha; 1,2,2}}
    + \sum_{ \substack{
        (j, k) \in J \\
        l > k/2
    } }
    S_{\alpha; j,k,l} \,
    \frac{\partial}{\partial S_{\alpha; j,k,l+1}}, \\
    \acute{D}^\vee & =
    r \, \frac{\partial}{\partial \acute{S}_{\alpha, e; 1,0,1}}
    + d \, \frac{\partial}{\partial \acute{S}_{\alpha, e; 1,2,2}}
    + e \, \frac{\partial}{\partial \acute{S}_{\alpha, e; {+},0,1}}
    + \sum_{ \substack{
        (j, k) \in \acute{J} \\
        l > k/2
    } }
    \acute{S}_{\alpha, e; j,k,l} \,
    \frac{\partial}{\partial \acute{S}_{\alpha, e; j,k,l+1}},
\end{align}
where $\alpha = (r,d)$.

Lemma~\ref{lem-xi-well-def} implies the following.

\begin{corollaryq}
    \label{cor-xi-is-pullback}
    Under the identifications~\eqref{eq-homol-free} and~\eqref{eq-pair-homol},
    we have
    \begin{equation}
        \xi = (i_\alpha)_* \, , \qquad 
        \acute{\xi} = (\acute{\imath}_{\alpha,e})_* \, .
    \end{equation}
    Dually, we have
    \begin{equation}
        \Xi = (i_\alpha)^* \, , \qquad 
        \acute{\Xi} = (\acute{\imath}_{\alpha,e})^* \, .
    \end{equation}
    In particular,
    the maps $\Xi$ and $\acute{\Xi}$ are ring homomorphisms.
\end{corollaryq}

\subsection{Wall-crossing for pair invariants}

We present some consequences of
Joyce's wall-crossing formulae developed in \cite{Joyce2021}.
They express pair invariants in terms of
sheaf invariants and pair invariants of lower rank,
providing a way to compute these invariants by induction on the rank.%
\footnote{These results are due to Dominic Joyce,
but they are not yet publicly available at the time of this writing.}

\begin{definition}
    Let notation be as in \S\ref{sect-background-pairs}.
    As in \cite[Definition~8.5]{Joyce2021},
    for each $s \in \bbR_{> 0}$\,,
    define a stability condition $\acute{\mu}^s$ on $\Cohb (X)^L$,
    by defining a function
    $\acute{\mu}^s \colon \acute{C} (X) \to \bbQ \cup \{ +\infty \}$,
    where $\bbQ \cup \{ +\infty \}$ is equipped with the obvious total order,
    as follows.
    For $(\alpha, e) = ((r, d), e) \in \acute{C} (X)$, define
    \begin{equation}
        \acute{\mu} ^s (\alpha, e) = \begin{cases}
            (d + s e) / r, & r > 0, \\
            + \infty , & r = 0.
        \end{cases}
    \end{equation}
    The stability condition $\acute{\mu}$ given by \eqref{eq-def-mu-acute}
    can be seen as the limit of $\acute{\mu}^s$ as $s \to 0$.
\end{definition}

\begin{theorem}
    \label{thm-wcf-1}
    Let $(r,d) \in K(X)$ with $r>1$. Then
    \upshape
    \begin{multline}
        \inv{\Mbss_{(r,d),1}} = \mspace{-18mu}
        \sum_{ \leftsubstack[5em]{
            \\
            & (r,d)=(r_0,d_0)+\cdots+(r_m,d_m), \\[-.6ex]
            & m \geq 1, \ r_i > 0 \text{ for all } i, \\[-.6ex]
            & \text{such that }
            d_0 / r_0 < d_1 / r_1 \leq \cdots \leq d_m / r_m
        } } \mspace{-9mu}
        \frac{(-1)^{m-1}}{\prod_{i=1}^l (a_i - a_{i-1})!} \cdot
        \bigl[ \bigl[ \dotsc \bigl[ \bigl[
            \inv{\Mbss_{(r_0,d_0),1}} \, , \ 
            \inv{\Mss_{(r_1,d_1)}} \bigr] , \\[-8ex] 
            \shoveright{\inv{\Mss_{(r_2,d_2)}} \bigr] , \dotsc \bigr] , \ } \\
            \inv{\Mss_{(r_m,d_m)}}
        \bigr],
    \end{multline}
    \itshape
    where in each term of the sum, the numbers
    $0 = a_0 < \cdots < a_l = m$ are defined such that
    for any $0 < i < m$, we have
    $d_i/r_i < d_{i+1}/r_{i+1}$ if and only if $i = a_j$ for some $0 < j < l$.
    Only finitely many terms in the sum are non-zero.
\end{theorem}

\begin{proof}
    As in \cite[Theorem~8.24]{Joyce2021},
    for any $s, \tilde{s} \in \bbR_{>0}$\,, we have
    \begin{align*}
        \numberthis
        \label{eq-wcf-orig}
        & \inv{\Mbss_{(r,d),1} (\acute{\mu}^{\tilde{s}})} = \mspace{-9mu}
        \sum_{ \leftsubstack[5em]{
            \\[-2ex]
            & (r,d)=(r_1,d_1)+\cdots+(r_m,d_m), \\[-.6ex]
            & 1=e_1+\cdots+e_m, \\[-.6ex]
            & m \geq 1, \ e_i \geq 0, \ r_i > 0 \text{ for all } i, \\[-.6ex]
            & \Mbss_{(r_i,d_i),e_i} (\acute{\mu}^s)
            \neq \varnothing \text{ for all } i
        } } \mspace{-9mu}
        \tilde{U} \bigl(
            ((r_1, d_1), e_1), \dotsc,
            ((r_m, d_m), e_m);
            \acute{\mu}^s, \acute{\mu}^{\tilde{s}}
        \bigr) \cdot {} \\
        & \hspace{1em} \bigl[ \bigl[ \dotsc \bigl[ 
            \inv{\Mbss_{(r_1,d_1),e_1} (\acute{\mu}^s)} \, , \ 
            \inv{\Mbss_{(r_2,d_2),e_2} (\acute{\mu}^s)} \bigr] , \dotsc \bigr] , \ 
            \inv{\Mbss_{(r_m,d_m),e_m} (\acute{\mu}^s)}
        \bigr],
    \end{align*}
    where $\tilde{U} ({\cdots})$ are combinatorial coefficients
    defined in \cites[Theorem~3.12]{Joyce2021}.
    We choose $s < 1/r$, so that
    $\acute{\mu}^s$ and $\acute{\mu}$ are indistinguishable
    for all terms that appear in \eqref{eq-wcf-orig}.
    Also, note that if $\tilde{s} > (r N + d) / (r - 1)$,
    where $N$ is as in Notation~\ref{ntn-n},
    then $\smash{\Mbss_{(r,d),1} (\acute{\mu}^{\tilde{s}})} = \varnothing$.
    This is because for any pair $(E, V, \rho)$ in the class $((r, d), 1)$,
    either $\rho = 0$, in which case
    the pair is destabilized by the sub-object $(0, V, 0)$,
    or the image of $\rho$ is a rank $1$ subsheaf of $E$,
    which is of degree at least $-N$, and by our choice of $\tilde{s}$,
    this destabilizes $(E, V, \rho)$.
    We choose $\tilde{s} > (r - 1) (r N + d)$ for the argument below to work.
    
    As in \cite[Definition~3.10]{Joyce2021}, 
    for $\alpha_1 = (r_1, d_1)$, $\dotsc$, $\alpha_m = (r_m, d_m)$, 
    $\beta = (r_0, d_0)$, such that $r_i > 0$ for all $i$ and
    $(r,d) = \alpha_1 + \cdots + \alpha_m + \beta$,
    and for $0 \leq j \leq m$, there is another combinatorial coefficient
    \begin{align*}
        \numberthis
        & U \bigl( 
            (\alpha_1, 0), \dotsc,
            (\alpha_j, 0), (\beta, 1),
            (\alpha_{j+1}, 0), \dotsc,
            (\alpha_m, 0);
            \acute{\mu}^s, \acute{\mu}^{\tilde{s}}
        \bigr) = \\
        & \hspace{1em} \begin{cases}
            \dfrac{(-1)^{m-j}}{
                \prod_{i=1}^l (a_i - a_{i-1})! \cdot
                \prod_{i=1}^{l'} (a'_i - a'_{i-1})!
            }, & \begin{aligned}[t]
                & \text{if } d_1 / r_1 \geq \cdots \geq d_j / r_j > d_0 / r_0 \text{ and} \\
                & \quad d_0 / r_0 < d_{j+1} / r_{j+1} \leq \cdots \leq d_m / r_m \, ,
            \end{aligned} \\[4ex]
            0, & \text{otherwise},
        \end{cases}
    \end{align*}
    where the numbers
    $0 = a_0 < \cdots < a_l = j = a'_0 < \cdots < a'_{l'} = m$
    are defined such that
    for any $0 < i < j$, we have
    $d_i/r_i > d_{i+1}/r_{i+1}$ if and only if $i = a_k$ for some $0 < k < l$;
    for any $j < i < m$, we have
    $d_i/r_i < d_{i+1}/r_{i+1}$ if and only if $i = a'_k$ for some $0 < k < l'$.
    The choice of $\tilde{s}$ ensures that
    whenever $\smash{\Mbss_{\beta,1}} \neq \varnothing$,
    we always have
    $\acute{\mu}^{\tilde{s}} ((\alpha_1, 0) + \cdots + (\alpha_i, 0)) <
    \acute{\mu}^{\tilde{s}} ((\alpha_{i+1}, 0) + \cdots + (\alpha_m, 0) + (\beta, 1))$
    for all $0 \leq i \leq j$, and
    $\acute{\mu}^{\tilde{s}} ((\alpha_1, 0) + \cdots + (\alpha_i, 0) + (\beta, 1)) >
    \acute{\mu}^{\tilde{s}} ((\alpha_{i+1}, 0) + \cdots + (\alpha_m, 0))$
    for all $j \leq i \leq m$,
    since in this case,
    the slope of $\beta$ is at least $-N$,
    and its rank is at most $r-1$,
    so that the slope of the $\alpha_i$ cannot exceed
    $d + N (r - 1)$,
    but $\acute{\mu}^{\tilde{s}} ((r,d),1) = (d+\tilde{s})/r > d+N(r-1)$.
    
    As a result, one can take
    \begin{align*}
        \numberthis
        & \tilde{U} \bigl( 
            (\alpha_1, 0), \dotsc,
            (\alpha_j, 0), (\beta, 1),
            (\alpha_{j+1}, 0), \dotsc,
            (\alpha_m, 0);
            \acute{\mu}^s, \acute{\mu}^{\tilde{s}}
        \bigr) = \\
        & \hspace{8em} \begin{cases}
            \dfrac{(-1)^m}{
                \prod_{i=1}^l (a_i - a_{i-1})!
            }, & \begin{aligned}[t]
                & \text{if } j = 0 \text{ and} \\
                & \quad d_0 / r_0 < d_1 / r_1 \leq \cdots \leq d_m / r_m \, ,
            \end{aligned} \\[4ex]
            0, & \text{otherwise},
        \end{cases}
    \end{align*}
    where $a_0, \dotsc, a_l$ are as in the statement of the theorem.
    
    The wall-crossing formula \eqref{eq-wcf-orig} now reads
    \begin{multline}
        0 = \mspace{-18mu}
        \sum_{ \leftsubstack[5em]{
            \\
            & (r,d)=(r_0,d_0)+\cdots+(r_m,d_m), \\[-.6ex]
            & m \geq 0, \ r_i > 0 \text{ for all } i, \\[-.6ex]
            & \text{such that }
            d_0 / r_0 < d_1 / r_1 \leq \cdots \leq d_m / r_m
        } } \mspace{-9mu}
        \frac{(-1)^m}{\prod_{i=1}^l (a_i - a_{i-1})!} \cdot
        \bigl[ \bigl[ \dotsc \bigl[ \bigl[
            \inv{\Mbss_{(r_0,d_0),1} (\acute{\mu})} \, , \ 
            \inv{\Mbss_{(r_1,d_1),0} (\acute{\mu})} \bigr] , \\[-8ex] 
            \shoveright{\inv{\Mbss_{(r_2,d_2),0} (\acute{\mu})} \bigr] , \dotsc \bigr] , \ } \\
            \inv{\Mbss_{(r_m,d_m),0} (\acute{\mu})}
        \bigr],
    \end{multline}
    which proves the theorem.
\end{proof}

\begin{definition}
    \label{def-three-stab-cond}
    Let notation be as in \S\ref{sect-background-pairs}.
    As in Joyce--Song~\cite[\S13.1]{JoyceSong2012},
    for each $p \in \bbQ$, let $\Cohb (X)^L_{\smash{p}} \subset \Cohb (X)^L$
    be the full subcategory consisting of pairs $(E, V, \rho)$
    with $E$ semistable of slope $p$ or $E = 0$.
    Let $\acute{C} (X)_p \subset \acute{C} (X)$ be the subset of classes
    that contain objects of $\Cohb (X)^L_{\smash{p}}$\,. Then
    \begin{equation}
        \acute{C} (X)_p =
        \{ ((r, d), e) \mid r > 0, \ d/r = p, \ e \geq 0 \} \cup
        \{ ((0, 0), e) \mid e > 0 \}.
    \end{equation}
    
    Fix $p \in \bbQ$. 
    As in Joyce--Song~\cite[\S13.2]{JoyceSong2012}, for each $t \in \bbR$, 
    we define a stability condition $\acute{\mu}_{p, t}$
    on the abelian category $\Cohb (X)^L_{\smash{p}}$ as follows.
    For a class $(\alpha, e) = ((r, d), e) \in \acute{C} (X)_p$, define
    \begin{equation}
        \acute{\mu}_{p, t} (\alpha, e) = t e / r
        \in \bbQ,
    \end{equation}
    where $\bbQ$ is equipped with its standard total order.
    
    Note that the stability conditions $\acute{\mu}_{p,t}$ for $t > 0$
    are equivalent to the restriction of $\acute{\mu}$ to $\Cohb (X)^L_{\smash{p}}$\,.
    The stability condition $\acute{\mu}_{p,0}$ is trivial,
    and the stability conditions $\acute{\mu}_{p,t}$ for $t < 0$
    are mutually equivalent,
    so we have defined essentially three different stability conditions here.
\end{definition}

One can verify~\cites[Assumptions~5.1--5.3]{Joyce2021}
for the family of stability conditions $\acute{\mu}_{p,t}$
on the category $\Cohb (X)^L_{\smash{p}}$ as follows.
Assumptions~5.1 and~5.2~(b)--(h)
can be proved similarly to~\cite[\S\S8.2.2--8.2.3, \S8.3.1]{Joyce2021},
where the permissible classes are those $(\alpha, e) = ((r, d), e)$
with $r > 0$ and $e = 0, 1$;
Assumption~5.2~(a) follows from
Joyce--Song~\cite[Lemma~13.2]{JoyceSong2012};
Assumption~5.3~(a) is automatic,
and (b) follows from the beginning of~\cite[\S13.3]{JoyceSong2012}.

Therefore, the wall-crossing formulae
\cites[Theorem~5.9]{Joyce2021} apply
when one changes between the stability conditions
$\acute{\mu}_{p,t}$ for different~$t$.

\begin{theorem}
    \label{thm-wcf-2}
    Let $(r,d) \in K(X)$ with $r>0$. Then
    \upshape
    \begin{multline}
        \inv{\Mbss_{(r,d),1}} =
        \sum_{ \leftsubstack{
            \\[-1ex]
            & (r,d) = (r_1,d_1) + \cdots + (r_m,d_m), \\[-.6ex]
            & m \geq 1, \ r_i > 0, \ 
            d_i / r_i = d / r \text{ for all } i
        } }
        \frac {(-1)^m} {m!} \cdot
        \bigl[ \bigl[ \dotsc \bigl[ \bigl[
            \upe^{((0,0),1)}, \ 
            \inv{\Mss_{(r_1,d_1)}} \bigr], \\[-5.75ex]
            \shoveright{\inv{\Mss_{(r_2,d_2)}} \bigr] , \dotsc \bigr] , \ } \\
            \inv{\Mss_{(r_m,d_m)}}
        \bigr] .
    \end{multline}
\end{theorem}

\begin{proof}
    In Definition~\ref{def-three-stab-cond}, take $p = d/r$.
    We apply \cites[Theorem~5.9]{Joyce2021}
    to wall-cross from $t=-1$ to $t=1$.
    That is, we set $\tau = \acute{\mu}_{p,-1}$
    and $\tilde{\tau} = \acute{\mu}_{p,1}$ in the cited theorem.
    
    We have
    $\smash{\Mbss_{(r,d),1} (\acute{\mu}_{p,1})} =
    \smash{\Mbss_{(r,d),1} (\acute{\mu})}$.
    Also, note that
    $\smash{\Mbss_{(r',d'),1} (\acute{\mu}_{p,-1})} = \varnothing$
    for all $(r', d') \neq (0, 0)$,
    since for any $(E, V, \rho) \in \Cohb (X)^L_{(r',d'),1}$\,,
    the sub-object $(E, 0, 0)$ is always destabilizing.
    On the other hand, 
    $\smash{\Mbss_{(0,0),1} (\acute{\mu}_{p,-1})} = \{ (0, \bbC, 0) \}$
    consists of a single point, and
    $\smash{\Mbss_{(r',d'),0} (\acute{\mu}_{p,-1})} =
    \smash{\Mss_{(r',d')} (\mu)}$ for all $(r', d') \neq (0, 0)$.
    Setting $\alpha = (r, d)$,
    The wall-crossing formula reads
    \begin{multline}
        \inv{\Mbss_{\alpha,1}} =
        \sum_{ \leftsubstack[4em]{
            \\[-2ex]
            & (\alpha,1) =
            (0,1) + (\alpha_1,0) + \cdots + (\alpha_m,0), \\[-.6ex]
            & m \geq 1, \ \alpha_i \in C (X)_p \text{ for all } i, \\[-.6ex]
            & 0 \leq j \leq m
        } }
        \tilde{U} \bigl(
            (\alpha_1, 0), \dotsc,
            (\alpha_j, 0), (0, 1),
            (\alpha_{j+1}, 0), \dotsc,
            (\alpha_m, 0);
            \acute{\mu}_{p,-1}, \acute{\mu}_{p,1}
        \bigr) \cdot {} \\
        \bigl[ \bigl[ \dotsc \bigl[ \bigl[ \bigl[ \bigl[ \dotsc \bigl[
            \inv{\Mss_{\alpha_1}} \, , \,
            \inv{\Mss_{\alpha_2}} \bigr] , \, \dotsc \bigr], \,
            \inv{\Mss_{\alpha_j}} \bigr] , \,
            \inv{\Mbss_{0,1} (\acute{\mu}_{p,-1})} \bigr] , \\
            \inv{\Mss_{\alpha_{j+1}}} \bigr] , \, \dotsc \bigr], \,
            \inv{\Mss_{\alpha_m}}
        \bigr] .
    \end{multline}
    The coefficients $\tilde{U} ({\cdots})$
    were computed in
    \cite[Propositions~13.8 and~13.10]{JoyceSong2012},
    where it was shown that one could take
    \begin{equation}
        \tilde{U} ({\cdots}) = \begin{cases}
            \dfrac{(-1)^m}{m!}, & j = 0, \\
            0, & \text{otherwise}.
        \end{cases}
    \end{equation}
    The theorem then follows from the fact that
    $\smash{\inv{\Mbss_{0,1} (\acute{\mu}_{p,-1})}} =
    \upe^{((0, 0), 1)}$,
    since this moduli space is a single point.
\end{proof}

\section{Counting rank zero sheaves}
\label{sect-rank-0}

Let $X$ be a smooth projective curve over $\mathbb{C}$,
and let $d > 0$ be an integer.
The goal of this section is to compute the invariants
$\smash{\inv{\Mss_{(0,d)}}}$
counting $0$-dimensional sheaves on $X$,
and the main result will be presented as
Theorem~\ref{thm-inv-mss0d-main} below.
The author is unsure whether
these invariants are related to existing results in the literature.
These invariants will also serve as the base case
for computing invariants for elliptic curves
in~\S\ref{sect-elliptic-curves} below.

Let $X^{[d]}$ denote the Hilbert scheme of
$0$-dimensional length $d$ subschemes of $X$.
Each closed point of $X^{[d]}$ corresponds to an unordered $d$-tuple
$[x_1, \dotsc, x_d]$, where each $x_i$ is a point of $X$.

We have a symmetrization map
\begin{equation}
    \sym \colon X^d \to X^{[d]},
\end{equation}
which sends $(x_1, \dotsc, x_d) \in X^d$
to $[x_1, \dotsc, x_d] \in X^{[d]}$,
where $[x_1, \dotsc, x_d]$ denotes an unordered $d$-tuple.
The map $\sym$ is a $d!$-fold branched covering of $X^{[d]}$.

For $\bix = [x_1, \dotsc, x_d] \in X^{[d]}$,
define a sheaf $\calO_{\bix} \in \Coh(X)$ by the quotient
\begin{equation}
    0 \to
    \calO_X (-x_1 - \cdots - x_d) \longrightarrow
    \calO_X \longrightarrow
    \calO_{\bix} \to 0.
\end{equation}
In other words, $\calO_{\bix} = \calO_D$
where $D = x_1 + \cdots + x_d$ as a divisor of $X$.

We have a natural identification
\begin{equation}
    \label{eq-mbss0d1-equals-hilb}
    \Mbss_{(0,d),1} \simeq X^{[d]},
\end{equation}
where the left-hand side is the moduli stack of
semistable pairs (see \S\ref{sect-background-pairs}).
This is because for a rank zero sheaf $\calF$,
a pair $\rho \colon L \to \calF$ is semistable
if and only if $\rho$ does not factor through a proper subsheaf of $\calF$.
This means that $\calF \simeq \calO_{\bix}$ for some $\bix \in X^{[d]}$,
and that $\rho$ takes non-zero values at each point in $\bix$.
Note that all semistable pairs in $\calF$ are isomorphic,
and they do not have non-trivial automorphisms.
This establishes \eqref{eq-mbss0d1-equals-hilb}.

On the other hand, the stack $\smash{\Mss_{(0,d)}}$
is not as nice as $\smash{\Mbss_{(0,d),1}}$.
It is a genuine stack, since the sheaves $\calO_{\bix}$
have non-trivial automorphisms when $d > 1$
(where `non-trivial' means `other than a scalar multiplication').
Moreover, $\Mss_{(0,d)}$ also contains sheaves other than
those of the form $\calO_{\bix}$,
since all rank zero sheaves are semistable.
For example, for $x \in X$, the sheaf $\calO_x^{\oplus d}$
is in $\Mss_{(0,d)}$, but it is not involved in any semistable pair.

Let
\begin{equation}
    \label{eq-def-mbss0d1-projection}
    \Pi \colon X^{[d]} \simeq \Mbss_{(0,d),1} \to \Mss_{(0,d)}
\end{equation}
be the projection,
where the first isomorphism is by \eqref{eq-mbss0d1-equals-hilb}.

Let $\calE \in \Coh (X \times X^{[d]})$ be the universal sheaf
whose restriction to $X \times \{ \bix \}$ gives $\calO_{\bix}$.
In other words, we have $\calE = \calO_D$,
where $D \subset X \times X^{[d]}$
is the divisor consisting of all $(x_0, \bix)$
such that $x_0 \in \bix$.

Let $E = (\pr_2)_* (\calE)$ as a sheaf on $X^{[d]}$,
where $\pr_2 \colon X \times X^{[d]} \to X^{[d]}$ is the projection.
Then $E \to X^{[d]}$ is a vector bundle of rank $d$
whose fibre over $\bix \in X^{[d]}$ can be identified with the space of
sections of the sheaf $\calO_{\bix}$ on $X$.

\begin{lemma}
    \label{lem-rel-tangent-exact}
    We have an exact sequence
    \begin{equation}
        0 \to
        \calO_{X^{[d]}} \longrightarrow
        E \longrightarrow
        \bbT_{\Mbss_{(0,d),1} / \Mss_{(0,d)}}
        \to 0
    \end{equation}
    of vector bundles over $X^{[d]},$
    where the last non-zero term is
    the relative tangent bundle of\/
    $\smash{\Mbss_{(0,d),1}} \simeq X^{[d]}$ over $\Mss_{(0,d)}$.
\end{lemma}

\begin{proof}
    Consider the diagram
    \begin{equation}
        \begin{tikzcd}[
            row sep={3em,between origins},
            column sep={4.5em,between origins}
        ]
            & \bbP (E) \ar[rr] \ar[dd]
                \ar[rrdd, phantom, start anchor=-45, pos=.1, "\ulcorner"]
            && \Mbpl[, \neq 0]_{(0,d),1} \ar[dd] \\
            P \ar[ur, hook, "j"]
                \ar[rr, pos=.45, "q", crossing over] \ar[dd, "p"']
                \ar[rrdd, phantom, start anchor=-45, pos=.1, "\ulcorner"]
            && \Mbss_{(0,d),1} \ar[ur, hook] \\
            & X^{[d]} \ar[rr, pos=.3, "\Pi"]
            && \Mpl_{(0,d)} \\
            X^{[d]} \ar[ur, equals] \ar[rr, "\Pi"]
                \ar[uu, dashed, bend left, shift left, "i"]
            && \Mss_{(0,d)} \ar[ur, equals] \ar[from=uu, crossing over]
        \end{tikzcd}
    \end{equation}
    where $\bbP(E)$ is the projectivization of $E$,
    and $P \subset \bbP(E)$ is the open subset
    whose fibre over $\bix \in X^{[d]}$
    consists of those sections of $\calO_{\bix}$
    that do not vanish at any of the points in~$\bix$;
    $\smash{\Mpl_{(0,d)}}$ denotes the projective linear moduli stack
    of sheaves (not complexes) of class $(0,d)$,
    and $\smash{\Mbpl[, \neq 0]_{(0,d),1}}$
    denotes the projective linear moduli stack
    of pairs (again, not complexes) of class $((0,d),1)$,
    such that the structure map of the pair is non-zero.
    
    Consider the section $i \colon X^{[d]} \to P$
    given by the section $1$ of $\calO_{\bix}$ at each $\bix \in X^{[d]}$.
    Then the bundle $i^* \bbT_{P / X^{[d]}}$ over $X^{[d]}$
    fits into an exact sequence
    \begin{equation}
        \label{eq-exact-rel-tan}
        0 \to
        \calO_{X^{[d]}} \longrightarrow
        E \longrightarrow
        i^* \bbT_{P / X^{[d]}}
        \to 0,
    \end{equation}
    since for each $\bix \in X^{[d]}$,
    each element $s \in E_{\bix}$
    regarded as a section of $\calO_{\bix}$
    corresponds to an element of
    $(\bbT_{P / X^{[d]}})_{i(\bix)} \simeq
    (\bbT_{\bbP (E) / X^{[d]}})_{j \circ i(\bix)}$
    defined by the one-parameter family
    $t \mapsto 1 + t s$ of elements of $\bbP (E_{\bix})$.
    This gives a map $E \to i^* \bbT_{P / X^{[d]}}$
    whose kernel at $\bix$ consists of scalar multiples
    of the section $1$ of $\calO_{\bix}$.
    
    Note that if we identify $\Mbss_{(0,d),1}$ with $X^{[d]}$
    in the canonical way,
    then under this identification, we have $p = q$
    and $q \circ i = \id_{X^{[d]}}$.
    Since $\bbT_{P / X^{[d]}} \simeq
    q^* \smash{\bbT_{\Mbss_{(0,d),1} / \Mss_{(0,d)}}}$,
    the result follows from \eqref{eq-exact-rel-tan}.
\end{proof}

Define
\begin{align}
    \label{eq-e-prime}
    \calE' & = (\id_X \times \sym)^* (\calE), \\
    E' & = (\pr_2)_* (\calE'),
\end{align}
where $\calE'$ is a sheaf on $X \times X^d$,
and $E'$ a sheaf on $X^d$. Equivalently, we take
\begin{equation}
    \calE' = \calO_{\Delta_{0,1} + \cdots + \Delta_{0,d}} \, ,
\end{equation}
where $\Delta_{0,i} \subset X \times X^d$
is the divisor consisting of all $(x_0, x_1, \dotsc, x_d)$
with $x_0 = x_i$.

\begin{lemma}
    \label{lem-bundle-filtration}
    There is a filtration
    \begin{equation}
        0 = E_0 \subset
        E_1 \subset \cdots \subset
        E_d = E',
    \end{equation}
    where each $E_i$ is a rank~$i$ vector bundle on $X^d$, such that
    \begin{equation}
        E_i / E_{i-1} \simeq
        \calO_{X^d} (-\Delta_{1, i} - \cdots - \Delta_{i-1,i})
    \end{equation}
    for $i = 1, \dotsc, d$,
    where $\Delta_{j,i} \subset X^d$ is the divisor consisting of points
    $(x_1, \dotsc, x_d) \in X^d$ with $x_j = x_i$.
\end{lemma}

\begin{proof}
    For $i = 0, \dotsc, d$, 
    define sheaves $\calE_i \subset \calE'$ and $E_i \subset E$ by
    \begin{align}
        \calE_i & =
        \calO_{\Delta_{0,1} + \cdots + \Delta_{0,i}} \, , \\
        E_i & = (\pr_2)_* (\calE_i),
    \end{align}
    where $\Delta_{0,j} \subset X \times X^d$ are divisors as above,
    and $\pr_2 \colon X \times X^d \to X^d$ is the projection.
    Since $H^1 (X; \calO_{\bix}) = 0$ for all $\bix \in X^{[i]}$,
    we have $\bbR^1 (\pr_2)_* (\calE_i) = 0$ for all $i$,
    so that
    \begin{align*}
        E_i / E_{i-1}
        & \simeq (\pr_2)_* (\calE_i / \calE_{i-1})
        \\
        & \simeq (\pr_2)_* \bigl(
            \calO (-\Delta_{0,1} - \cdots - \Delta_{0,i-1}) /
            \calO (-\Delta_{0,1} - \cdots - \Delta_{0,i})
        \bigr) \\
        & \simeq (\pr_2)_* \bigl(
            \calO_{\Delta_{0,i}}
            (-\Delta_{0,1} - \cdots - \Delta_{0,i-1})
        \bigr) \\
        & \simeq \calO_{X^d} (-\Delta_{1,i} - \cdots - \Delta_{i-1,i}),
        \numberthis
    \end{align*}
    where the last step follows from the fact that
    $\pr_2 \colon \Delta_{0,i} \to X^d$ is an isomorphism.
\end{proof}

\begin{lemma}
    \label{lem-c-d-minus-1}
    One has
    \begin{equation}
        c_d (E') = 0, \quad
        c_{d-1} (E') = (-1)^{d-1} (d-1)! \PD ([\Delta]),
    \end{equation}
    where $\PD$ stands for `Poincar\'e dual',
    and $\Delta = \{ (x, \dotsc, x) \mid x \in X \} \subset X^d$
    is the diagonal.
\end{lemma}

\begin{proof}
    By Lemma~\ref{lem-bundle-filtration},
    the total Chern class of $E'$ is given by
    \begin{align*}
        c (E')
        & = \prod_{i=1}^d c \bigl(
            \calO_{X^d} (-\Delta_{1,i} - \cdots - \Delta_{i-1,i}) \bigr) \\
        & = \prod_{i=2}^d {} \bigl(
            1 - \PD([\Delta_{1,i}] + \cdots + [\Delta_{i-1,i}]) \bigr),
        \numberthis
    \end{align*}
    so that $c_d (E') = 0$ by dimensional reason, and
    \begin{align*}
        c_{d-1} (E')
        & = (-1)^{d-1} \prod_{i=2}^d {}
            \PD([\Delta_{1,i}] + \cdots + [\Delta_{i-1,i}]) \\
        & = (-1)^{d-1} \PD \biggl( {} \bigcap_{i=2}^d {}
            ([\Delta_{1,i}] + \cdots + [\Delta_{i-1,i}]) \biggr) \\
        & = (-1)^{d-1} (d-1)! \PD([\Delta]),
        \numberthis
    \end{align*}
    as each term in the expansion of the product
    gives a copy of $[\Delta]$.
\end{proof}

Recall that by \eqref{eq-h-mpl-is-quot-by-imd}, we have
$H_\bullet (\smash{\Mpl_{(0,d)}}; \bbQ) \simeq H_\bullet (\calM_{(0,d)}; \bbQ) / \im D_{(0,d)}$.
Let us consider the element $S_{1,0,1} \in H^2 (\calM_{(0,d)}; \bbQ)$.
Since it vanishes when paired with $\im D_{(0,d)}$, it defines a cohomology class
\begin{equation}
    \label{eq-def-s101-pl}
    S^\pl_{1,0,1} \in H^2 (\Mpl_{(0,d)}; \bbQ).
\end{equation}
Similarly, for $j, j' = 1, \dotsc, 2g$, we have well-defined cohomology classes
\begin{equation}
    \label{eq-def-sj11-pl}
    S^\pl_{j,1,1} \in H^1 (\Mpl_{(0,d)}; \bbQ), \quad
    S^\pl_{j,1,1} \, S^\pl_{j',1,1} \in H^2 (\Mpl_{(0,d)}; \bbQ).
\end{equation}

\begin{lemma}
    \label{lem-pullback-s-01}
    We have
    \begin{equation}
        \label{eq-sym-star-pi-star-s101}
        \sym^* \, \Pi^* (S^\pl_{1,0,1}) =
        \sum_{i=1}^d {} (\pr_i)^* \epsilon_{1,2}
    \end{equation}
    in $H^2 (X^d; \bbQ)$,
    where $\Pi$ is the map \eqref{eq-def-mbss0d1-projection},
    $S_{1,0,1}$ is defined by \eqref{eq-def-s101-pl},
    $\pr_i \colon X^d \to X$ is the $i$-th projection,
    and $\epsilon_{1,2} \in H^2 (X)$ is the canonical generator.
    Moreover, for $1 \leq j \leq 2g$, we have
    \begin{equation}
        \label{eq-sym-star-pi-star-sj11}
        \sym^* \, \Pi^* (S^\pl_{j,1,1}) =
        \pm \sum_{i=1}^d {} (\pr_i)^* \epsilon_{j \pm g,1} \, ,
    \end{equation}
    where the sign is \textnormal{`$+$'} if $j \leq g$,
    or \textnormal{`$-$'} otherwise.
\end{lemma}

\begin{proof}
    By definition, one has
    \begin{equation}
        S_{1,0,1} = \ch_1 (\calU_{(0,d)}) \setminus e_{1,0}
    \end{equation}
    in $H^2 (\calM_{(0,d)}; \bbQ)$, where
    $\calU_{(0,d)} \to X \times \calM_{(0,d)}$ is the universal complex,
    and $e_{1,0} \in H_0 (X; \bbQ)$ is the class of a point.
    
    Let $\Pi' \colon X^{[d]} \to \calM_{(0,d)}$
    be the map classifying the sheaf $\calE$. 
    Then $\Pi = \Pi^\pl \circ \Pi'$,
    where $\Pi^\pl \colon \calM_{(0,d)} \to \Mpl_{(0,d)}$
    is the projection.
    Note that we have $(\Pi^\pl)^* \, S^\pl_{1,0,1} = S_{1,0,1}$.
    
    Since
    $\calE' \simeq (\id_X \times (\sym \circ \Pi'))^* (\calU_{(0,d)})$,
    where $\calE'$ is given by \eqref{eq-e-prime},
    we have
    \begin{equation}
        \sym^* \, \Pi^* (S^\pl_{1,0,1}) =
        \sym^* \, \Pi'^* (S_{1,0,1}) =
        \ch_1 (\calE') \setminus [\pt]_X.
    \end{equation}
    By the definition of $\calE'$, we have an exact sequence
    \begin{equation}
        0 \to \calO_{X \times X^d} (-\Delta_{0,1} - \cdots - \Delta_{0,d})
        \longrightarrow \calO_{X \times X^d}
        \longrightarrow \calE' \to 0,
    \end{equation}
    which gives
    \begin{align*}
        \ch_1 (\calE')
        & = -c_1 \bigl(
            \calO_{X \times X^d} (-\Delta_{0,1} - \cdots - \Delta_{0,d}) \bigr) \\
        & = \PD([\Delta_{0,1}] + \cdots + [\Delta_{0,d}]),
        \numberthis
    \end{align*}
    and since $\PD([\Delta_{0,i}]) \setminus [\pt]_X = (\pr_i)^* \epsilon_{1,2}$,
    we have proved \eqref{eq-sym-star-pi-star-s101}.
    
    Similarly, for $1 \leq j \leq 2g$, we have
    \begin{equation}
        \sym^* \, \Pi^* (S^\pl_{j,1,1})
        = \ch_1 (\calE') \setminus e_{j,1} \, ,
    \end{equation}
    where $S^\pl_{j,1,1}$ is given by \eqref{eq-def-sj11-pl}.
    But $\PD([\Delta_{0,i}]) \setminus e_{j,1} =
    \pm (\pr_i)^* \epsilon_{j \pm g, 1}$,
    where the sign is `$+$' if $j \leq g$, or `$-$' otherwise.
    This proves \eqref{eq-sym-star-pi-star-sj11}.
\end{proof}

\begin{lemma}
    \label{lem-lie-mss0d-vanish}
    The Lie subalgebra
    \[
        \bigoplus_{d>0}
        \check{H}_0 (\Mpl_{(0,d)}; \bbQ)
        \ \subset \ 
        \check{H}_0 (\Mpl_{>0}; \bbQ)
    \]
    is abelian.
\end{lemma}

\begin{proof}
    Note that $\smash{\check{H}_0 (\Mpl_{(0,d)}; \bbQ)}$
    is the quotient of $\smash{\hat{H}_2 (\calM_{(0,d)}; \bbQ)}$
    by $\im D_{(0,d)}$.
    Since 
    \begin{equation}
        \hat{H}_2 (\calM_{(0,d)}; \bbQ) \simeq
        \upe^{(0,d)} \cdot \langle
            s_{1,0,1} \, , \ 
            s_{1,2,2} \, , \ 
            s_{j,1,1} \, s_{j',1,1} :
            1 \leq j < j' \leq 2g
        \rangle_{\bbQ} \, ,
    \end{equation}
    for $\upe^{(0,d)} \cdot p \in \smash{\hat{H}_2 (\calM_{(0,d)}; \bbQ)}$
    and $\upe^{(0,d')} \cdot p' \in \smash{\hat{H}_2 (\calM_{(0,d')}; \bbQ)}$,
    the vertex operation can be simplified to
    \begin{multline}
        \label{eq-y-explicit-rank-zero}
        Y (\upe^{(0,d)} \cdot p, \ z)
        (\upe^{(0,d')} \cdot p') =
        \upe^{(0,d+d')} \cdot {} \\
        \exp (z D_{(0, d)})
        \exp \Bigl(
            (2-2g) \, z^{-2} \,
            \frac{\partial}{\partial s_{1,0,1}}
            \frac{\partial}{\partial s'_{1,0,1}}
        \Bigr)
        \bigl( p (s_{jkl}) \cdot p' (s'_{jkl}) \bigr) \,
        \big|_{s_{jkl} = s'_{jkl}} \ .
    \end{multline}
    Let $a$ and $a'$ be the coefficient of $s_{1,0,1}$
    in $p$ and $p'$, respectively.
    Then the $z = 0$ residue of \eqref{eq-y-explicit-rank-zero}
    equals
    \begin{equation}
        \upe^{(0,d+d')} \cdot
        (2-2g) a a' \,
        D_{(0,d)} 1 =
        \upe^{(0,d+d')} \cdot
        (2-2g) a a' \, \frac{d}{d+d'} \, D_{(0,d+d')} 1,
    \end{equation}
    which lies in $\im D_{(0,d+d')}$,
    and hence,
    is zero in $\smash{\check{H}_0 (\Mpl_{(0,d+d')}; \bbQ)}$.
\end{proof}

\begin{theorem}
    \label{thm-inv-mss0d-main}
    For any integer $d > 0$, we have
    \begin{equation}
        \inv{ \Mss_{(0,d)} } =
        \upe^{(0,d)} \cdot
        (-1)^{d-1}
        \biggl(
            \frac{1}{d} \, s_{1,0,1} +
            \sum_{j=1}^g s_{j,1,1} \, s_{j+g,1,1}
        \biggr)
        + \im D_{(0,d)}.
    \end{equation}
\end{theorem}

\begin{proof}
    Since
    \begin{align}
        \label{eq-h2m0d}
        H_2 (\calM_{(0,d)}; \bbQ) & = \langle
            s_{1,0,1} \, , s_{1,2,2} \, ,
            s_{j,1,1} \, s_{j',1,1} \mid
            1 \leq j < j' \leq 2g
        \rangle_{\bbQ} \ , \\
        \label{eq-dh0m0d}
        D (H_0 (\calM_{(0,d)}; \bbQ)) & = \langle
            s_{1,2,2}
        \rangle_{\bbQ} \ ,
    \end{align}
    and $H_2 (\Mpl_{(0,d)}; \bbQ)$
    is the quotient of~\eqref{eq-h2m0d} by~\eqref{eq-dh0m0d},
    we must have
    \begin{equation}
        \inv{ \Mss_{(0,d)} } =
        a \, s_{1,0,1} + \sum_{j<j'} b_{j,j'} \,
        s_{j,1,1} \, s_{j',1,1} + \im D_{(0, d)}
    \end{equation}
    for some coefficients $a, b_{j,j'} \in \bbQ$.
    
    By \cite[Theorem~7.63]{Joyce2021}
    and Lemma~\ref{lem-rel-tangent-exact},
    we have
    \begin{equation}
        \label{eq-inv-mss0d-equals-cap-c-top}
        \inv{ \Mss_{(0,d)} } =
        \frac{1}{d} \, \Pi_* ([X^{[d]}] \cap c_{d-1} (E)),
    \end{equation}
    where all the Lie brackets in the cited theorem vanish
    by Lemma~\ref{lem-lie-mss0d-vanish}.
    
    Pairing both sides of \eqref{eq-inv-mss0d-equals-cap-c-top}
    with the cohomology class $S^\pl_{1,0,1}$, we have
    \begin{align*}
        a
        & = \frac{1}{d} \bigl<
            S^\pl_{1,0,1}, \ 
            \Pi_* ([X^{[d]}] \cap c_{d-1} (E))
        \bigr> \\
        & = \frac{1}{d} \bigl<
            \Pi^* S^\pl_{1,0,1}, \ 
            [X^{[d]}] \cap c_{d-1} (E)
        \bigr> \\
        & = \frac{1}{d} \bigl<
            \Pi^* S^\pl_{1,0,1} \cup c_{d-1} (E), \ 
            [X^{[d]}]
        \bigr> \\
        & = \frac{1}{d \cdot d!} \bigl<
            \sym^* \, \Pi^* S^\pl_{1,0,1} \cup c_{d-1} (E'), \ 
            [X^d]
        \bigr> \\
        & = \frac{(-1)^{d-1}}{d^2}
        \int_{\Delta \subset X^d} 
        \sum_{i=1}^d {} (\pr_i)^* \epsilon_{1,2}
        \\
        & = \frac{(-1)^{d-1}}{d}, 
        \numberthis
    \end{align*}
    where we used the fact that $\sym_* ([X^d]) = d! \, [X^{[d]}]$,
    and we used Lemma~\ref{lem-c-d-minus-1}
    and Lemma~\ref{lem-pullback-s-01}.
    Similarly, we obtain that for $1 \leq j \leq g$,
    \begin{align*}
        b_{j,j+g} & =
        \frac{(-1)^{d-1}}{d^2}
        \int_{\Delta \subset X^d}
        \sum_{i=1}^d \sum_{i'=1}^d {}
        (\pr_i)^* \epsilon_{j,1} \cup
        (\pr_{i'})^* \epsilon_{j+g,1} \\
        & = (-1)^{d-1}, 
        \numberthis
    \end{align*}
    and that all the other $b_{j,j'}$ are zero.
\end{proof}

\section{Counting line bundles}
\label{sect-rank-1}

Let $X$ be a smooth, projective curve over $\bbC$.
In this section, we compute the rank~$1$ invariants
$\inv{\Mss_{(1,d)}}$ and $\inv{\Mbss_{(1,d),1}}$ for all $d$.
The former is just the fundamental class of
a component of the Picard variety of $X$,
but we express it in terms of the variables $s_{j,k,l}$ explicitly
in Theorem~\ref{thm-fund-mss-1-d} below.
The pair invariants coincide with virtual fundamental classes
of the pair moduli stacks, and are computed explicitly
in Theorem~\ref{thm-fund-mbss-1-d-1} below.

Let $J (X) = \Mss_{(1,0)}$ be the Jacobian variety of $X$,
and let $\calU_{(1,0)} \to X \times J (X)$
be the universal line bundle
such that $\calU_{(1,0)} |_{\{ x_0 \} \times J(X)}$
is the trivial line bundle over $J(X)$ for some chosen $x_0 \in X$.

A classical result of Narasimhan--Seshadri~\cite[Theorem~2]{narasimhan-seshadri}
states that stable holomorphic bundles on $X$ of degree zero
arise from irreducible unitary representations
of the fundamental group of $X$,
and isomorphic bundles correspond to equivalent representations.
(See also \cite[\S8]{AtiyahBott1983}.)
This provides a homeomorphism
\begin{equation}
    J (X) \simeq
    \operatorname{Rep} (\pi_1 (X), \mathrm{U} (1)),
\end{equation}
where $\operatorname{Rep} (\pi_1 (X), \mathrm{U} (1))$
is the space of $1$-dimensional unitary representations
of $\pi_1 (X)$. 
Since $H_1 (X; \bbZ)$ is the abelianization of $\pi_1 (X)$,
we have an identification
\begin{equation}
    \label{eq-jacobian-equals-torus}
    J (X) \simeq
    \operatorname{Rep} (H_1 (X; \bbZ), \mathrm{U} (1)) \simeq
    \mathrm{U} (1)^{2g}.
\end{equation}
For $j = 1, \dotsc, 2g$, let
\begin{equation}
    \label{eq-def-gamma-j}
    \gamma_j = \pr_j^* (\omega)
    \quad \in H^1 (J (X); \bbZ),
\end{equation}
where $\pr_j \colon J (X) \to \mathrm{U} (1)$
is induced by the element $e_{j,1} \in H_1 (X; \bbZ)$
via the isomorphism \eqref{eq-jacobian-equals-torus},
and $\omega \in H^1 (\mathrm{U} (1); \bbZ)$
is the Poincar\'e dual of a point.

Now, for any $d \in \bbZ$,
if we fix a line bundle $L_d \to X$ of degree $d$,
then taking the tensor product with $L_d$ gives an isomorphism
\begin{equation}
    \label{eq-jacobian-equals-mss1d}
    {-} \otimes L_d \colon
    J (X) \simeq \Mss_{(1,0)}
    \longrightarrow \Mss_{(1,d)}.
\end{equation}
Let $\calU_{(1,d)} \to X \times \Mss_{(1,d)}$
denote the universal line bundle of degree $d$,
obtained as
\begin{equation}
    \label{eq-constr-of-u1d}
    \calU_{(1,d)} = (\id_X \times \varphi)^* \, \calU_{(1,0)} 
    \otimes \pr_1^* (L_d),
\end{equation}
where $\varphi$ is the inverse of the isomorphism
\eqref{eq-jacobian-equals-mss1d},
and $\pr_1 \colon X \times \Mss_{(1,d)} \to X$
is the projection.

The Chern class of~$\calU_{(1,d)}$ can be computed
by a standard argument to be
\begin{equation}
    \label{eq-c1-u1d}
    c_1 (\calU_{(1,d)}) =
    -\sum_{j=1}^{2g} \epsilon_{j,1} \boxtimes \gamma_{j} +
    d \epsilon_{1,2} \boxtimes 1,
\end{equation}
where $\gamma_j \in H^1 (\Mss_{(1,d)}; \bbZ)$
is defined by identifying $\Mss_{(1,d)}$ with $J(X)$
via \eqref{eq-jacobian-equals-mss1d}.
Using this,
we can relate the classes $\gamma_j$
with the classes $S_{j,k,l}$ as follows.

\begin{lemma}
    \label{lem-pb-of-sjkl-to-jacobian}
    Let $i \colon J (X) \hookrightarrow \Mpl_{(1,d)}$ denote the inclusion,
    where we have identified $J (X) \simeq \Mss_{(1,d)}$ via
    \eqref{eq-jacobian-equals-mss1d}. 
    Let $\Xi$ be the map defined in \textnormal{\S\ref{sect-def-xi}}.
    Then
    \begin{align}
        i^* \Xi (S_{j,1,1}) & =
        -\gamma_j \, , \\
        i^* \Xi (S_{1,2,2}) & =
        -\sum_{j=1}^g \gamma_j \, \gamma_{j+g} \, ,
    \end{align}
    and $i^* \Xi (S_{j,k,l}) = 0$ for all other $(j, k, l)$ with
    $(j, k) \in J$ and $l > k/2$.
\end{lemma}

\begin{proof}
    By~\eqref{eq-c1-u1d}, we have
    \begin{align}
        \ch_1 (\calU_{(1,d)}) & =
        c_1 (\calU_{(1,d)}) =
        -\sum_{j=1}^{2g} \epsilon_{j,1} \boxtimes \gamma_{j} +
        d \epsilon_{1,2} \boxtimes 1, \\
        \ch_2 (\calU_{(1,d)}) & =
        \frac{1}{2} c_1 (\calU_{(1,d)})^2 =
        -\epsilon_{1,2} \boxtimes \sum_{j=1}^g \gamma_j \, \gamma_{j+g} \, , \\
        \ch_l (\calU_{(1,d)}) & =
        0 \quad (l \geq 3).
    \end{align}
    The result then follows from the fact that
    \begin{equation}
        i^* \Xi (S_{j,k,l}) =
        \ch_l (\calU_{(1,d)}) \setminus e_{j,k}
    \end{equation}
    for all $j, k, l$.
\end{proof}

To simplify notations, write
\begin{align}
    \label{eq-def-rho}
    \rho (w) & = \exp \biggl(
        \sum_{l=1}^\infty \frac{(-1)^l}{l!} \, w^l \, s_{{+},0,l}
    \biggr), \\
    \label{eq-def-sigma}
    \sigma (x) & =
    \prod_{j=1}^g {} (x + s_{j,1,1} \, s_{j+g,1,1}),
\end{align}
where $w$ is a formal variable of degree $-2$,
and $x$ is a formal variable of degree $2$.

\begin{theorem}
    \label{thm-fund-mss-1-d}
    For any integer $d$,
    the fundamental class $\fund{\Mss_{(1, d)}}$ is given by
    \[
        \xi (\fund{\Mss_{(1, d)}}) =
        \upe^{(1, d)} \cdot \sigma (-s_{1,2,2}),
    \]
    where $\xi$ is the map defined in \textnormal{\S\ref{sect-def-xi}},
    and $\sigma$ is as in \eqref{eq-def-sigma}.
\end{theorem}

\begin{proof}
    By Lemma~\ref{lem-pb-of-sjkl-to-jacobian},
    the pairing of any monomial in the variables $S_{j,k,l}$
    with the homology class $\xi (\fund{\Mss_{(1,d)}})$
    will be zero, unless that monomial only involves
    the variables $S_{j,1,1}$ and $S_{1,2,2}$\,,
    and that for $j = 1, \dotsc, g$,
    the variable $S_{j,1,1}$ is present
    if and only if $S_{j+g,1,1}$ is also present.
    This means that $\xi (\fund{\Mss_{(1,d)}})$
    is a polynomial in the variables
    $\sigma_j = s_{j,1,1} \, s_{j+g,1,1}$ and $s_{1,2,2}$\,,
    where $j = 1, \dotsc, g$.
    
    Write $\delta_j = \gamma_j \, \gamma_{j+g}$ for $j = 1, \dotsc, g$.
    Then the coefficient of
    $s_{1,2,2}^m \, \sigma_1^{\smash{m_1}} \cdots \sigma_g^{\smash{m_g}}$
    in $\xi (\fund{\Mss_{(1, d)}})$
    is $1/m!$ times the pairing of
    $S_{1,2,2}^m \, (S_{1,1,1} \, S_{1+g,1,1})^{m_1} \cdots
    (S_{g,1,1} \, S_{2g,1,1})^{m_g}$
    with $\xi (\fund{\Mss_{(1, d)}})$,
    which is $1/m!$ times the
    coefficient of $\delta_1 \cdots \delta_g$ in
    $(-\delta_1 - \cdots - \delta_g)^m \,
    \delta_1^{\smash{m_1}} \cdots \delta_g^{\smash{m_g}}$,
    which is only non-zero when $m_j \in \{0, 1\}$ for all $j$
    and when $m + m_1 + \cdots + m_g = g$. In this case,
    the coefficient in question is
    $(-1)^m \, m!$, so that
    \begin{align*}
        \xi (\fund{\Mss_{(1, d)}}) & =
        \upe^{(1, d)} \cdot \sum_{ \substack{
            m + m_1 + \cdots + m_g = g \\
            m \geq 0, \ m_j \in \{0, 1\}
        } }
        (-1)^m \, s_{1,2,2}^m \,
        \sigma_1^{\smash{m_1}} \cdots \sigma_g^{\smash{m_g}} \\
        & = \upe^{(1, d)} \cdot \sigma (-s_{1,2,2}).
        \numberthis
    \end{align*}
\end{proof}

\begin{theorem}
    \label{thm-fund-mbss-1-d-1}
    For any integer $d$,
    the virtual fundamental class
    $\virt{\Mbss_{(1, d), 1}}$ is given by
    \begin{equation}
        \acute{\xi} (\virt{\Mbss_{(1, d), 1}}) =
        \upe^{((1, d), 1)} \cdot \res_w \biggl(
            \frac{1}{w^{\nu+d}} \,
            \rho (w) \,
            \sigma \Bigl( \frac{1}{w} - s_{1,2,2} \Bigr)
        \biggr),
    \end{equation}
    where $\acute{\xi}$ is the map defined in \textnormal{\S\ref{sect-def-xi}},
    $\rho$ and $\sigma$ are given by
    \eqref{eq-def-rho} and \eqref{eq-def-sigma},
    $\nu$ is as in Notation~\textnormal{\ref{ntn-n}},
    and $w$ is a formal variable of degree $-2$.
\end{theorem}

\begin{proof}
    By Theorem~\ref{thm-wcf-2}, we have
    \begin{equation}
        \virt{\Mbss_{(1, d), 1}} =
        -\bigl[ \upe^{((0,0),1)}, \ \fund{\Mss_{(1,d)}} \bigr].
    \end{equation}
    Using the expression \eqref{eq-vertex-algebra-pairs}
    for the vertex operation for pairs,
    and using Theorem~\ref{thm-fund-mss-1-d},
    we see that
    \begin{align*}
        & \acute{Y} (\upe^{((0,0),1)}, \ -w) (\fund{\Mss_{(1,d)}}) \\
        = {} &
        \upe^{((1,d),1)} \cdot \frac{1}{w^{\nu+d}}
        \exp \biggl[ -w \biggl(
            s_{{+},0,1} +
            \sum_{l=1}^{\infty} s_{{+},0,l+1} \,
            \frac{\partial}{\partial s_{{+},0,l}}
        \biggr) \biggr] 
        \exp \biggl( -\frac{1}{w} \, \frac{\partial}{\partial s_{1,2,2}} \biggr) \ 
        \sigma (-s_{1,2,2}) \\
        = {} &
        \upe^{((1,d),1)} \cdot \frac{1}{w^{\nu+d}} \,
        \rho (w) \, \sigma \Bigl( \frac{1}{w} - s_{1,2,2} \Bigr),
        \numberthis
    \end{align*}
    where we used the fact that
    \begin{equation}
        \exp \biggl[ -w \biggl(
            s_{{+},0,1} +
            \sum_{l=1}^{\infty} s_{{+},0,l+1} \,
            \frac{\partial}{\partial s_{{+},0,l}}
        \biggr) \biggr] \ f =
        \rho (w) \, f
    \end{equation}
    for any polynomial $f$ not involving the variables $s_{{+},0,l}$,
    and we used the identity
    \begin{equation}
        \exp \Bigl( a \, \frac{\partial}{\partial x} \Bigr)
        \, \sigma (x) = \sigma (a + x)
    \end{equation}
    for formal variables $a, x$ of degree~$2$.
\end{proof}

\section{Counting higher rank bundles}
\label{sect-higher-rank}

In this section, we state the main results of this paper,
Theorems~\ref{thm-inv-mss-main-as-reg-sum} and
\ref{thm-inv-mss-main},
which compute the invariants $\inv{\Mss_{(r,d)}}$ for $r \geq 2$.
The explicit expressions will involve the \emph{regularized sum},
which is a way to assign finite values to divergent series,
and will be rigorously defined in \S\ref{sect-regularized-sum} below.
Theorem~\ref{thm-inv-mss-main-as-reg-sum}
expresses the invariants $\inv{\Mss_{(r,d)}}$
as a regularized sum
of iterated Lie brackets of the rank~$1$ invariants $\inv{\Mss_{(1,d)}}$
computed in \S\ref{sect-rank-1} above,
where the Lie brackets come from Joyce's vertex algebra structure
discussed in \S\ref{sect-background-sheaves}.
Theorem~\ref{thm-inv-mss-main} gives an explicit formula
for these invariants, in terms of the variables $s_{j,k,l}$
defined in \S\ref{sect-background-sheaves}.

We also consider the fixed determinant version of these invariants
in \S\ref{sect-fixed-determinant},
and study the special case when $r = 2$ in \S\ref{sect-rank-2}.
In \S\ref{sect-comparison},
we compare our formula for intersection pairings
to results of Witten~\cite{Witten1992}
and Jeffrey--Kirwan~\cite{jeffrey-kirwan-2}.

\subsection{The main results}
\label{sect-main-results}

Before stating the main results,
we give a few preparatory definitions.

\begin{definition}[$F_r$, $\acute{F}_r$, $F_r^+$ and $\acute{F}_r^+$]
    \label{def-field-inf-many-variables}
    Let $r > 0$ be an integer.
    We define a field $F_r$ of meromorphic functions in
    infinitely many complex variables
    \[
        z_0, \dotsc, z_{r-1}; \quad
        s_{j,k,l} \ \ ((j, k) \in J, \ k \text{ even}, \ l > k/2),
    \]
    as follows.
    For each integer $L \geq 0$,
    let $R_{r,L}$ be the ring of holomorphic functions
    in finitely many variables $z_0, \dotsc, z_{r-1}$ and
    $s_{j,k,l}$ for $l \leq L$.
    When $L_1 > L_2$, we have a projection $R_{r,L_1} \to R_{r,L_2}$
    defined by setting $s_{j,k,l} = 0$ for $L_2 < l \leq L_1$. Define
    \begin{equation}
        R_r = \varprojlim_{L} R_{r,L}\,.
    \end{equation}
    Then $R_r$ is an integral domain.
    We define 
    \begin{equation}
        F_r = \text{fractional field of } R_r \, .
    \end{equation}
    
    Similarly, we define a field $\acute{F}_r$ of meromorphic functions in
    infinitely many complex variables
    \[
        w, z_0, \dotsc, z_{r-1}; \quad
        s_{j,k,l} \ \ ((j, k) \in \acute{J}, \ k \text{ even}, \ l > k/2),
    \]
    by the same procedure as above.
    
    Define an $F_r$-vector space $F_r^+$
    by adding the odd variables $s_{j,1,l}$ to $F_r$.
    Precisely, we define
    \begin{equation}
        F_r^+ = \Bigl(
            \varprojlim_L {} ( R_{r,L} \otimes_{\bbC} \Omega_L )
        \Bigr) \otimes_{R_r} F_r,
    \end{equation}
    where
    \begin{equation}
        \Omega_L =
        {\wedge}^\bullet \, \langle s_{j,1,l} : (j, 1) \in J, \ 0 < l \leq L \rangle_{\bbC}
    \end{equation}
    is the exterior algebra of the $\bbC$-vector space
    spanned by the elements $s_{j,1,l}$.
    We will often work with regularized sums in $F_r^+$
    over the field $F_r$.
    
    Define an $\acute{F}_r$-vector space $\acute{F}_r^+$
    by a similar process, adding in the variables $s_{j,1,l}$.
\end{definition}

In \S\ref{sect-regularized-sum},
we will define the \emph{regularized sum}.
This will involve a finite-dimensional $\bbQ$-vector space $V$,
a subset $\Lambda \subset V$ of the form
\begin{equation}
    \Lambda = \Bigl\{
        \bix_0 + \sum_{i=1}^m a_i \, \bie_i \Bigm|
        a_i \in \bbZ
    \Bigr\},
\end{equation}
with $\bix_0, \bie_1, \dotsc, \bie_m \in V$
and $\bie_1, \dotsc, \bie_m$ linearly independent.
We call $\Lambda$ a \emph{lattice}.
For each $\bix \in \Lambda$, we are given a meromorphic function
$f (\bix; \biz) \in F_r$.
Under certain assumptions, we can define the regularized sum
\begin{equation}
    \sumbar_{\bix \in \Lambda} f (\bix; \biz) \in F_r ,
\end{equation} 
as in Definition~\ref{def-regularized-sum},
and this process will be a key ingredient in our result.

For a subset $\Sigma \subset \Lambda$,
we may also write 
\begin{equation}
    \sumbar_{\bix \in \Sigma} f (\bix; \biz),
\end{equation}
which means that we are taking a regularized sum over $\Lambda$,
but we assign $f (\bix; \biz) = 0$ to all $\bix \notin \Sigma$.
Note that this will not be defined unless
the conditions in Definition~\ref{def-regularized-sum} are satisfied;
this roughly means that $\Sigma$ should be a \emph{sector}
in the sense of Definition~\ref{def-cone-lattice}.

\begin{definition} [Regularized sum of iterated Lie brackets]
    \label{def-reg-sum-of-residues}
    Let $V$ be a $\bbQ$-vector space of finite dimension,
    and let $\Lambda \subset V$ be a lattice as above.
    For a function $f \colon \Lambda \to F_r$,
    we will use the notation
    \begin{equation}
        \label{eq-reg-sum-res}
        \sumbar_{\bix \in \Lambda} \res_{z_{i_1}} \circ \cdots \circ \res_{z_{i_k}} (f (\bix; \biz)) :=
        \res_{z_{i_1}} \circ \cdots \circ \res_{z_{i_k}} \biggl( \sumbar_{\bix \in \Lambda} f (\bix; \biz) \biggr),
    \end{equation}
    whenever the right-hand side is defined.
    This is an abuse of notation,
    as the regularized sum on the left-hand side is often not defined,
    and when it is defined,
    it can be different from the right-hand side.
    But we always use the right-hand side of~\eqref{eq-reg-sum-res}
    as the definition of this expression.

    Recall that if $W$ is a vertex algebra,
    then the quotient space $W / \im D$ carries a Lie algebra structure,
    where $D$ is the translation operator, using the Lie bracket
    \begin{equation}
        [A + \im D, B + \im D] = \res_z Y (A, z) \, B + \im D
        = -\res_z (B, z) \, A + \im D 
    \end{equation}
    in $W / \im D$, where $A, B \in W$.
    
    We use~\eqref{eq-reg-sum-res} to define regularized sums of
    iterated Lie brackets coming from a vertex algebra.
    Given elements $A_{0, \bix}, \dotsc, A_{n, \bix} \in W$
    for each $\bix \in \Lambda$, by abuse of notation, we define
    \begin{multline}
        \label{eq-reg-sum-va}
        \sumbar_{\bix \in \Lambda} \ 
        [ [ \dotsc [ A_{0, \bix} + \im D, \, A_{1, \bix} + \im D ], \dotsc ],
        A_{n, \bix} + \im D ]
        := \\
        (-1)^n \cdot
        \sumbar_{\bix \in \Lambda} {}
        \res_{z_n} \circ \cdots \circ \res_{z_1}
        \bigl(
        Y (A_{n, \bix}, z_n) \cdots Y (A_{1, \bix}, z_1) \, A_{0, \bix} \bigr) + \im D
    \end{multline}
    as an element of $W / \im D$,
    where the right-hand side (excluding `$+ \, \im D$') is interpreted
    in the sense of~\eqref{eq-reg-sum-res}.
    Of course, this is not defined unless the conditions
    in Definition~\ref{def-regularized-sum} are satisfied.

    Note that despite the notation,
    the sum does depend on the choices of $A_{i, \bix}$
    as a representative in the coset $A_{i, \bix} + \im D$,
    since for instance,
    a bad choice would lead to the sum being undefined.
    Therefore, whenever one writes a regularized sum of this form,
    a choice of representatives must be made,
    although this choice is sometimes obvious from the context.

    There is also a graded, or super, version of the above,
    where a suitable sign is inserted when odd elements are involved.
    This will be the version used below,
    but we will not need the signs, 
    as the elements $A_{i, \bix}$ will be even
    in all cases that we will consider below.
\end{definition}

\begin{remark}
    \label{rmk-reg-sum-lie-brack}
    In Definition~\ref{def-reg-sum-of-residues},
    there is an ambiguity that $[A + \im D, B + \im D]$
    can mean either $\res_z Y (A, z) \, B$ or $-{\res_z Y (B, z) \, A}$,
    so an expression $\sumbar_{\bix} [A_{\bix} + \im D, B_{\bix} + \im D]$
    can have two possible interpretations;
    we chose the latter one in~\eqref{eq-reg-sum-va}.
    However, we argue that we can often choose either one of them,
    and obtain the same results.
    
    By Frenkel--Ben-Zvi~\cite[Proposition~3.2.5]{FrenkelBenZvi2004},
    for any vertex algebra $W$ and any $A, B \in W$, we have 
    \begin{equation}
        \label{eq-va-comm}
        Y (A, z) \, B = \upe^{z D} \, Y (B, -z) \, A.
    \end{equation}
    But in our case,
    $\upe^{z D}$ is a field automorphism of $F_r$ up to a constant in $F_r$,
    and commutes with the regularized sum by
    Lemma~\ref{lem-reg-sum-field-aut}.
    Thus, we are free to choose between the two interpretations.

    This can be generalized to iterated Lie brackets as well,
    using multi-variable versions of~\eqref{eq-va-comm}
    which can be obtained using the commutativity~\eqref{eq-va-comm}
    and associativity of vertex algebras.
    This means that we can freely apply anti-commutativity 
    and the Jacobi identity of the Lie bracket
    inside the regularized sum.
\end{remark}

Here is the first version of our main result,
computing the invariants $\inv{\Mss_{(r,d)}}$
and the pair invariants $\inv{\Mbss_{(r,d),1}}$ when $r > 0$.
The proof will be given in \S\ref{sect-proof-of-main-theorem},
but we will give a sketch proof here.

\begin{theorem}
    \label{thm-inv-mss-main-as-reg-sum}
    Let $(r, d) \in K(X)$ with $r > 0$.
    \begin{enumerate}
        \item \label{item-inv-mss-main-as-reg-sum}
        We have
        \upshape
        \begin{align*}
            \numberthis
            \label{eq-inv-mss-main-as-reg-sum}
            & \inv{\Mss_{(r,d)}} = 
            \frac{1}{r} \cdot {} \\[1ex]
            & \sumbar_{ \leftsubstack{
                \\[-2ex]
                & d = d_0 + \cdots + d_{r-1} \\[-.6ex]
                & (d_0 + \cdots + d_{i-1})/i \leq d/r, \ i = 1, \dotsc, r-1 \\[-.6ex]
                & \text{with $m$ equalities} 
            } } \hspace{-1.5em}
            \frac{1}{m+1} \cdot
            \bigl[ \bigl[ \dotsc \bigl[
                \fund{\Mss_{(1,d_0)}} \, , \ 
                \fund{\Mss_{(1,d_1)}} \bigr] , \ 
                \dotsc \bigr] , \ 
                \fund{\Mss_{(1,d_{r-1})}}
            \bigr].
        \end{align*}
        \itshape
        In particular, if $r$ and $d$ are coprime, then
        \upshape
        \begin{align*}
            \label{eq-inv-mss-coprime-as-reg-sum}
            \numberthis
            & \fund{\Mss_{(r,d)}} = 
            \frac{1}{r} \cdot {} \\[1ex]
            & \hspace{1em}
            \sumbar_{ \leftsubstack{ 
                \\[-2ex]
                & d = d_0 + \cdots + d_{r-1} \\[-.6ex]
                & (d_0 + \cdots + d_{i-1})/i < d/r, \ i = 1, \dotsc, r-1
            } } \mspace{-9mu}
            \bigl[ \bigl[ \dotsc \bigl[
                \fund{\Mss_{(1,d_0)}} \, , \ 
                \fund{\Mss_{(1,d_1)}} \bigr] , \ 
                \dotsc \bigr] , \ 
                \fund{\Mss_{(1,d_{r-1})}}
            \bigr].
        \end{align*}
        \itshape
        
        \item \label{item-inv-mbss-main-as-reg-sum}
        We have
        \begin{align*}
            \numberthis 
            \label{eq-inv-mbss-main-as-reg-sum}
            & \inv{\Mbss_{(r,d),1}} =
            \\[1ex] & \hspace{1em}
            \sumbar_{ \leftsubstack{ 
                \\[-3ex]
                & d = d_0 + \cdots + d_{r-1} \\[-.6ex]
                & (d_0 + \cdots + d_{i-1})/i < d/r, \ i = 1, \dotsc, r-1
            } } \hspace{-.5em}
            \bigl[ \bigl[ \dotsc \bigl[
                \inv{\Mbss_{(1,d_0),1}} \, , \ 
                \fund{\Mss_{(1,d_1)}} \bigr] , \ 
                \dotsc \bigr] , \ 
                \fund{\Mss_{(1,d_{r-1})}}
            \bigr].
        \end{align*}
    \end{enumerate}
\end{theorem}

Here, a few details hidden behind the notations have to be mentioned.
The regularized sums are to be interpreted
as in Definition~\ref{def-reg-sum-of-residues},
using the vertex algebras $\hat{H}_\bullet (\calM; \bbQ)$
and $\hat{H}_\bullet (\Mb; \bbQ)$ in
\S\S\ref{sect-background-sheaves}--\ref{sect-background-pairs}.
The choices of representatives in cosets of $\im D$ or $\im \acute{D}$,
as mentioned in Definition~\ref{def-reg-sum-of-residues},
are given by applying the maps $\xi$ and $\acute{\xi}$
in \S\ref{sect-def-xi}.
The regularized sums~\eqref{eq-inv-mss-main-as-reg-sum}
and~\eqref{eq-inv-mss-coprime-as-reg-sum}
are taken in the space $F_r^+$ in Definition~\ref{def-field-inf-many-variables},
although the field $F$ in Definition~\ref{def-regularized-sum}
is taken to be $F_r$.
We will see (Lemma~\ref{lem-z-multiplicative})
that the terms in the regularized sums
lie in a one-dimensional $F_r$-subspace of $F_r^+$,
and Definition~\ref{def-regularized-sum} can be adapted to this case.
Similarly, the regularized sum~\eqref{eq-inv-mbss-main-as-reg-sum}
is taken in $\acute{F}_r^+$,
with the field $\acute{F}_r$ being used in Definition~\ref{def-regularized-sum}.

\begin{remark}
    The regularized sums in Theorem~\ref{thm-inv-mss-main-as-reg-sum}
    cannot be replaced by ordinary sums,
    even if the latter might exist.
    For example, when $(r,d) = (2,1)$ and $g = 1$,
    we can take $d_0 = -i$ and $d_1 = 1+i$,
    and sum over $i = 0,1,2,\dotsc$, giving
    \[
        0 + 0 + 0 + \cdots,
    \]
    which does not produce the expected answer
    $-s_{1,2,2}/2$.
    When $g = 2$, if we ignore the odd variables $s_{j,1,l}$ for simplicity,
    then the sum gives
    \begin{equation}
        \Bigl( -\frac{1}{4} s_{1,2,2}^5 \Bigr) +
        \Bigl( -\frac{3}{4} s_{1,2,2}^5 \Bigr) +
        \Bigl( -\frac{5}{4} s_{1,2,2}^5 \Bigr) + \cdots ,
    \end{equation}
    while the expected answer is
    \begin{equation}
        \frac{1}{8} s_{1,0,2} \, s_{1,2,2}^3
        - \frac{1}{48} s_{1,2,2}^5
        + \frac{1}{16} s_{1,2,2}^3 \, s_{1,2,3}
        - \frac{1}{4} s_{1,2,2} \, s_{1,2,3}^2
        + \frac{11}{48} s_{1,2,2}^2 \, s_{1,2,4} \, .
    \end{equation}
    Note that many new terms ($s_{1,0,2} \, s_{1,2,2}^3$\,, etc.)\ 
    arose in the answer, while they were not present in any of the summands.
    These examples show that it is necessary to interpret the regularized sums 
    in the sense of Definition~\ref{def-reg-sum-of-residues},
    and not as a value associated to the collection of the individual terms.
\end{remark}

\begin{remark}
    The results in Theorem~\ref{thm-inv-mss-main-as-reg-sum},
    expressing the rank~$r$ invariants in terms of the rank~$1$ invariants,
    are a result of an inductive process,
    where in each step, we express the rank~$r$ invariants
    in terms of invariants of rank $< r$.
    The inductive step is stated in
    Theorem~\ref{thm-reduction-main} below.
\end{remark}

We provide a sketch proof of
Theorem~\ref{thm-inv-mss-main-as-reg-sum},
to demonstrate some key ideas
involved in the full proof,
which will be given in \S\ref{sect-proof-of-main-theorem}.

\newcommand{\rankthreeminidiagram}[3][]{%
    \begin{tikzpicture}[baseline={(0,.2)}]
        \fill [black!10] (0:0) -- (#2:.8) arc (#2:#3:.8);
        \draw [black!30, very thick] (#2:.8) -- (0:0) -- (#3:.8);
        \draw [-stealth] (0:0) -- (-150:1);
        \draw [-stealth] (0:0) -- (-30:1);
        \draw [-stealth] (0:0) -- (90:1);
        \node at (-150:1.2) {$d_1$};
        \node at (-30:1.2) {$d_2$};
        \node at (90:1.2) {$d_0$}; #1
    \end{tikzpicture}%
}

\begin{bproof}[Sketch Proof of Theorem~\ref{thm-inv-mss-main-as-reg-sum}]
    We demonstrate the proof in the case $(r, d) = (3, 1)$,
    where it is possible to use pictures to demonstrate some equations.
    Also, as $r$ and $d$ are coprime in this case,
    some combinatorial arguments are simplified.
    
    We start by computing the pair invariant
    $\smash{\inv{\Mbss_{(3,1),1}}}$\,.
    By Theorem~\ref{thm-wcf-1}, one has
    \begin{align*}
        \numberthis
        \label{eq-inv-mbss311-wcf}
        \inv{\Mbss_{(3,1),1}} = {}
        & \sum_{ \substack{ d_{01} + d_2 = 1 \\ d_{01}/2 < d_2 } } {}
        \bigl[
            \inv{\Mbss_{(2,d_{01}),1}} \, , \ 
            \inv{\Mss_{(1,d_2)}}
        \bigr] \\
        & \hspace{2em} {} 
        \mathmakebox[1em][r]{+} \enspace
        \mathmakebox[4em][c]{ \sum_{ \substack{ d_0 + d_{12} = 1 \\ d_0 < d_{12}/2 } } }
        \bigl[
            \inv{\Mbss_{(1,d_0),1}} \, , \ 
            \inv{\Mss_{(2,d_{12})}}
        \bigr] \\
        & \hspace{2em} {}
        \mathmakebox[1em][r]{-} \enspace
        \mathmakebox[4em][c]{ \sum_{ \substack{ d_0 + d_1 + d_2 = 1 \\ d_0 < d_1 < d_2 } } }
        \bigl[ \bigl[
            \inv{\Mbss_{(1,d_0),1}} \, , \ 
            \inv{\Mss_{(1,d_1)}} \bigr] , \ 
            \inv{\Mss_{(1,d_2)}}
        \bigr] \\
        & \hspace{2em} {}
        \mathmakebox[1em][r]{ {} - \frac{1}{2} } \enspace
        \mathmakebox[4em][c]{ \sum_{ \substack{ d_0 + d_1 + d_2 = 1 \\ d_0 < d_1 = d_2 } } }
        \bigl[ \bigl[
            \inv{\Mbss_{(1,d_0),1}} \, , \ 
            \inv{\Mss_{(1,d_1)}} \bigr] , \ 
            \inv{\Mss_{(1,d_2)}}
        \bigr].
    \end{align*}
    We expand the rank~$2$ pair invariants using Theorem~\ref{thm-wcf-1} again,
    and expand the rank~$2$ sheaf invariants using the inductive hypothesis
    (we assume that the theorem is true when $r = 2$),
    so that each term is expressed solely in terms of the rank~$1$ sheaf and pair invariants.
    Then, we use the Jacobi identity to move the pair invariants
    into the innermost layer of the iterated Lie brackets.
    After this process, some terms will cancel out,
    and the final result is
    \begin{equation}
        \label{eq-inv-mbss311-expanded-wcf}
        \inv{\Mbss_{(3,1),1}} = \\
        \sumbar_{(d_0, d_1, d_2) \in \Sigma} {}
        \bigl[ \bigl[
            \inv{\Mbss_{(1,d_0),1}} \, , \ 
            \inv{\Mss_{(1,d_1)}} \bigr] , \ 
            \inv{\Mss_{(1,d_2)}}
        \bigr],
    \end{equation}
    where $\Sigma =
    \{ (d_0, d_1, d_2) \in \bbZ^3 \mid d_0 + d_1 + d_2 = 1, \ 
    d_0 < 1/3, \ d_0 + d_1 < 2/3 \}$,
    as shown in Figure~\ref{fig-sum-region-of-inv-mbss311},
    where the origin is taken to be $(1/3, 1/3, 1/3)$.
    This proves \ref{item-inv-mbss-main-as-reg-sum}.
    
    \begin{figure}
        \centering
        \begin{tikzpicture}
            \fill [black!10] (0:0) -- (-120:1.9) arc (-120:0:1.9);
            \draw [black!30, very thick] (0:1.9) -- (0:0) -- (-120:1.9);
            \draw [-stealth] (0:0) -- (-150:2.2);
            \draw [-stealth] (0:0) -- (-30:2.2);
            \draw [-stealth] (0:0) -- (90:2.2);
            \node at (-150:2.5) {$d_1$};
            \node at (-30:2.5) {$d_2$};
            \node at (90:2.5) {$d_0$};
            \foreach \i in {0, ..., 3} {
                \foreach \j in {0, ..., 3} {
                    \fill (-120:0.5*\j+0.5/3)+(0:0.5*\i+1/3) circle (.04);
                }
            }
            \node at (-60:1.5) {$\Sigma$};
        \end{tikzpicture}
        \caption{Summation region of \eqref{eq-inv-mbss311-expanded-wcf}}
        \label{fig-sum-region-of-inv-mbss311}
    \end{figure}
    
    To simplify notation,
    we write \eqref{eq-inv-mbss311-expanded-wcf} pictorially as
    \begin{equation}
        \label{eq-inv-mbss311-pict}
        \inv{\Mbss_{(3,1),1}} \ = \ 
        \rankthreeminidiagram{-120}{0} \ .
    \end{equation}
    
    On the other hand, we apply Theorem~\ref{thm-wcf-2},
    which says in this case
    \begin{equation}
        \inv{\Mbss_{(3,1),1}} =
        -\bigl[ \upe^{((0,0),1)}, \ \inv{\Mss_{(3,1)}} \bigr].
    \end{equation}
    Taking the Lie bracket with $\upe^{((0,0),1)}$
    is invertible for homology classes of $\Mbpl_{(r,d),0} \simeq \Mpl_{(r,d)}$
    with $r \nu + d > 0$, where $\nu$ is as in Notation~\ref{ntn-n},
    and can be chosen so that $3 \nu + 1 > 0$.
    Therefore, to prove \ref{item-inv-mss-main-as-reg-sum},
    we only need to show that
    \begin{multline}
        \label{eq-inv-mbss311-equals-sum}
        3 \cdot \inv{\Mbss_{(3,1),1}} = \\
        - \sumbar_{(d_0, d_1, d_2) \in \Sigma} {}
        \bigl[ \upe^{((0,0),1)}, \ 
            \bigl[ \bigl[
                \inv{\Mss_{(1,d_0)}} \, , \ 
                \inv{\Mss_{(1,d_1)}} \bigr] , \ 
                \inv{\Mss_{(1,d_2)}}
            \bigr]
        \bigr].
    \end{multline}
    Again, we use the Jacobi identity to move the term $\upe^{((0,0),1)}$
    to the innermost layer of the iterated Lie brackets,
    which shows that the right-hand side of
    \eqref{eq-inv-mbss311-equals-sum} equals
    \begin{align}
        \smash{\sumbar_{(d_0, d_1, d_2) \in \Sigma} {}} &
        \smash{\bigg(}
            \bigl[ \bigl[
                \inv{\Mbss_{(1,d_0),1}} \, , \ 
                \inv{\Mss_{(1,d_1)}} \bigr] , \ 
                \inv{\Mss_{(1,d_2)}}
            \bigr] \notag \\[1ex]
            & \hspace{2em} {} - \bigl[ \bigl[
                \inv{\Mbss_{(1,d_1),1}} \, , \ 
                \inv{\Mss_{(1,d_0)}} \bigr] , \ 
                \inv{\Mss_{(1,d_2)}}
            \bigr] \notag \\[1ex]
            & \hspace{2em} {} - \bigl[ \bigl[
                \inv{\Mbss_{(1,d_2),1}} \, , \ 
                \inv{\Mss_{(1,d_0)}} \bigr] , \ 
                \inv{\Mss_{(1,d_1)}}
            \bigr] \notag \\[1ex]
            & \hspace{2em} {} + \bigl[ \bigl[
                \inv{\Mbss_{(1,d_2),1}} \, , \ 
                \inv{\Mss_{(1,d_1)}} \bigr] , \ 
                \inv{\Mss_{(1,d_0)}}
            \bigr]
        \smash{\bigg)},
        \label{eq-inv-mbss311-sum-expanded}
    \end{align}
    where we used the fact from Theorem~\ref{thm-wcf-2} that
    \begin{equation}
        -\bigl[ \upe^{((0,0),1)}, \ \inv{\Mss_{(1,d)}} \bigr] =
        \inv{\Mbss_{(1,d),1}}
    \end{equation}
    for all $d$.
    Pictorially, \eqref{eq-inv-mbss311-sum-expanded} equals
    \begin{align*}
        & \rankthreeminidiagram{-120}{0}
        - \rankthreeminidiagram{-60}{60}
        - \rankthreeminidiagram{0}{120}
        + \rankthreeminidiagram{60}{180} \\
        = {}
        & \rankthreeminidiagram{-120}{0}
        + \rankthreeminidiagram{-120}{-60}
        + \rankthreeminidiagram{-60}{0}
        + \rankthreeminidiagram{-120}{0} \\
        = {}
        & 3 \cdot \rankthreeminidiagram{-120}{0} \ ,
        \numberthis
        \label{eq-pict-comb-main}
    \end{align*}
    where we used the relations
    \begin{equation}
        \rankthreeminidiagram{-120}{60} =
        \rankthreeminidiagram{-60}{120} =
        \rankthreeminidiagram{0}{180} = 0,
    \end{equation}
    which follow from~\eqref{eq-reg-sum-degen}.
    Now, \eqref{eq-pict-comb-main} together with \eqref{eq-inv-mbss311-pict}
    proves \eqref{eq-inv-mbss311-equals-sum}.
\end{bproof}

Here is the second version of our main result,
where we explicitly compute
regularized sums in Theorem~\ref{thm-inv-mss-main-as-reg-sum}.

\begin{theorem}
    \label{thm-inv-mss-main}
    Let $(r,d) \in K(X)$ with $r>0$.
    \begin{enumerate}
        \item 
        We have
        \upshape
        \begin{multline}
            \label{eq-main-explicit}
            \xi \bigl( \inv{\Mss_{(r,d)}} \bigr) =
            \upe^{(r,d)} \cdot
            \res_{z_{r-1}} \circ \cdots \circ \res_{z_1} 
            \Biggl\{
                \frac{(-1)^{(g-1)r(r-1)/2+(r-1)(d-1)}}{\displaystyle
                    r \cdot \mspace{-9mu}
                    \prod_{0 \leq i < j \leq r-1} \mspace{-9mu}
                    (z_i - z_j)^{2g-2}
                } \cdot {} \\
                \sumbar_{ \leftsubstack{
                    \\[-2ex]
                    & d = d_0 + \cdots + d_{r-1} \\[-.6ex]
                    & (d_0 + \cdots + d_{i-1})/i \leq d/r, \ i = 1, \dotsc, r-1 \\[-.6ex]
                    & \text{with $m$ equalities} 
                } }
                \frac{1}{m+1} \cdot
                \prod_{i=0}^{r-1} {}
                \bigl[
                    \exp \bigl( \tilde{z}_i \,
                        D_{(1, d_i)}
                    \bigr) \, \sigma (-s_{1,2,2})
                \bigr]
            \Biggr\} \Bigg|_{z_0=0}\,,
        \end{multline}
        \itshape
        where
        $\xi$ is the map defined in \textnormal{\S\ref{sect-def-xi}},
        $D_{(1,d_1)}$ is given by \eqref{eq-def-d-r-d},
        $\sigma$ is defined by \eqref{eq-def-sigma},
        and we write
        $\tilde{z}_i = z_i - (z_0 + \cdots + z_{r-1})/r$.
        Explicitly, we have
        \upshape
        \begin{multline}
            \label{eq-main-explicit-summed}
            \xi \bigl( \inv{\Mss_{(r,d)}} \bigr) =
            \upe^{(r,d)} \cdot
            \res_{z_{r-1}} \circ \cdots \circ \res_{z_1}
            \Biggl\{
                \frac{(-1)^{(g-1)r(r-1)/2+(r-1)(d-1)}}{\displaystyle
                    r \cdot \mspace{-9mu}
                    \prod_{0 \leq i < j \leq r-1} \mspace{-9mu}
                    (z_i - z_j)^{2g-2}
                } \cdot {} \\
                \sum_{ \leftsubstack[7em]{
                    \\[-2ex]
                    & 0 \leq m \leq \gcd(r,d)-1 \\[-.6ex]
                    & 1 \leq i_1 < \cdots < i_m \leq r-1 \\[-.6ex]
                    & \text{such that } i_k d / r \in \bbZ \text{ for all } k
                } }
                \frac{(-1)^m}{m+1} \cdot 
                \frac{1}{\displaystyle
                    \prod_{ \leftsubstack[4em]{
                        & 1 \leq i \leq r-1 \\[-.6ex]
                        & i \neq i_k \text{ for any } k
                    } } \biggl[ 
                        1 - \exp \biggl(
                            \sum_{l=1}^{\infty}
                            \frac{\tilde{z}_i^l - \tilde{z}_{i-1}^l}{l!} \,
                            s_{1,2,l+1}
                        \biggr)
                    \biggr]
                } \cdot {} \\
                \prod_{i=0}^{r-1} {}
                \Bigl[
                    \exp \Bigl( \tilde{z}_i \,
                        D_{\left(1, \,
                            \lfloor\mspace{-3mu} \frac{(i+1)d}{r} \mspace{-3mu}\rfloor -
                            \lfloor\mspace{-3mu} \frac{id}{r} \mspace{-3mu}\rfloor
                        \right)}
                    \Bigr) \, \sigma (-s_{1,2,2})
                \Bigr]
            \Biggr\} \Bigg|_{z_0=0}\,.
        \end{multline}
        \itshape
    
        \item We have
        \begin{multline}
            \label{eq-main-pairs-explicit}
            \acute{\xi} \bigl( \inv{\Mbss_{(r,d),1}} \bigr) =
            \upe^{((r,d),1)} \cdot
            \res_{z_{r-1}} \circ \cdots \circ \res_{z_1} \circ \res_w \\
            \shoveleft{
            \Biggl\{ 
                \frac{(-1)^{(g-1)r(r-1)/2+(r-1)(d-1)}}{\displaystyle
                    \prod_{0 \leq i < j \leq r-1} \mspace{-9mu}
                    (z_i - z_j)^{2g-2} \cdot
                    \prod_{i=0}^{r-1} {} (w + z_i)^{\nu + d_i}
                } \cdot
                \rho (w) \cdot {} } \\
                \sumbar_{ \leftsubstack[6em]{
                    \\[-2ex]
                    & d = d_0 + \cdots + d_{r-1} \\[-.6ex]
                    & (d_0 + \cdots + d_{i-1})/i < d/r, \ i = 1, \dotsc, r-1
                } }
                \prod_{i=0}^{r-1} {}
                \Bigl[
                    \exp \bigl( \tilde{z}_i \,
                        D_{(1, d_i)}
                    \bigr) \, \sigma \Bigl( \frac{1}{w+z_i}-s_{1,2,2} \Bigr)
                \Bigr]
            \Biggr\} \Bigg|_{z_0=0}\,,
        \end{multline}
        \itshape
        where
        $\acute{\xi}$ is the map defined in \textnormal{\S\ref{sect-def-xi}},
        $\sigma$ and $\rho$ are defined by \eqref{eq-def-sigma} and \eqref{eq-def-rho},
        $\nu$ is as in Notation~\textnormal{\ref{ntn-n}},
        and we write
        $\tilde{z}_i = z_i - (z_0 + \cdots + z_{r-1})/r$.
        Explicitly, we have
        \upshape
        \begin{multline}
            \label{eq-main-pairs-explicit-summed}
            \acute{\xi} \bigl( \inv{\Mbss_{(r,d),1}} \bigr) =
            \upe^{((r,d),1)} \cdot
            \res_{z_{r-1}} \circ \cdots \circ \res_{z_1} \circ \res_w \\
            \shoveleft{
            \Biggl\{ 
                \frac{(-1)^{(g-1)r(r-1)/2+(r-1)(d-1)}}{\displaystyle
                    \prod_{0 \leq i < j \leq r-1} \mspace{-9mu}
                    (z_i - z_j)^{2g-2} \cdot
                    \prod_{i=0}^{r-1} {} (w + z_i)^{\nu + d_i}
                } \cdot
                \rho (w) \cdot {} } \\
                \frac{1}{\displaystyle
                    \prod_{i=1}^{r-1} {} \biggl[ 
                        1 - \exp \biggl(
                            \sum_{l=1}^{\infty}
                            \frac{\tilde{z}_i^l - \tilde{z}_{i-1}^l}{l!} \,
                            s_{1,2,l+1}
                        \biggr)
                    \biggr]
                } \cdot {} \\
                \prod_{i=0}^{r-1} {}
                \Bigl[
                    \exp \Bigl( \tilde{z}_i \,
                        D_{\left(1, \,
                            \lceil\mspace{-3mu} \frac{(i+1)d}{r} \mspace{-3mu}\rceil -
                            \lceil\mspace{-3mu} \frac{id}{r} \mspace{-3mu}\rceil
                            + \delta_{i,r-1} - \delta_{i,0} 
                        \right)}
                    \Bigr) \, \sigma \Bigl( \frac{1}{w+z_i}-s_{1,2,2} \Bigr)
                \Bigr]
            \Biggr\} \Bigg|_{z_0=0}\,,
        \end{multline}
        \itshape
        where $\delta_{i,j}$ is the Kronecker delta.
    \end{enumerate}
\end{theorem}

The proof will be given in \S\ref{sect-proof-of-main-theorem}.
Note that the expressions
\eqref{eq-main-explicit-summed} and \eqref{eq-main-pairs-explicit-summed}
do not involve any regularized sums.

It is possible to extract intersection pairings on $\Mss_{(r,d)}$
from the expression \eqref{eq-main-explicit-summed}.

\begin{theorem}
    \label{thm-inv-mss-pairing}
    Let $(r,d) \in K(X)$ with $r>0$. Suppose that
    \begin{itemize}
        \item 
            $( m_l )_{l \geq 2}$
            is a family of non-negative integers,
            of which only finitely many are non-zero.
        \item 
            $( p_{j,l} )_{1 \leq j \leq 2g, \, l \geq 1}$
            is a family of integers that are either $0$ or $1$,
            of which only finitely many are non-zero.
        \item 
            $( \alpha_l )_{l \geq 2}$
            is a sequence of formal variables.
    \end{itemize}
    Then we have the intersection pairing
    \upshape
    \begin{multline}
        \label{eq-main-pairing}
        \int_{\inv{\Mss_{(r,d)}}}
        \Xi \biggl[
            \prod_{l=2}^{\infty} S_{1,0,l}^{m_l} \cdot
            \prod_{j=1}^{2g} \prod_{l=1}^{\infty} S_{j,1,l}^{p_{j,l}} \cdot
            \exp \biggl( \sum_{l=2}^{\infty} \alpha_l \, S_{1,2,l} \biggr)
        \biggr] = \\
        \shoveleft{ \quad \res_{z_{r-1}} \circ \cdots \circ \res_{z_1}
        \Biggl\{
            \frac{(-1)^{(g-1)r(r-1)/2+(r-1)(d-1)}}{\displaystyle
                r \cdot \mspace{-9mu}
                \prod_{0 \leq i < j \leq r-1} \mspace{-9mu}
                (z_i - z_j)^{2g-2}
            } \cdot {} } \\
            \Biggl[ \ 
                \sum_{ \leftsubstack[7em]{
                    \\[-2ex]
                    & 0 \leq m \leq \gcd(r,d)-1 \\[-.6ex]
                    & 1 \leq i_1 < \cdots < i_m \leq r-1 \\[-.6ex]
                    & \text{such that } i_k d / r \in \bbZ \text{ for all } k
                } }
                \frac{(-1)^m}{m+1} \cdot 
                \frac{1}{\displaystyle
                    \prod_{ \leftsubstack[4em]{
                        & 1 \leq i \leq r-1 \\[-.6ex]
                        & i \neq i_k \text{ for any } k
                    } } \biggl[ 
                    1 - \exp \biggl(
                        \sum_{l=1}^{\infty}
                        \frac{\tilde{z}_i^l - \tilde{z}_{i-1}^l}{l!} \,
                        \alpha_{l+1}
                    \biggr)
                \biggr]
                }
            \Biggr] \cdot {} \\
            \prod_{l=2}^{\infty} {}
            \biggl( \sum_{i=0}^{r-1} {} \frac{\tilde{z}_i^l}{l!} \biggr)^{m_l} \cdot
            \exp \biggl( 
                \sum_{i=0}^{r-1} {} \Bigl(
                    \Bigl\lfloor \frac{(i+1)d}{r} \Bigr\rfloor - 
                    \Bigl\lfloor \frac{id}{r} \Bigr\rfloor
                \Bigr)
                \sum_{l=1}^{\infty} \frac{\tilde{z}_i^l}{l!} \,
                \alpha_{l+1}
            \biggr) \cdot {} \\
            \sum_{(i_{j,l})} {} (-1)^* \cdot
            \biggl(
                \prod_{ \substack{ (j,l) \text{ with} \\ p_{j,l} = 1 } }
                \frac{\tilde{z}_{i_{j,l}}^{l-1}}{(l-1)!}
            \biggr) \cdot
            \prod_{i=0}^{r-1} {} \biggl(
                -\sum_{l=0}^{\infty}
                \frac{\tilde{z}_i^l}{l!} \, \alpha_{l+2}
            \biggr)^{g-e_i}
        \Biggr\} \Bigg|_{z_0=0}\,,
    \end{multline}
    \itshape
    where $\Xi$ is the map defined in \textnormal{\S\ref{sect-def-xi}},
    and $\tilde{z}_i = z_i - (z_0 + \cdots + z_{r-1})/r$.
    The sum on the last row is over all possible ways to choose
    $0 \leq i_{j,l} \leq r-1$ for each $(j,l)$ with $p_{j,l} = 1$,
    such that for all $0 \leq i \leq r-1$ and all $1 \leq j \leq g$,
    either $i \neq i_{j,l}$ and $i \neq j_{j+g,l'}$ for all $l, l'$,
    or $i = i_{j,l}$ and $i = i_{j+g},l'$ for a unique pair $(l, l')$.
    For each $i$, define $e_i$ to be the number of $j$ in the second case above. 
    The sign $(-1)^*$ is defined as follows.
    For $(j,l)$ and $(j',l')$ with $p_{j,l} = p_{j',l'} = 1$,
    define $(j,l) \preceq_1 (j',l')$ if $j < j'$, or $j = j'$ and $l \leq l'$.
    Define $(j,l) \preceq_2 (j',l')$ if $i_{j,l} < i_{j',l'}$,
    or $i_{j,l} = i_{j',l'}$ and $(j \operatorname{mod} g) < (j' \operatorname{mod} g)$,
    or $i_{j,l} = i_{j',l'}$ and $j'=j+g$ or $j$.
    Then, $(-1)^*$ is the sign of the permutation
    on the set of all $(j,l)$ with $p_{j,l} = 1$
    that intertwines the total orders $\preceq_1$ and~$\preceq_2$.
    
    In particular, we have
    \upshape
    \begin{multline}
        \label{eq-main-pairing-simpler}
        \int_{\inv{\Mss_{(r,d)}}}
        \Xi \biggl[
            \prod_{l=2}^{\infty} S_{1,0,l}^{m_l} \cdot
            \exp ( \alpha \, S_{1,2,2} )
        \biggr] = \\
        \shoveleft{ \quad \res_{z_{r-1}} \circ \cdots \circ \res_{z_1}
        \Biggl\{
            \frac{(-1)^{(g-1)r(r-1)/2+(r-1)(d-1)} \cdot (-\alpha)^{rg}}
            {\displaystyle
                r \cdot \mspace{-9mu}
                \prod_{0 \leq i < j \leq r-1} \mspace{-9mu}
                (z_i - z_j)^{2g-2}
            } \cdot {} } \\
            \Biggl[ \ 
                \sum_{ \leftsubstack[7em]{
                    \\[-2ex]
                    & 0 \leq m \leq \gcd(r,d)-1 \\[-.6ex]
                    & 1 \leq i_1 < \cdots < i_m \leq r-1 \\[-.6ex]
                    & \text{such that } i_k d / r \in \bbZ \text{ for all } k
                } }
                \frac{(-1)^m}{m+1} \cdot 
                \frac{1}{\displaystyle
                    \prod_{ \leftsubstack[4em]{
                        & 1 \leq i \leq r-1 \\[-.6ex]
                        & i \neq i_k \text{ for any } k
                    } } \bigl( 
                    1 - \exp (\alpha (z_i - z_{i-1}))
                \bigr)
                }
            \Biggr] \cdot {} \\
            \prod_{l=2}^{\infty} {}
            \biggl( \sum_{i=0}^{r-1} {} \frac{\tilde{z}_i^l}{l!} \biggr)^{m_l} \cdot
            \exp \biggl( 
                \alpha \sum_{i=0}^{r-1} 
                \tilde{d_i} \, z_i
            \biggr)
        \Biggr\} \Bigg|_{z_0=0}\,,
    \end{multline}
    \itshape
    where $\tilde{d}_i = \lfloor (i+1)d/r \rfloor - \lfloor id/r \rfloor - d/r$,
    and $\alpha$ is a formal variable.
    Note that $\tilde{d}_i$ only depends on $d \operatorname{mod} r$,
    and so does the pairing \eqref{eq-main-pairing-simpler}.
\end{theorem}

\begin{proof}
    We use the fact that for a polynomial $f \in \bbC [x]$, one has
    \[
        \exp \Bigl( a \, \frac{\partial}{\partial x} \Bigr) \, f (x) \Big|_{x=0} =
        f (a)
    \]
    for any $a \in \bbC$.
    Since we know that the result of the repeated residue
    is a polynomial in the $s_{j,k,l}$ variables,
    we may apply this fact to the variables $s_{1,2,l}$\,,
    which means that the pairing is obtained by
    replacing $s_{1,2,l}$ with $\alpha_l$ in \eqref{eq-main-explicit-summed},
    and then taking the coefficient of the monomial
    $\smash{\prod_{l=2}^{\infty} s_{1,0,l}^{m_l} \cdot
    \prod_{j=1}^{2g} \prod_{l=1}^{\infty} s_{j,1,l}^{p_{j,l}}}$.
    
    Writing $d_i = \lfloor (i+1)d/r \rfloor - \lfloor id/r \rfloor$, we have
    \begin{align*}
        & \exp (\tilde{z}_i \, D_{(1,d_i)}) \,
        \sigma (-s_{1,2,2}) \\
        = {} &
        \exp \biggl( \sum_{l=1}^{\infty}
            \frac{\tilde{z}_i^l}{l!} 
            \, ( s_{1,0,l} + d_i \, s_{1,2,l+1} )
        \biggr) \, \exp (\tilde{z}_i \, \Delta) \, \sigma (-s_{1,2,2}) \\
        = {} &
        \exp \biggl( \sum_{l=1}^{\infty}
            \frac{\tilde{z}_i^l}{l!} 
            \, s_{1,0,l}
        \biggr) \,
        \exp \biggl( d_i \sum_{l=1}^{\infty}
            \frac{\tilde{z}_i^l}{l!} 
            \, s_{1,2,l+1}
        \biggr) \cdot {} \\
        & \hspace{2em}
        \exp (\tilde{z}_i \, \Delta) \, 
        \biggl( \,
            \sum_{ \substack{ 0 \leq e \leq g \\ 1 \leq j_1 < \cdots < j_e \leq g } }
            \mspace{-12mu}
            (-s_{1,2,2})^{g-e} \, s_{j_1,1,1} \, s_{j_1+g,1,1} \cdots s_{j_e,1,1} \, s_{j_e+g,1,1}
        \biggr) ,
    \end{align*}
    where $\Delta = D_{(0,0)}$\,.
    The last line of the expression 
    (starting with $\exp (\tilde{z}_i \, \Delta)$)
    can be expanded as
    \[
        \sum_{ \substack{ 0 \leq e \leq g \\ 1 \leq j_1 < \cdots < j_e \leq g } }
        \biggl(
            -\sum_{l=0}^{\infty}
            \frac{\tilde{z}_i^l}{l!} \, s_{1,2,l+2}
        \biggr)^{g-e} \cdot
        \prod_{h=1}^e {}
        \biggl[
            \biggl(
                \sum_{l=0}^{\infty}
                \frac{\tilde{z}_i^l}{l!} \, s_{j_h,1,l+1}
            \biggr)
            \biggl(
                \sum_{l=0}^{\infty}
                \frac{\tilde{z}_i^l}{l!} \, s_{j_h+g,1,l+1}
            \biggr)
        \biggr].
    \]
    If we take the coefficient of 
    $\smash{\prod_{j=1}^{2g} \prod_{l=1}^{\infty} s_{j,1,l}^{p_{i,j,l}}}$,
    we obtain
    \[
        (-1)^* \cdot
        \biggl(
            -\sum_{l=0}^{\infty}
            \frac{\tilde{z}_i^l}{l!} \, s_{1,2,l+2}
        \biggr)^{g-e} \cdot
        \prod_{(j,l) \text{ with } p_{i,j,l} = 1} {}
        \frac{\tilde{z}_i^{l-1}}{(l-1)!},
    \]
    where the sign $(-1)^*$ is caused by
    changing the order of the odd variables $s_{j,1,l}$\,,
    provided that the following condition is satisfied:
    for all $1 \leq j \leq g$,
    either $p_{i,j,l} = p_{i,j+g,l'} = 0$ for all $l,l'$,
    or $p_{i,j,l} = p_{i,j+g,l'} = 1$ for unique $l,l'$,
    and $e$ is the number of $j$ in the second case.
    Otherwise, the coefficient is $0$.
    From this observation,
    we obtain \eqref{eq-main-pairing} by setting
    $p_{i,j,l} = 1$ if and only if $i_{j,l} = i$.
\end{proof}

\subsection{The fixed determinant invariants}
\label{sect-fixed-determinant}

Next, using the results of the previous section,
we define and study the invariants $\inv{\Mssfd_{(r,d)}}$
for the fixed determinant moduli space $\Mssfd_{(r,d)}$.

For $(r,d) \in K(X)$ with $r > 0$, consider the map
\begin{equation}
    \label{eq-def-det}
    \det^\pl \colon \Mpl_{(r,d)} \longrightarrow \Mss_{(1,d)} ,
\end{equation}
which takes a coherent sheaf to its determinant line bundle.
Let $\Mssfd_{(r,d)}$ be a fibre of the restriction
of this map to $\Mss_{(r,d)}$.

We define the invariant
\begin{align}
    \label{eq-def-inv-mssfd}
    \inv{\Mssfd_{(r,d)}} & =
    \inv{\Mss_{(r,d)}} \cap
    (\det^\pl)^* (\omega) \\
    & \in H_{2(r^2-1)(g-1)} (\Mpl_{(r,d)}; \bbQ), \notag
\end{align}
where $\omega \in H^{2g} (\Mss_{(1,d)}; \bbQ)$
is the Poincar\'e dual of a point.

If $r$ and $d$ are coprime,
then $\Mssfd_{(r,d)}$ is a smooth, projective variety, and 
\begin{equation}
    i_* \fund{\Mssfd_{(r,d)}} =
    \inv{\Mssfd_{(r,d)}},
\end{equation}
where $i \colon \Mssfd_{(r,d)} \hookrightarrow \Mpl_{(r,d)}$
is the inclusion.

\begin{lemma}
    \label{lem-fd-factor}
    In \eqref{eq-def-inv-mssfd}, we may take
    \begin{equation}
        \omega = \Xi \biggl( {}
            \prod_{j=1}^{g} S_{j,1,1} \, S_{j+g,1,1}
        \biggr).
    \end{equation}
\end{lemma}

\begin{proof}
    By Lemma~\ref{lem-pb-of-sjkl-to-jacobian},
    we see that $\Xi (S_{j,1,1}) = -\gamma_j$\,,
    where $\gamma_j$ is as in that Lemma.
    The result then follows.
\end{proof}

\begin{lemma}
    \label{lem-det-star}
    Consider the map $\det^\pl$ in
    \eqref{eq-def-det}. Then
    \begin{align}
        (\det^\pl)^* \, \Xi (S_{j,1,1}) & = \Xi (S_{j,1,1}), \\
        (\det^\pl)^* \, \Xi (S_{1,2,2}) & =
        \displaystyle \Xi \biggl(
            -\sum_{j=1}^g S_{j,1,1} \, S_{j+g,1,1}
        \biggr),
    \end{align}
    and $(\det^\pl)^* \, \Xi (S_{j,k,l}) = 0$
    for all other $(j,k,l)$.
\end{lemma}

\begin{proof}
    Consider also the map
    \begin{align}
        \det \colon \calM_{(r,d)}
        & \longrightarrow \calM_{(1,d)}
    \end{align}
    sending a complex to its determinant line bundle.
    Then we have
    \begin{align*}
        \det^* S_{j,k,l} & =
        \det^* \ch_l (\calU_{(1,d)}) \setminus e_{j,k}
        \\
        & = \ch_l (\det \calU_{(r,d)}) \setminus e_{j,k} \, ,
        \numberthis
    \end{align*}
    where $\calU_{(r,d)} \to X \times \calM_{(r,d)}$
    is the universal perfect complex.
    Since
    \begin{equation}
        \ch (\det \calU_{(r,d)}) =
        \exp \ch_1 (\calU_{(r,d)}),
    \end{equation}
    we see that
    \begin{equation}
        \label{eq-det-star-sjkl}
        \det^* S_{j,k,l} =
        \frac{1}{l!} \biggl(
            \sum_{(j',k') \in J}
            \epsilon_{j',k'} \boxtimes S_{j',k',1}
        \biggr)^l \setminus e_{j,k} \, ,
    \end{equation}
    and in particular, $\det^* S_{1,0,1} = S_{1,0,1}$.
    By Corollary~\ref{cor-xi-is-pullback}, we have
    \begin{equation}
        (\det^\pl)^* \circ \Xi =
        \Xi \circ \det^*,
    \end{equation}
    where we set $\Xi (S_{1,0,1}) = 0$.
    Applying $\Xi$ to both sides of
    \eqref{eq-det-star-sjkl}, we obtain
    \begin{equation}
        \label{eq-det-pl-star-sjkl}
        (\det^\pl)^* \, \Xi (S_{j,k,l}) =
        \frac{1}{l!} \biggl(
            \sum_{ \substack{ (j',k') \in J \\ (j',k') \neq (1,0) } }
            \epsilon_{j',k'} \boxtimes \Xi (S_{j',k',1})
        \biggr)^l \setminus e_{j,k} \, ,
    \end{equation}
    which implies the desired results.
\end{proof}

\begin{theorem}
    \label{thm-inv-mssfd}
    For $(r, d) \in K (X)$ such that $r > 0$, we have
    \begin{equation}
        \xi \bigl( \inv{\Mssfd_{(r,d)}} \bigr) =
        \biggl(
            \prod_{j=1}^g
            \frac{\partial}{\partial s_{j+g,1,1}}
            \frac{\partial}{\partial s_{j,1,1}}
        \biggr) \,
        \xi \bigl( \inv{\Mss_{(r,d)}} \bigr) \, .
    \end{equation}
\end{theorem}

\begin{proof}
    This follows from Lemmas~\ref{lem-fd-factor}
    and~\ref{lem-det-star}.
\end{proof}

\begin{theorem}
    \label{thm-inv-mssfd-pairing}
    Let $(r,d) \in K(X)$ with $r>0$.
    Using the notations in Theorem~\textnormal{\ref{thm-inv-mss-pairing},}
    we have the intersection pairing
    \begin{multline}
        \label{eq-main-fd-pairing}
        \int_{\inv{\Mssfd_{(r,d)}}}
        \Xi \biggl[
            \prod_{l=2}^{\infty} S_{1,0,l}^{m_l} \cdot
            \prod_{j=1}^{2g} \prod_{l=2}^{\infty} S_{j,1,l}^{p_{j,l}} \cdot
            \exp \biggl( \sum_{l=2}^{\infty} \alpha_l \, S_{1,2,l} \biggr)
        \biggr] = \\
        \int_{\inv{\Mss_{(r,d)}}}
        \Xi \biggl[
            \prod_{l=2}^{\infty} S_{1,0,l}^{m_l} \cdot
            \prod_{j=1}^{g} {} (S_{j,1,1} \, S_{j+g,1,1}) \cdot
            \prod_{j=1}^{2g} \prod_{l=2}^{\infty} S_{j,1,l}^{p_{j,l}} \cdot
            \exp \biggl( \sum_{l=2}^{\infty} \alpha_l \, S_{1,2,l} \biggr)
        \biggr],
    \end{multline}
    where $\Xi$ is the map defined in \textnormal{\S\ref{sect-def-xi}},
    and the right-hand side is given by
    Theorem~\textnormal{\ref{thm-inv-mss-pairing}}.
    In particular, we have
    \upshape
    \begin{multline}
        \label{eq-main-fd-pairing-simpler}
        \int_{\inv{\Mssfd_{(r,d)}}}
        \Xi \biggl( \prod_{l=2}^{\infty} S_{1,0,l}^{m_l} \cdot
        \exp (\alpha \, S_{1,2,2}) \biggr) = \\
        \shoveleft{ \quad \res_{z_{r-1}} \circ \cdots \circ \res_{z_1}
        \Biggl\{
            \frac{(-1)^{(g-1)r(r-1)/2+(r-1)(d-1)} \cdot r^{g-1} \cdot (-\alpha)^{(r-1)g}}
            {\displaystyle
                \prod_{0 \leq i < j \leq r-1} \mspace{-9mu}
                (z_i - z_j)^{2g-2}
            } \cdot {} } \\
            \Biggl[ \ 
                \sum_{ \leftsubstack[8em]{
                    \\[-2ex]
                    & 0 \leq m \leq \gcd(r,d)-1 \\[-.6ex]
                    & 1 \leq i_1 < \cdots < i_m \leq r-1 \\[-.6ex]
                    & \text{such that } i_k d / r \in \bbZ \text{ for all } k
                } }
                \frac{(-1)^m}{m+1} \cdot 
                \frac{1}{\displaystyle
                    \prod_{ \leftsubstack[4em]{
                        & 1 \leq i \leq r-1 \\[-.6ex]
                        & i \neq i_k \text{ for any } k
                    } } \bigl(
                        1 - \exp (\alpha (z_i - z_{i-1}))
                    \bigr)
                }
            \Biggr] \cdot {} \\
            \prod_{l=2}^{\infty} {}
            \biggl( \sum_{i=0}^{r-1} {} \tilde{z}_i^l \biggr)^{m_l} \cdot
            \exp \biggl( \alpha \sum_{i=0}^{r-1} \tilde{d}_i \, z_i \biggr)
        \Biggr\} \Bigg|_{z_0=0}\,,
    \end{multline}
    \itshape
    where $\tilde{d}_i = \lfloor (i+1)d/r \rfloor - \lfloor id/r \rfloor - d/r$,
    and $\alpha$ is a formal variable.
    Note that $\tilde{d}_i$ only depends on $d \operatorname{mod} r$,
    and so does the pairing \eqref{eq-main-fd-pairing-simpler}.
\end{theorem}

\begin{proof}
    The first formula follows from Theorem~\ref{thm-inv-mssfd}.
    The second one follows from the first one
    together with Theorem~\ref{thm-inv-mss-pairing}.
\end{proof}

\begin{definition}
    \label{def-symp-vol}
    Let $(r,d) \in K(X)$ with $r>0$.
    We define the \emph{symplectic volume} of $\Mssfd_{(r,d)}$ by
    \begin{equation}
        \label{eq-def-symp-vol}
        \mathrm{vol} (\Mssfd_{(r,d)}) =
        \int_{\inv{\Mssfd_{(r,d)}}} \exp (-S_{1,2,2}).
    \end{equation}
    This is explicitly given by
    \eqref{eq-main-fd-pairing-simpler} with $\alpha = -1$.

    In particular, when $r$ and $d$ are coprime and $g > 1$,
    this coincides with the usual notion of the symplectic volume
    with symplectic form $-S_{1,2,2}$\,,
    or $f_2$ in the notation of Jeffrey--Kirwan~\cite{jeffrey-kirwan-2}.
\end{definition}

Some numeric values of the symplectic volume when $g = 2$
for small $r$ and $d$ was given in \S\ref{sect-intro}.

\subsection{Results for rank 2}
\label{sect-rank-2}

Next, we apply the results in the previous sections
in the case when $r = 2$,
where simpler formulae can be obtained.

\begin{theorem}
    \label{thm-inv-mss2d}
    One has
    \upshape
    \begin{multline}
        \label{eq-inv-mss2d-after-sum}
        \xi \bigl( \inv{\Mss_{(2,d)}} \bigr) =
        \upe^{(2,d)} \cdot
        \res_z \Biggl\{
            \frac{(-1)^{g+d}} {2 \, z^{2g-2}} \cdot 
            \frac{1}{\displaystyle
                1 -
                \exp \biggl(
                    \sum_{l \geq 1 \text{ odd}}
                    \frac{z^l}{2^{l-1} \, l!} 
                    \, s_{1,2,l+1}
                \biggr)
            } \cdot {} \\
            \biggl(
                \exp \Bigl( -\frac{1}{2} z D_{(1,\lfloor d/2 \rfloor)} \Bigr) \,
                \sigma (-s_{1,2,2})
            \biggr) \cdot
            \biggl(
                \exp \Bigl( \frac{1}{2} z D_{(1,\lceil d/2 \rceil)} \Bigr) \,
                \sigma (-s_{1,2,2})
            \biggr) 
        \Biggr\} .
    \end{multline}
\end{theorem}

\begin{proof}
    If $d$ is odd, this is precisely the expression given by
    Theorem~\ref{thm-inv-mss-main}.
    If $d$ is even, we have omitted the term with $m = 1$ in
    \eqref{eq-main-explicit-summed}.
    But this term becomes zero after taking the residue,
    as it represents the Lie bracket
    \begin{equation}
        -\frac{1}{2}
        \bigl[ \inv{\Mss_{(1,d/2)}} \, , \ \inv{\Mss_{(1,d/2)}} \bigr] = 0.
    \end{equation}
\end{proof}

One can also obtain formulae for intersection pairings.

\begin{theorem}
    \label{thm-inv-mss2d-explicit}
    Let $( m_l )_{l \geq 2}$ be a sequence of non-negative integers,
    of which only finitely many are non-zero,
    and let $h = \sum_l l m_l$.
    \begin{enumerate}
        \item
            If $h \leq 2g-2$ and
            $m_l = 0$ whenever $l$ is odd, then
            \begin{multline}
                \label{{eq-pairing-mss-2-d}}
                \int_{\inv{\Mss_{(2,d)}}}
                \Xi \biggl(
                    \prod_{l=2}^\infty S_{1,0,l}^{m_l} \cdot S_{1,2,2}^{4g-3-h}
                \biggr)
                = \\
                \frac {(-1)^{g+d-1} \, (4g-3-h)!}
                {2^{h-\sum_l m_l+1} \prod_l l!^{m_l} \cdot (2g-2-h)!} \,
                B_{2g-2-h} \Bigl(
                    \Bigl\{ \frac{d}{2} \Bigr\}
                \Bigr),
            \end{multline}
            where $B_n (x)$ is the $n$-th Bernoulli polynomial,
            and $\{ d/2 \}$ denotes the fractional part of $d/2$.
            Otherwise, this pairing will be zero.
        \item
            If $h \leq 2g-2$ and
            $m_l = 0$ whenever $l$ is odd, then
            \begin{multline}
                \label{eq-pairing-mssfd-2-d}
                \int_{\inv{\Mssfd_{(2,d)}}}
                \Xi \biggl(
                    \prod_{l=2}^\infty S_{1,0,l}^{m_l} \cdot S_{1,2,2}^{3g-3-h}
                \biggr)
                = \\
                \frac {(-1)^{d-1} \, (3g-3-h)!}
                {2^{h-\sum_l m_l+1-g} \prod_l l!^{m_l} \cdot (2g-2-h)!} \,
                B_{2g-2-h} \Bigl(
                    \Bigl\{ \frac{d}{2} \Bigr\}
                \Bigr).
            \end{multline}
            Otherwise, this pairing will be zero.
    \end{enumerate}
\end{theorem}

\begin{proof}
    For (i), by Theorem~\ref{thm-inv-mss-pairing},
    the pairing equals \eqref{eq-main-pairing-simpler} with $\alpha = 1$,
    multiplied by $(4g-3-h)!$.
    Theorem~\ref{thm-inv-mss2d} means that we can omit
    the term with $m = 1$ in \eqref{eq-main-pairing-simpler}.
    Writing $z$ for $z_1$, and
    substituting with $\tilde{z}_0 = -z/2$, $\tilde{z}_1 = z/2$,
    $\tilde{d}_0 = -\{ d/2 \}$ and $\tilde{d}_1 = \{ d/2 \}$, we obtain
    \begin{multline}
        \label{eq-inv-mss2d-pairing-1}
        \int_{\inv{\Mss_{(2,d)}}}
        \Xi \biggl[
            \prod_{l=2}^{\infty} S_{1,0,l}^{m_l} \cdot
            S_{1,2,2}^{4g-3-h}
        \biggr] = \\
        \res_{z} \Biggl\{
            \frac{(-1)^{g+d} \, (4g-3-h)!}
            {\displaystyle
                2 \, z^{2g-2} \, (1 - \exp z)
            } \cdot
            \prod_{l \geq 2 \text{ even}} {}
            \Bigl( \frac{z^l}{2^{l-1} \, l!} \Bigr)^{m_l} \cdot
            \exp \Bigl( 
                \Bigl\{ \frac{d}{2} \Bigr\} \, z
            \Bigr)
        \Biggr\},
    \end{multline}
    under the assumption that $m_l = 0$ whenever $l$ is odd.
    By the definition of the Bernoulli polynomial,
    we have the power series expansion
    \begin{equation}
        \frac{z \exp (t z)}{1 - \exp z} =
        -\sum_{n=0}^{\infty} B_n (t) \, \frac{z^n}{n!}
    \end{equation}
    for any $t \in \bbC$. Taking $t = \{ d/2 \}$,
    this reduces the right-hand side of \eqref{eq-inv-mss2d-pairing-1} to
    \begin{align*}
        & \res_{z} \Biggl\{
            \frac{(-1)^{g+d-1} \, (4g-3-h)!}
            {\displaystyle
                2 \, z^{2g-2}
            } \cdot
            \prod_{l \geq 2 \text{ even}} {}
            \Bigl( \frac{z^l}{2^{l-1} \, l!} \Bigr)^{m_l} \cdot
            \sum_{n=0}^{\infty} B_n \Bigl( \Bigl\{ \frac{d}{2} \Bigr\} \Bigr) \,
            \frac{z^{n-1}}{n!}
        \Biggr\} \\
        = {} &
        \frac {(-1)^{g+d-1} \, (4g-3-h)!}
        {2^{h-\sum_l m_l+1} \prod_l l!^{m_l} \cdot (2g-2-h)!} \,
        B_{2g-2-h} \Bigl(
            \Bigl\{ \frac{d}{2} \Bigr\}
        \Bigr).
        \numberthis
    \end{align*}
    
    The proof of (ii) is similar,
    where we use \eqref{eq-main-fd-pairing-simpler}
    instead of \eqref{eq-main-pairing-simpler}.
\end{proof}

\begin{remark}
    \label{rem-inv-mss21-equiv-to-jeffrey-kirwan}
    In Theorem~\ref{thm-inv-mss2d-explicit},
    If we take $d=1$ and $m_l = 0$ for $l \neq 2$ in (ii),
    using that
    \begin{equation}
        \frac{1}{n!} \, B_n \Bigl( \frac{1}{2} \Bigr) =
        \frac{(-1)^{n/2} \, (2^{n-1}-1)} {2^{2n-2} \, \uppi^n} \, \zeta (n)
    \end{equation}
    for $n \geq 0$ even, where $\zeta$ is the Riemann zeta function,
    we obtain the intersection pairing
    \begin{align*}
        \int_{\Mssfd_{(2,1)}}
        a_2^k \cdot \frac{f_2^{3g-3-2k}}{(3g-3-2k)!}
        & =
        \frac{2^{2g-3-2k}-1}{2^{3g-5-2k} \, \uppi^{2g-2-2k}} \,
        \zeta(2g-2-2k)
        \\ & =
        \frac{1}{2^{g-2} \, \uppi^{2g-2-2k}}
        \sum_{n=1}^{\infty} {} \frac{(-1)^{n+1}}{n^{2g-2-2k}}
        \numberthis
    \end{align*}
    when $k < g-1$,
    as given in \cite[(2.9)]{jeffrey-kirwan-2},
    where
    \begin{align}
        a_2 & = \Xi (-S_{1,0,2}), \\
        f_2 & = \Xi \biggl( -S_{1,2,2} + \sum_{j=1}^g S_{j,1,1} \, S_{j+g,1,1} \biggr).
    \end{align}
    Note that the variables $S_{j,1,1}$ for $j = 1, \dotsc, 2g$
    are irrelevant here,
    as the fixed determinant invariants in this case
    do not involve the variables $s_{j,1,1}$\,.
\end{remark}

\subsection{Comparison with previous results}
\label{sect-comparison}

We compare our results with
Jeffrey--Kirwan~\cite[Theorem~9.12]{jeffrey-kirwan-2},
who gave intersection pairings on $\Mssfd_{(r,d)}$
in the case when $r$ and $d$ are coprime.
We expect our formula to be equivalent to theirs in the coprime case,
and we provide a proof in the simpler case
of the symplectic volume, as in Definition~\ref{def-symp-vol}.

To relate our results with Jeffrey--Kirwan's,
we first need to apply a transformation to our formula
in Theorem~\ref{thm-inv-mss-main-as-reg-sum}.

\begin{lemma}
    \label{lem-inv-mss-inverted}
    For any $(r, d) \in K(X)$ with $r > 0$, we have
    \upshape
    \begin{multline}
        \label{eq-inv-mss-inverted}
        \inv{\Mss_{(r,d)}} = 
        \frac{1}{r} \cdot {} \\
        \sumbar_{ \leftsubstack{
            \\[-2ex]
            & d = d_0 + \cdots + d_{r-1} \\[-.6ex]
            & (d_0 + \cdots + d_{i-1})/i \leq d/r, \ i = 1, \dotsc, r-1 \\[-.6ex]
            & \text{with $m$ equalities} 
        } } \mspace{-9mu}
        \frac{1}{m+1} \cdot
        \bigl[
            \fund{\Mss_{(1,d_0)}} \, , \ 
            \bigl[ \fund{\Mss_{(1,d_1)}} \, , \ 
            \bigl[ \dotsc , \ 
            \fund{\Mss_{(1,d_{r-1})}}
        \bigr] \bigr] \bigr].
        \raisetag{3ex}
    \end{multline}
\end{lemma}

\begin{proof}
    By Theorem~\ref{thm-inv-mss-main-as-reg-sum}
    and~\eqref{eq-reg-sum-degen},
    we see that
    \begin{multline}
        \label{eq-inv-mss-inverted-1}
        \inv{\Mss_{(r,d)}} = 
        \frac{(-1)^{r-1}}{r} \cdot {} \\
        \sumbar_{ \leftsubstack{
            \\[-2ex]
            & d = d_0 + \cdots + d_{r-1} \\[-.6ex]
            & (d_0 + \cdots + d_{i-1})/i \geq d/r, \ i = 1, \dotsc, r-1 \\[-.6ex]
            & \text{with $m$ equalities} 
        } } \mspace{-9mu}
        \frac{1}{m+1} \cdot
        \bigl[ \bigl[ \dotsc \bigl[
            \fund{\Mss_{(1,d_0)}} \, , \ 
            \fund{\Mss_{(1,d_1)}} \bigr] , \ 
            \dotsc \bigr] , \ 
            \fund{\Mss_{(1,d_{r-1})}}
        \bigr]
        \raisetag{3ex}
    \end{multline}
    (note that the inequalities are inverted
    compared to Theorem~\ref{thm-inv-mss-main-as-reg-sum}~\ref{item-inv-mss-main-as-reg-sum}).
    This can be illustrated in the case $r = 3$ by the diagram
    \begin{equation}
        \rankthreeminidiagram{-120}{0} \ = \ 
        \rankthreeminidiagram{60}{180} \ ,
    \end{equation}
    where shaded regions show the domains of the regularized sums,
    and we used~\eqref{eq-reg-sum-degen} to obtain the relations
    \begin{equation}
        \rankthreeminidiagram{-120}{60} \ = \ 
        \rankthreeminidiagram{0}{180} \ = \ 0.
    \end{equation}
    But the inequality $(d_0 + \cdots + d_{i-1})/i \geq d/r$
    is equivalent to $(d_i + \cdots + d_{r-1})/(r-i) \leq d/r$,
    so that if we replace $(d_0, \dotsc, d_{r-1})$
    by $(d_{r-1}, \dotsc, d_0)$ in \eqref{eq-inv-mss-inverted-1},
    we obtain \eqref{eq-inv-mss-inverted}.
    Note that we also used Remark~\ref{rmk-reg-sum-lie-brack}
    to apply the anti-commutativity of Lie brackets inside the regularized sum.
\end{proof}

\begin{lemma}
    \label{lem-chg-of-var}
    In the repeated residue expressions in
    Theorems~\textnormal{\ref{thm-inv-mss-main}, \ref{thm-inv-mss-pairing}}
    and~\textnormal{\ref{thm-inv-mssfd-pairing}},
    the residue operator $\res_{z_{r-1}} \circ \cdots \circ \res_{z_1}$
    can be replaced by
    \[
        (-1)^{r-1} \res_{y_1} \circ \cdots \circ \res_{y_{r-1}} \, ,
    \]
    where $y_i = z_{i-1} - z_i$\,.
\end{lemma}

\begin{proof}
    Each summand of \eqref{eq-main-explicit}
    is a vertex operation
    \begin{equation}
        \label{eq-voa-change-variables-1}
        (-1)^{r-1} \, Y (A_{r-1}, z_{r-1}) \cdots Y (A_0, z_0) \, 1 ,
    \end{equation}
    where we write $A_i = \smash{\fund{\Mss_{(1,d_i)}}}$,
    so that setting $z_0 = 0$ gives
    \begin{equation}
        (-1)^{r-1} \, Y (A_{r-1}, z_{r-1}) \cdots Y (A_1, z_1) \, A_0 \, ,
    \end{equation}
    and taking the repeated residue
    $\res_{z_{r-1}} \circ \cdots \circ \res_{z_1}$
    of this expression
    gives the iterated Lie bracket
    $[ [ \dotsc [ A_0, A_1 ], \dotsc ], A_{r-1}]$.
    
    Using the associativity property of the vertex algebra,
    we may rewrite \eqref{eq-voa-change-variables-1} as
    \begin{equation}
        (-1)^{r-1} \, 
        Y ( Y ( \cdots
            Y ( Y (A_{r-1}, z_{r-1} - z_{r-2}) \, 
            A_{r-2}, z_{r-2} - z_{r-3} ) \cdots A_1,
            z_1 - z_0
        ) \, A_0, z_0 ) \, 1 ,
    \end{equation}
    so that setting $z_0 = 0$ gives
    \begin{equation}
        (-1)^{r-1} \, 
        Y ( \cdots
            Y ( Y (A_{r-1}, z_{r-1} - z_{r-2}) \, 
            A_{r-2}, z_{r-2} - z_{r-3} ) \cdots A_1,
            z_1
        ) \, A_0 \, ,
    \end{equation}
    and applying the repeated residue operator
    $(-1)^{r-1} \res_{y_1} \circ \cdots \circ \res_{y_{r-1}}$ gives
    the iterated Lie bracket
    $[ A_0, [ \dotsc [ A_{r-2}, A_{r-1} ] \dotsc ] ]$.
    
    The result now follows from Lemma~\ref{lem-inv-mss-inverted}.
\end{proof}

We can now obtain an alternative formula for the symplectic volume
of $\Mssfd_{(r,d)}$\,.

\begin{theorem}
    For $(r, d) \in K(X)$ with $r > 0$ and $\gcd (r, d) = 1$,
    the symplectic volume of $\Mssfd_{(r,d)}$ is given by
    \begin{equation}
        \label{eq-symp-vol}
        \int_{\Mssfd_{(r,d)}} \mspace{-6mu}
        \exp (f_2) =
        \res_{y_{1}} \circ \cdots \circ \res_{y_{r-1}}
        \left. \left\{
            \frac{\displaystyle
                (-1)^{(g-1)r(r-1)/2} \, r^{g-1} \cdot
                \exp \biggl( -\sum_{i=0}^{r-1} \tilde{d}_i \, z_i \biggr)
            }
            {\displaystyle
                \prod_{0 \leq i < j \leq r-1} \mspace{-9mu}
                (z_i - z_j)^{2g-2} \cdot
                \prod_{i=1}^{r-1} {} ( \exp y_i - 1 )
            }
        \right\} \right|_{z_0=0},
        \raisetag{-1ex}
    \end{equation}
    where $f_2 = \Xi (-S_{1,2,2})$,
    and $\tilde{d}_i = \lfloor (i+1)d/r \rfloor - \lfloor id/r \rfloor - d/r$.
\end{theorem}

\begin{proof}
    This follows from taking $\alpha = -1$ in \eqref{eq-main-fd-pairing-simpler},
    and applying Lemma~\ref{lem-chg-of-var}.
    The sign $(-1)^{r-1}$ in Lemma~\ref{lem-chg-of-var}
    accounts for changing $1 - \exp y_i$ to $\exp y_i - 1$ on the denominator.
\end{proof}

\begin{remark}
    Jeffrey and Kirwan \cite[Theorem~8.1]{jeffrey-kirwan-2}
    expressed the symplectic volume of $\Mssfd_{(r,d)}$ as 
    \begin{equation}
        \label{eq-symp-vol-jk}
        \res_{y_{1}} \circ \cdots \circ \res_{y_{r-1}}
        \left\{
            \frac{\displaystyle
                (-1)^{(g-1)r(r-1)/2} \cdot
                \sum_{w \in W_{r-1}} \exp \langle \llbr w \tilde{c} \rrbr , \biz \rangle
                \int_{T_r^{2g}} \upe^{\omega}
            }
            {\displaystyle
                r! \cdot
                \prod_{0 \leq i < j \leq r-1} \mspace{-9mu}
                (z_i - z_j)^{2g-2} \cdot
                \prod_{i=1}^{r-1} {} ( \exp y_i - 1 )
            }
        \right\}.
    \end{equation}
    This is equivalent to our formula \eqref{eq-symp-vol}, since we have
    $\llbr w \tilde{c} \rrbr =
    -(\tilde{d}_0, \dotsc, \tilde{d}_{r-1})$
    for any $w \in W_{r-1}$,
    and $|W_{r-1}| = (r-1)!$,
    and we have that $\int_{T_r^{2g}} \upe^{\omega} = r^g$.
    
    As shown in \cite{jeffrey-kirwan-2},
    Witten's \cites{Witten1992} formulae for intersection pairings
    are equivalent to \eqref{eq-symp-vol-jk}
    via the transformation \eqref{eq-szenes-sum-iter-res-2}.
    But we have just shown that \eqref{eq-symp-vol-jk}
    can be transformed to our regularized sum formulae.
    This suggests that the transformation \eqref{eq-szenes-equals-regularized}
    could be the link between Witten's and our formulae.
\end{remark}

\section{Results for elliptic curves}
\label{sect-elliptic-curves}

In this section, let $X$ be an elliptic curve over $\bbC$,
and we compute the invariants
$\smash{\inv{\Mss_{(r,d)}}}$ for~$X$.
Although this was already done in \S\ref{sect-higher-rank},
we use a different method here,
involving a Fourier--Mukai transform on the category $\Db (X)$,
together with results from \S\ref{sect-rank-0}.
This also serves as a verification of the results in \S\ref{sect-higher-rank}.
The main result is Theorem~\ref{thm-ell-main}.

\subsection{The Fourier--Mukai transform}

Let $X$ be an elliptic curve over $\bbC$.
Let $J (X)$ be the Jacobian variety of $X$,
which is also the dual abelian variety of $X$,
and is isomorphic to $X$ via the Abel--Jacobi map.
As in Huybrechts~\cite[\S9]{huybrechts},
we have the \emph{Poincar\'e bundle}
\begin{equation}
    P \longrightarrow X \times J (X),
\end{equation}
which is identical to the universal line bundle
$\calU_{(1,0)}$ defined in \S\ref{sect-rank-1},
with $x_0$ being the unit element of $X$.
One defines the \emph{Fourier--Mukai transform}
with respect to the Poincar\'e bundle by
\begin{align}
    \label{eq-def-fm-transform}
    \Phi \colon \Db (X) & \longrightarrow \Db (J (X)), \\
    \calF & \longmapsto
    (\pr_2)_* \, (\pr_1^* (\calF) \otimes P), \notag
\end{align}
where $\pr_1$ and $\pr_2$ are the projections
from $X \times J (X)$ to its factors,
and we are using the derived versions
of the pushforward, pullback and tensor product functors.
By \cite[Proposition~9.19]{huybrechts},
$\Phi$ is an equivalence of triangulated categories.

Note that since $X$ is isomorphic to $J(X)$,
$\Phi$ can be seen as an auto-equivalence
of the triangulated category $\Db (X)$. 

\begin{lemma}
    \label{lem-fm-acting-on-cohomology}
    Let $\Phi \colon \calM \to \calM$
    be the isomorphism of moduli stacks
    induced by the equivalence \eqref{eq-def-fm-transform}.
    Then $\Phi$ maps $\calM_{(r,d)}$
    isomorphically to $\calM_{(d,-r)}$, and we have
    \begin{align}
        \Phi^* (S_{d,-r;\,1,0,l}) & =
        S_{r,d;\,1,2,l+1}\,, \\
        \Phi^* (S_{d,-r;\,1,1,l}) & =
        S_{r,d;\,2,1,l}\,, \\
        \Phi^* (S_{d,-r;\,2,1,l}) & =
        -S_{r,d;\,1,1,l}\,, \\
        \Phi^* (S_{d,-r;\,1,2,l}) & =
        -S_{r,d;\,1,0,l-1}\,.
    \end{align}
\end{lemma}

\begin{proof}
    For clarity,
    let $\calM'$ denote the moduli stack of objects in
    $\Db ( J(X) )$.
    Let $\epsilon'_{j,k}$ denote the 
    generators of the cohomology ring of $J (X)$,
    induced by $\epsilon_{j,k}$
    via the Abel--Jacobi isomorphism.
    We see that $\epsilon'_{j,1}$ is equal to the element $\gamma_j$
    defined in \eqref{eq-def-gamma-j}.
    Let $S'_{\alpha;\,j,k,l}$ denote the
    generators of the cohomology ring of $\calM'$,
    defined using $\epsilon'_{j,k}$\,.
    
    Since $\Phi$ is an isomorphism of stacks,
    it must map $\calM_{(r,d)}$ isomorphically
    to some $\calM_{(r',d')}$.
    Let $\calU'_{(r',d')} \to J(X) \times \calM'_{(r',d')}$
    be the universal complex.
    Let $\pr_{ij}$ denote the projection from
    $X \times J(X) \times \calM_{(r,d)}$
    to the product of its $i$-th and $j$-th factors.
    Then, by~\eqref{eq-c1-u1d}, we have
    \begin{align*}
        &
        \ch \bigl( (\id_{J(X)} \times \Phi)^* \, (\calU'_{(r',d')}) \bigr) 
        \\ = {} &
        \int_X
        \ch (\pr_{13}^* (\calU_{(r,d)})) \cdot
        \ch (\pr_{12}^* (P))
        \\ = {} &
        \int_X {} \pr_{13}^*
        \biggl(
            \sum_{ \substack{ (j,k) \in J \\ l \geq k/2 } }
            \epsilon_{j,k} \boxtimes S_{r,d;\,j,k,l}
        \biggr) \cdot \pr_{12}^* \bigl(
            1 \boxtimes 1
            -\epsilon_{1,1} \boxtimes \epsilon'_{1,1}
            -\epsilon_{2,1} \boxtimes \epsilon'_{2,1}
            -\epsilon_{1,2} \boxtimes \epsilon'_{1,2}
        \bigr) 
        \\ = {} &
        1 \boxtimes \sum_{l\geq1} S_{r,d;\,1,2,l}
        +\epsilon'_{1,1} \boxtimes \sum_{l\geq1} S_{r,d;\,2,1,l}
        -\epsilon'_{2,1} \boxtimes \sum_{l\geq1} S_{r,d;\,1,1,l}
        -\epsilon'_{1,2} \boxtimes \sum_{l\geq0} S_{r,d;\,1,0,l}\,,
        \numberthis
    \end{align*}
    and the result follows from this computation.
    In particular, we have $\Phi^* (S_{r',d';\,1,0,0}) = d$
    and $\Phi^* (S_{r',d';\,1,2,1}) = -r$,
    so that $r' = d$ and $d' = -r$.
\end{proof}

Consider the Bridgeland~\cite{bridgeland} stability condition $(Z, \calP)$
on $\Db (X)$ given by
$Z (r,d) = -d + \upi r$,
so that the heart $\Coh (X)$
is spanned by $\calP ((0,1])$ under extensions.

For $(r,d) \in K(X)$ and $\phi \in \bbR$ such that
$Z (r,d)$ is a positive multiple of $\upe^{\upi \phi}$,
let
\begin{equation*}
    \calP (\phi)_{(r,d)} \subset \calP (\phi)
\end{equation*}
be the full subcategory spanned by objects
of rank $r$ and degree $d$.

\begin{lemma}
    \label{lem-ss-neg-deg-no-sections}
    Let $\calF \to X$ be a semistable sheaf of negative degree.
    Then $\calF$ does not admit any non-zero global sections.
\end{lemma}

\begin{proof}
    Suppose that $s$ is such a non-zero section.
    Then $s$ generates a subsheaf of $\calF$
    of rank $1$ and degree $0$,
    contradicting the definition of semistability.
\end{proof}

\begin{lemma}
    \label{lem-fm-acts-as-rotation}
    For any $(r,d) \in K(X)$ and $\phi \in \bbR$, one has
    \begin{equation}
        \Phi (\calP (\phi)_{(r,d)}) =
        \calP (\phi - 1/2)_{(d,-r)} \, .
    \end{equation}
\end{lemma}

\begin{proof}
    First, let us show that $\Phi (\calP (\phi)_{(r,d)}) \subset
    \calP (\phi - 1/2)_{(d,-r)}$\,.
    Since $\Phi$ commutes with shifting,
    we may assume that $(r,d) \in C(X)$ and $\phi \in (0, 1]$.
    Also, we may assume that $\calP (\phi)_{(r,d)} \neq \varnothing$.
    
    If $\phi \in (1/2, 1)$, then $r, d > 0$,
    and any sheaf $\calF \in \calP (\phi)_{(r,d)}$ is a vector bundle.
    For any $[L] \in J(X)$,
    by Lemma~\ref{lem-ss-neg-deg-no-sections}, one has
    \begin{equation}
        \label{eq-h0-f-ox-l}
        H^0 (X, (\calF \otimes L)^\vee) = 0.
    \end{equation}
    By Serre duality, we have
    \begin{equation}
        \label{eq-coh-dual}
        \Phi (\calF)^\vee [-1] \simeq
        (\pr_2)_* \, \bigl(
            (\pr_1^* (\calF) \otimes P)^\vee
        \bigr),
    \end{equation}
    where $\pr_1$ and $\pr_2$ are the projections
    from $X \times J (X)$ to its factors,
    and we are using the derived version
    of the pushforward and dual functors.
    By~\eqref{eq-h0-f-ox-l},
    the right-hand side of~\eqref{eq-coh-dual} is concentrated in degree~$1$,
    so that $\Phi (\calF)$ is concentrated in degree $0$.
    That is, it is equivalent to a sheaf, which we denote by $\calF'$.
    By Lemma~\ref{lem-fm-acting-on-cohomology},
    $\calF'$ is a sheaf of rank $d$ and degree $-r$.
    Since $\Phi$ is an equivalence,
    if $\calF$ is indecomposable, then so is $\calF'$,
    so that by \cite[Appendix~A]{tu},
    $\calF'$ is semistable.
    In general, if $\calF$ decomposes
    into the sum of several indecomposable semistable sheaves of the same slope,
    then $\calF'$ will also be the sum of several semistable sheaves of the same slope.
    Therefore, we see that $\Phi (\calP (\phi)_{(r,d)}) \subset
    \calP (\phi - 1/2)_{(d,-r)}$\,.
    
    If $\phi = 1$, then $r = 0$,
    and the above argument still works with some modifications.
    Now $(\pr_1^* (\calF) \otimes P)^\vee$
    is in degree $1$,
    so the right-hand side of~\eqref{eq-coh-dual}
    is in degrees~$[1, 2]$, and $\Phi (\calF)$ in $[-1, 0]$,
    so it has to be in degree $0$.
    
    If $\phi \in (0, 1/2)$, then $d < 0$.
    For any sheaf $\calF \in \calP (\phi)_{(r,d)}$\,,
    and for any $[L] \in J(X)$,
    by Lemma~\ref{lem-ss-neg-deg-no-sections}, one has
    \begin{equation}
        H^0 (X, \calF \otimes L) = 0,
    \end{equation}
    so that $\Phi (\calF)$ is concentrated in degree $1$.
    Suppose that $\Phi (\calF) \simeq \calF'[-1]$
    for some sheaf $\calF'$. Then, similarly to the previous cases,
    we see that $\calF'$ is a semistable sheaf
    of rank $-d$ and degree $r$, so that
    $\Phi (\calP (\phi)_{(r,d)}) \subset
    \calP (\phi - 1/2)_{(d,-r)}$\,.
    
    Finally, if $\phi = 1/2$, then $d = 0$.
    Consider the Jordan--H\"older filtration of $\calF$.
    Then its quotients are $r$ line bundles of degree $0$,
    by \cite[Appendix~A]{tu}.
    This means that if $[L] \in J(X)$
    is not the dual of one of these $r$ line bundles,
    then $H^0 (X, \calF \otimes L) = 0$.
    Thus, $H^0 (X \times U, \pr_1^* (\calF) \otimes P) = 0$
    for any open set $U \subset J (X)$,
    so that $\bbR^0 (\pr_2)_* (\pr_1^* (\calF) \otimes P) = 0$
    and $\Phi (\calF)$ is in degree~$1$,
    so $\Phi (\calF) \in \calP (0)_{(0, -r)}$.
    
    Now that we have shown that
    $\Phi (\calP (\phi)_{(r,d)}) \subset
    \calP (\phi - 1/2)_{(d,-r)}$
    in all cases, we wish to prove the inverse inclusion.
    But this is because by \cite[Proposition~9.19]{huybrechts},
    for any $\calF \in \Db (X)$,
    one has $\Phi^4 (\calF) \simeq \calF [-2]$,
    so that $\Phi^4 (\calP (\phi)_{(r,d)}) =
    \calP (\phi-2)_{(r,d)}$\,.
\end{proof}

\subsection{Invariants for elliptic curves}

In the following, we use the Fourier--Mukai transform
to compute the invariants $\inv{\Mss_{(r,d)}}$ for $X$.

Let $T = [1] \colon \Db ( X ) \to \Db ( X )$
denote the shifting operator,
and also $T \colon \Mpl \to \Mpl$ its induced automorphism of the moduli stack.

\begin{lemma}
    \label{lem-fm-preserves-inv}
    Let $(r,d) \in C(X)$. Then
    \begin{equation}
        \Phi_* \inv{\Mss_{(r,d)}} = \begin{cases}
            \inv{\Mss_{(d,-r)}} \, , & d > 0, \\
            T_* \inv{\Mss_{(-d,r)}} \, , & d \leq 0.
        \end{cases}
    \end{equation}
\end{lemma}

\begin{proof}
    This follows from Lemma~\ref{lem-fm-acts-as-rotation},
    since by the construction of $\inv{\Mss_{(r,d)}}$,
    it is invariant under an automorphism of the triangulated category
    preserving the Bridgeland stability condition.
\end{proof}

\begin{lemma}
    \label{lem-tensor-line-bundle-sjkl}
    Let $F \colon \calM_{(r,d)} \to \calM_{(r,d+r)}$
    be the map induced by tensoring with
    a degree~$1$ line bundle. Then
    \begin{align}
        F_* (\upe^{(r,d)} \cdot s_{1,0,l}) & =
        \upe^{(r,d+r)} \cdot (s_{1,0,l} + s_{1,2,l+1}), \\
        F_* (\upe^{(r,d)} \cdot s_{j,k,l}) & =
        \upe^{(r,d+r)} \cdot s_{j,k,l} \qquad ((j,k) \neq (1,0)).
    \end{align}
\end{lemma}

\begin{proof}
    Let $L$ denote the line bundle.
    One has
    \begin{align*}
        & \ch (\calU_{(r,d)} \otimes \pr_1^* (L))
        \\ = {} &
        \biggl(
            \sum_{ \substack{ (j,k) \in J \\ l \geq k/2 } }
            \epsilon_{j,k} \boxtimes S_{r,d;\,j,k,l}
        \biggr) \cdot \bigl(
            (1 + \epsilon_{1,2}) \boxtimes 1
        \bigr)
        \\ = {} &
        \sum_{ \substack{ (j,k) \in J \\ l \geq k/2 } }
        \epsilon_{j,k} \boxtimes S_{r,d;\,j,k,l}
        +
        \sum_{l \geq 0}
        \epsilon_{1,2} \boxtimes S_{r,d;\,1,0,l}\,.
        \numberthis
    \end{align*}
    Therefore,
    \begin{align}
        F^* (S_{r,d+r;\,1,2,l}) & =
        S_{r,d;1,2,l} + S_{r,d;1,0,l-1}, \\
        F^* (S_{r,d+r;\,j,k,l}) & =
        S_{r,d;j,k,l} \qquad ((j,k) \neq (1,2)),
    \end{align}
    and the result follows.
\end{proof}

\begin{theorem}
    \label{thm-ell-main}
    For any $(r, d) \in C (X)$, we have
    \begin{align*}
        \inv{\Mss_{(r,d)}} & =
        \upe^{(r,d)} \cdot
        (-1)^{(r-1)(d-1)} \Bigl(
            \frac{1}{d} \, s_{1,0,1} +
            s_{1,1,1} \, s_{2,1,1}
        \Bigr) + \im D_{(r,d)} \\
        & =
        \upe^{(r,d)} \cdot
        (-1)^{(r-1)(d-1)} \Bigl(
            -\frac{1}{r} \, s_{1,2,2} +
            s_{1,1,1} \, s_{2,1,1}
        \Bigr) + \im D_{(r,d)} \, ,
        \numberthis
    \end{align*}
    where the first expression applies whenever $d \neq 0$,
    and the second applies whenever $r \neq 0$.
\end{theorem}

\begin{proof}
    First, we see that the two formulae are equivalent when
    both $r$ and $d$ are non-zero, since in this case, one has
    $D_{(r,d)} 1 = r \, s_{1,0,1} + d \, s_{1,2,2}$.
    
    Let us prove the formula by induction on $r$.
    For $r = 0$, the formula is true by Theorem~\ref{thm-inv-mss0d-main}.
    If $r > 0$, then by the induction hypothesis,
    the formula is true for $(r', r) \in C (X)$
    with $r' = 0, \dotsc, r-1$.
    By Lemma~\ref{lem-fm-acting-on-cohomology}
    and Lemma~\ref{lem-fm-preserves-inv},
    the formula is true for $(r, -r')$ for $r'$ as above.
    Finally,
    since tensoring by a line bundle of degree $\pm1$
    identifies $\Mss_{(r,d)}$ with $\Mss_{(r,d\pm r)}$,
    by Lemma~\ref{lem-tensor-line-bundle-sjkl},
    we see that the formula is true for all $(r,d)$.
\end{proof}

\begin{corollary}
    \label{cor-ell-fd}
    For any $(r, d) \in C (X)$, we have
    \begin{equation}
        \inv{\Mssfd_{(r,d)}} =
        \upe^{(r,d)} \cdot
        (-1)^{(r-1)(d-1)}.
    \end{equation}
\end{corollary}

\begin{proof}
    This follows from Theorem~\ref{thm-ell-main}
    and Theorem~\ref{thm-inv-mssfd}.
\end{proof}

Finally, we verify that the results of this section are consistent
with the results in~\S\ref{sect-higher-rank},
which were obtained by a different approach.

\begin{proposition}
    When $r > 0$,
    Theorem~\textnormal{\ref{thm-ell-main}}
    agrees with Theorem~\textnormal{\ref{thm-inv-mss-main}}
    specialized to the case $g = 1$.
\end{proposition}

\begin{proof}
    We apply Theorem~\ref{thm-inv-mss-main}
    to compute the invariant $\inv{\Mss_{(r,d)}}$.
    For degree reasons,
    we may set $s_{1,2,l} = 0$ for $l > 2$.
    This reduces \eqref{eq-main-explicit} to
    \begin{multline*}
        \xi \bigl( \inv{\Mss_{(r,d)}} \bigr) =
        \upe^{(r,d)} \cdot
        \res_{z_{r-1}} \circ \cdots \circ \res_{z_1}
        \Biggl\{
            \frac{(-1)^{(g-1)r(r-1)/2+(r-1)(d-1)}}{r} \cdot {} \\
            \sum_{ \leftsubstack[7em]{
                \\[-2ex]
                & 0 \leq m \leq \gcd(r,d)-1 \\[-.6ex]
                & 1 \leq i_1 < \cdots < i_m \leq r-1 \\[-.6ex]
                & \text{such that } i_k d / r \in \bbZ \text{ for all } k
            } }
            \frac{(-1)^m}{m+1} \cdot 
            \frac{1}{\displaystyle
                \prod_{ \leftsubstack[4em]{
                    & 1 \leq i \leq r-1 \\[-.6ex]
                    & i \neq i_k \text{ for all } k
                } } \bigl( 
                    1 - \exp ( (z_i - z_{i-1}) \, s_{1,2,2} )
                \bigr)
            } \cdot {} \\[-2ex]
            \prod_{i=0}^{r-1} {}
            \bigl(
                \sigma (-s_{1,2,2}) + O (\biz)
            \bigr)
        \Biggr\} \Bigg|_{z_0=0}\,,
        \numberthis
    \end{multline*}
    where $O (\biz)$ is a holomorphic function
    that is zero when $z_0 = \cdots = z_{r-1} = 0$.
    Since $1/(1 - \exp ( (z_i - z_{i-1}) \, s_{1,2,2} ))$
    has a simple pole along the set $z_i = z_{i-1}$,
    we see that the only term that contributes to the $(r-1)$-fold residue
    is the term with $m = 0$, so that
    \begin{align*}
        & \xi \bigl( \inv{\Mss_{(r,d)}} \bigr) \\
        = {} &
        \upe^{(r,d)} \cdot
        \res_{z_{r-1}} \circ \cdots \circ \res_{z_1}
        \biggl(
            \frac{(-1)^{(r-1)(d-1)} \cdot \sigma (-s_{1,2,2})^{r}}
            {r \cdot \prod_{ i=1 }^{ r-1 } \bigl( 
                1 - \exp ( (z_i - z_{i-1}) \, s_{1,2,2} )
            \bigr)} 
        \biggr) \bigg|_{z_0=0} \\
        = {} &
        \upe^{(r,d)} \cdot
        \frac{(-1)^{(r-1)(d-1)} \cdot \sigma(-s_{1,2,2})^r}
        {r \cdot (-s_{1,2,2})^{r-1}} \\
        = {} &
        \upe^{(r,d)} \cdot
        (-1)^{(r-1)(d-1)} \Bigl(
            -\frac{1}{r} \, s_{1,2,2} +
            s_{1,1,1} \, s_{2,1,1}
        \Bigr) ,
        \numberthis
    \end{align*}
    which is consistent with Theorem~\ref{thm-ell-main}.
\end{proof}

\appendix
\section{The regularized sum}
\label{sect-regularized-sum}

This appendix is devoted to introducing
the \emph{regularized sum},
which was used in the statements our main results
computing the higher rank invariants $\inv{\Mss_{(r,d)}}$\,.
We present the definition of the regularized sum
in Definition~\ref{def-regularized-sum} below,
and we discuss a relation to the work of Szenes~\cite{Szenes1998},
which provides a link connecting our formula for intersection pairings
to Witten's formula.

The regularized sum is a way to assign finite values
to divergent geometric series.
It is based on the idea to `incorrectly' extend the formula
\begin{equation}
    \sum_{n=0}^\infty a^n = \frac{1}{1-a}
\end{equation}
to the case $|a| \geq 1$,
so that we may write, for example, $1 + 2 + 4 + 8 + \cdots = -1$.

The regularized sum is denoted by $\sumbar$.
In the present work,
we work primarily with regularized sums of the form
\begin{equation}
    \sumbar_{\bix \in \Lambda} f (\bix; \biz),
\end{equation}
where $\Lambda$ is a lattice in $\bbQ^n$,
or a part of such a lattice,
and for each $\bix \in \Lambda$,
we have a meromorphic function $f (\bix; \biz)$ in 
the complex variables $\biz = (z_1, \dotsc, z_n)$,
depending on $\bix$ in a way described by Definition~\ref{def-regularized-sum} below.
For example, we have
\begin{equation}
    \sumbar_{n \in \bbZ} z^n = 0
\end{equation}
as a meromorphic function in $z$,
since $\sumbar_{n \geq 0} z^n = 1/(1 - z)$,
and $\sumbar_{n < 0} z^n = 1/(z - 1)$.
Here, the lattice being used is $\bbZ \subset \bbQ$.

\subsection{Definition}
\label{sect-def-regularized-sum}

We give the definition of the regularized sum
in Definition~\ref{def-regularized-sum} below.
The main ideas used in the definition
can be seen, perhaps more clearly, in Construction~\ref{cons-s-lambda-c},
since it is a basic case used to define the regularized sum.

Throughout the following, we fix a $\bbQ$-vector space~$V$
of dimension $n \geq 0$.

\begin{definition}
    \label{def-cone-lattice}
    We define several types of subsets of $V$ as follows.
    \begin{enumerate}
        \item
            A subset $C \subset V$ is
            called a \emph{polyhedral cone}, or simply a \emph{cone},
            if there exists an integer $m \geq 0$,
            and linear maps $f_1, \dotsc, f_m \colon V \to \bbQ$, such that
            \begin{equation}
                C = \{ \bix \in V \mid f_i (\bix) \geq 0
                \text{ for all } i \}.
            \end{equation}
            A cone $C \subset V$ is \emph{salient},
            if $C$ does not contain any $1$-dimensional subspace of~$V$,
            as in, for example, \cite[\S2.6.1]{Edwards1965}.

        \item
            A subset $C \subset V$ is
            called a \emph{convex polytope}, or simply a \emph{polytope},
            if there exists an integer $m \geq 0$,
            linear maps $f_1, \dotsc, f_m \colon V \to \bbQ$, 
            and numbers $c_1, \dotsc, c_m \in \bbQ$, such that
            \begin{equation}
                \label{eq-def-polytope}
                C = \{ \bix \in V \mid f_i (\bix) \geq c_i
                \text{ for all } i \}.
            \end{equation}
            Note that such a polytope need not be bounded.

        \item
            A subset $\Lambda \subset V$ is called a \emph{lattice},
            if there exists an integer $m \geq 0$,
            elements $\bix_0, \bie_1, \dotsc, \bie_m \in V$,
            with $\bie_1, \dotsc, \bie_m$ linearly independent, such that
            \begin{equation}
                \label{eq-def-lattice}
                \Lambda = \Bigl\{
                    \bix_0 + \sum_{i=1}^m a_i \, \bie_i \Bigm|
                    a_i \in \bbZ
                \Bigr\}.
            \end{equation}

        \item
            A subset $\Sigma \subset V$ is called a \emph{sector},
            if there exists a polytope $C \subset V$,
            and a lattice $\Lambda \subset V$
            such that $\Sigma = C \cap \Lambda$.

        \item
            A subset $\Sigma \subset V$ is called a
            \emph{simple sector},
            if there exists an integer $m \geq 0$,
            elements $\bix_0, \bie_1, \dotsc, \bie_m \in V$,
            with $\bie_1, \dotsc, \bie_m$ linearly independent, such that
            \begin{equation}
                \label{eq-def-simple-sector}
                \Sigma = \Bigl\{
                    \bix_0 + \sum_{i=1}^m a_i \, \bie_i \Bigm|
                    a_i \in \bbZ_{\geq0}
                \Bigr\}.
            \end{equation}
            In particular, a simple sector is always a sector.
            Note also that singletons are simple sectors.
    \end{enumerate}
\end{definition}

\begin{lemma}
    \label{lem-lattice-simple-sector}
    Every sector in $V$ can be written as
    a disjoint union of finitely many simple sectors.
\end{lemma}

\begin{proof}
    Suppose that $\Sigma = C \cap \Lambda$,
    with $C$ a polytope and $\Lambda$ a lattice.
    We may assume that the affine linear span of $\Sigma$ is $V$,
    since otherwise, we may replace $V$ by this span.
    
    First, we consider the case when $C$ has the form
    \begin{equation}
        \label{eq-simple-chamber}
        C = \Bigl\{
            \sum_{j=1}^n a_j \, \bif_j \Bigm|
            a_j \in \bbQ_{\geq0}
        \Bigr\}
    \end{equation}
    for linearly independent vectors $\bif_1, \dotsc, \bif_n \in V$.
    In this case, for $j = 1, \dotsc, n$,
    let $\tilde{\bif}_j$ be the least positive rational multiple of
    $\smash{\bif_j}$
    that is the difference of two elements of~$\Lambda$.
    Then there are only finitely many elements $\bix \in \Sigma$,
    such that none of the vectors $\bix - \smash{\tilde{\bif}_j}$,
    $j = 1, \dotsc, n$, lie in $\Sigma$.
    Let $\bix_1, \dotsc, \bix_k$ be all these elements, and let
    \begin{equation}
        \Sigma_i = \Bigl\{
            \bix_i + \sum_{j=1}^n a_j \, \tilde{\bif}_j \Bigm|
            a_j \in \bbZ_{\geq0}
        \Bigr\}.
    \end{equation}
    We see that $\Sigma = \Sigma_1 \cup \cdots \cup \Sigma_k$,
    and that the union is disjoint.
    
    Next, we consider the case when $C$ is a salient cone.
    Since every bounded polytope has a linear triangularization
    (we are allowed to add new vertices),
    $C$ can be decomposed into chambers of the form \eqref{eq-simple-chamber},
    which may share parts of their boundaries.
    By the proof of the previous case,
    for each such chamber $C'$,
    the intersection of the interior of $C'$ with $\Lambda$
    can be written as a disjoint union of simple sectors.
    Similarly, the lower dimensional strata can also be decomposed
    into simple sectors.

    The case when $C$ is a non-salient cone is analogous,
    since one can always cut $C$ into salient pieces using hyperplanes.

    Finally, we consider the general case.
    Use the notations in~\eqref{eq-def-polytope},
    and we proceed by induction on the dimension of~$V$.
    We divide into the following cases.

    \begin{enumerate}
        \item 
            \label{itm-pf-simple-decomp-1}
            The space
            $W = \mathrm{span} \bigl( \bigcap_i {} \{ f_i (\bix) \geq 0 \} \bigr) \subset V$
            is not equal to $V$.
            Then there are only finitely many elements $W + \bix \in V/W$,
            where $\bix \in V$, such that $(W + \bix) \cap \Sigma \neq \varnothing$.
            We may then apply the induction hypothesis to each of the slices
            $(W + \bix) \cap \Sigma$.
        \item
            Otherwise, we can choose
            $c'_1, \dotsc, c'_k \in \bbQ$ with $c'_i > c_i$ for all $i$,
            such that $C' = \bigcap_i {} \{ f_i (\bix) \geq c'_i \}$ is 
            a translation of a cone.
            Use the hyperplanes $\{ f_i (\bix) = c'_i \}$ to cut $C$ into pieces,
            with $C'$ being one of them.
            The lower dimensional strata intersected with $\Sigma$
            can be decomposed by the induction hypothesis.
            Each top dimensional stratum other than $C'$ 
            is bounded by a pair of parallel hyperplanes,
            and falls in the case~\ref{itm-pf-simple-decomp-1}.
    \end{enumerate}
\end{proof}

\begin{definition}
    \ 
    \begin{enumerate}
        \item 
            Let $\bbQ [V]$ be the group algebra of $V$.
            As a vector space, we have
            \begin{equation}
                \bbQ [V] \simeq
                \bigoplus_{\bix \in V} \bbQ \cdot \upe^{\bix},
            \end{equation}
            and the multiplication is given by
            $\upe^{\bix} \cdot \upe^{\biy} = \upe^{\bix + \biy}$
            for $\bix, \biy \in V$.
            Let $\bbQ (V)$ denote the fractional field of $\bbQ [V]$.
        \item 
            For a lattice $\Lambda \subset V$, 
            using the notation of \eqref{eq-def-lattice},
            we define linear subspaces
            \begin{alignat}{2}
                \bbQ [\Lambda] & =
                \upe^{\bix_0} \cdot \bbQ [ \upe^{\pm \bie_1}, \dotsc, \upe^{\pm \bie_m} ]
                && \subset \bbQ [V], \\
                \bbQ (\Lambda) & =
                \upe^{\bix_0} \cdot \bbQ ( \upe^{\bie_1}, \dotsc, \upe^{\bie_m} )
                && \subset \bbQ (V).
            \end{alignat}
            Note that these do not depend on the choice of
            $\bix_0$ and $\bie_i$.
            In addition, if $0 \in \Lambda$,
            then $\bbQ [\Lambda]$ is a subalgebra of $\bbQ [V]$,
            and $\bbQ (\Lambda)$ is a subfield of $\bbQ (V)$.
    \end{enumerate}
\end{definition}

The following construction is key to defining regularized sums,
as the elements $S_\Sigma$ defined below
are essentially the \emph{universal} regularized sums.

\begin{construction}
    \label{cons-s-lambda-c}
    Let $\Sigma \subset V$ be a simple sector,
    and let $\Lambda \subset V$ be the minimal lattice containing $\Sigma$.
    Using the notations in~\eqref{eq-def-simple-sector},
    we define an element
    \begin{equation}
        \label{eq-def-s-sigma-simple}
        S_{\Sigma} =
        \frac{\upe^{\bix_0}}{(1 - \upe^{\bie_1}) \dotsm (1 - \upe^{\bie_m})}
        \in \bbQ (\Lambda).
    \end{equation}
    Intuitively, this is what one should get by summing $\upe^{\bix}$
    for all $\bix \in \Sigma$.
    We shall temporarily denote this by
    \begin{equation}
        \label{eq-univ-reg-sum-simple}
        S_{\Sigma} = \sumbar_{\bix \in \Sigma} \upe^{\bix},
    \end{equation}
    and this expression will become valid after we define the regularized sum.
    The element $S_\Sigma$ is well-defined,
    since the vectors $\bix_0, \bie_1, \dotsc, \bie_m$
    are determined by $\Sigma$ up to a permutation of $\bie_1, \dotsc, \bie_m$.
    
    Now, let $\Sigma \subset V$ be an arbitrary sector,
    and let $\Lambda \subset V$ be the minimal lattice containing $\Sigma$.
    We define an element 
    \begin{equation}
        S_{\Sigma} \in \bbQ (\Lambda)
    \end{equation}
    by decomposing $\Sigma$ into disjoint simple sectors,
    $\Sigma = \Sigma_1 \cup \cdots \cup \Sigma_k$,
    as in Lemma~\ref{lem-lattice-simple-sector}, and setting
    \begin{equation}
        \label{eq-s-sigma-decomp}
        S_{\Sigma} = S_{\Sigma_1} + \cdots + S_{\Sigma_k} \, .
    \end{equation}
    Again, intuitively, we should have
    \begin{equation}
        \label{eq-univ-reg-sum-sector}
        S_{\Sigma} = \sumbar_{\bix \in \Sigma} \upe^{\bix},
    \end{equation}
    and this expression will become valid soon.
    We will prove in Lemma~\ref{lem-s-sigma-well-defined} below
    that the element $S_\Sigma$ does not depend on the choice
    of the decomposition~\eqref{eq-s-sigma-decomp}.
\end{construction}

\begin{lemma}
    \label{lem-s-sigma-well-defined}
    Let $\Sigma \subset V$ be a sector.
    Then the element $S_\Sigma$
    in Construction~\textnormal{\ref{cons-s-lambda-c}}
    is well-defined, i.e.~it does not depend on
    the decomposition~\eqref{eq-s-sigma-decomp}.
\end{lemma}

\begin{proof}
    We may assume $\Sigma = C \cap \Lambda$
    with $C$ a polytope and $\Lambda$ a lattice.

    First, we consider the case when $C$ is a salient cone.
    Let $f \colon V \to \bbQ$ be a linear map with
    $f (\bix) < 0$ for all $\bix \in C \setminus \{ 0 \}$.
    Then the sum $\sum_{\bix \in \Sigma} \upe^{f (\bix)}$
    converges in $\bbR$, and we have
    \begin{equation}
        \label{eq-pf-s-sigma-wd}
        \sum_{\bix \in \Sigma} \upe^{f (\bix)} =
        S_\Sigma \, \big|_{ \bix \mapsto f (\bix) \text{ for all } \bix \in V }
    \end{equation}
    by the definition of $S_\Sigma$\,,
    regardless of which decomposition of $\Sigma$ we are using.
    However, the space of such linear functions~$f$
    form a non-empty open subset of $V^\vee$.
    We regard the right-hand side of~\eqref{eq-pf-s-sigma-wd}
    as a rational function with $\bbR$-coefficients
    in the variables
    $x_1 = \upe^{f (\bie_1)}, \dotsc, x_n = \upe^{f (\bie_n)}$,
    with $\bie_1, \dotsc, \bie_n$ as in~\eqref{eq-def-lattice} for~$\Lambda$.
    For those $f$ such that~\eqref{eq-pf-s-sigma-wd} is valid,
    the set of all possible values $(x_1, \dotsc, x_n)$,
    as a subset of $\bbR^n$, is dense in a non-empty open set.
    This means that the right-hand side of~\eqref{eq-pf-s-sigma-wd}
    is uniquely determined
    as a rational function in the variables $x_1, \dotsc, x_n$,
    and $S_\Sigma$ is recovered from this rational function
    by replacing $x_i \mapsto \upe^{\bie_i}$ for all $i$.

    Note that the above argument also proves the following.
    For a lattice $\Lambda$, a salient cone $C$,
    and any subset $\Sigma \subset C \cap \Lambda$
    that can be written as a finite disjoint union of simple sectors,
    one can define the element $S_\Sigma$ similarly,
    and it does not depend on the choice of the decomposition of~$\Sigma$.

    For the general case, any two decompositions of $\Sigma$
    into simple sectors have a common refinement,
    which is a decomposition of $\Sigma$ into disjoint subsets
    $\Sigma_1, \dotsc, \Sigma_k$ of the form
    described in the previous paragraph.
    This follows from Lemma~\ref{lem-lattice-simple-sector},
    and the fact that the intersection of two polytopes is again a polytope.
    It follows that the element $S_\Sigma$ must be the sum of all $S_{\Sigma_i}$,
    no matter which of the two decompositions we choose.
\end{proof}

We mention an important property of the element $S_\Sigma$.

\begin{lemma}
    \label{lem-degen-cone}
    Let $\Sigma \subset V$ be a sector.
    If $\Sigma = \Sigma + \bix$ for some non-zero vector $\bix \in V$, then
    \begin{equation} 
        S_{\Sigma} = 0.
    \end{equation}
\end{lemma}

\begin{proof}
    By the construction of $S_\Sigma$, we have
    \begin{equation}
        S_\Sigma = S_{\Sigma + \bix} = \upe^{\bix} \cdot S_{\Sigma} \, ,
    \end{equation}
    so that $S_{\Sigma} = 0$.
\end{proof}

\begin{definition}
    Let $F$ be a field.
    A \emph{polytope function} on~$V$ with coefficients in $F$
    is a function $c \colon V \to A$ of the form
    \begin{equation}
        \label{eq-def-polytope-function}
        c = \sum_{i=1}^k a_i \, \chi_{C_i}\,,
    \end{equation}
    where $a_i \in F$,
    and $\chi_{C_i}$ is the characteristic function
    of a polytope $C_i$. Let
    \begin{equation}
        \PF (V; F) = \{ c \colon V \to A \mid c \text{ is a polytope function} \}
    \end{equation}
    be the $F$-algebra of polytope functions.

    For a lattice $\Lambda \subset V$, define a linear map
    \begin{equation}
        \label{eq-def-s-lambda}
        S_\Lambda \colon \PF (V; \bbQ) \longrightarrow \bbQ (\Lambda)
    \end{equation}
    by letting $S_\Lambda (\chi_C) = S_{C \cap \Lambda}$ for any polytope $C$,
    where $S_{C \cap \Lambda}$ is the element defined in Construction~\ref{cons-s-lambda-c}.
    Intuitively, for $c \in \PF (V; \bbQ)$, 
    the element $S_\Lambda (c)$ is the universal regularized sum
    \begin{equation}
        \label{eq-univ-reg-sum-pf}
        S_\Lambda (c) =
        \sumbar_{\bix \in \Lambda} c (\bix) \cdot \upe^{\bix},
    \end{equation}
    and this will become valid after the regularized sum is defined.

    Note that the element $S_\Lambda (c)$
    only depends on the restriction $c|_{\Lambda} \colon \Lambda \to \bbQ$.
\end{definition}

\begin{definition}
    \label{def-regularized-sum}
    Suppose that we are given the following data.
    \begin{itemize}
        \item 
            A $\bbQ$-vector space $V$ of finite dimension.
        \item 
            A lattice $\Lambda \subset V$.
            Write $\Lambda = \bix_0 + \Lambda_0$ for some $\bix_0 \in \Lambda$ and $\Lambda_0 \subset V$.
        \item
            A field $F$ of characteristic zero.
        \item
            A function $f \colon \Lambda \to F$,
            such that there exist functions
            $c, g \colon \Lambda \to F$ with $f = c \cdot g$,
            such that 
            \begin{itemize}
                \item
                    $c \colon \Lambda \to F$
                    is a polytope function.
                \item 
                    $g \colon \Lambda \to F^\times$
                    is a \emph{geometric series} in $F$,
                    that is, there exists a group homomorphism
                    $e \colon \Lambda_0 \to F^\times$, such that
                    \begin{equation}
                        \label{eq-geom-series}
                        g (\bix) = e (\bix - \bix_0) \cdot g (\bix_0)
                    \end{equation}
                    for all $\bix \in \Lambda$.
            \end{itemize}
    \end{itemize}
    We will define the \emph{regularized sum}
    of the function $f$ on $\Lambda$, denoted by
    \begin{equation}
        \sumbar_{\bix \in \Lambda} f (\bix)
        \in F \cup \{ \infty \},
    \end{equation}
    by the procedure below.

    Define a $\bbQ$-algebra homomorphism
    \begin{equation}
        e_* \colon \bbQ [\Lambda_0] \longrightarrow F
    \end{equation}
    by sending $\upe^{\bix}$ to $e (\bix)$ for all $\bix \in \Lambda_0$.
    Let $\mathfrak{p}_e = (e_*)^{-1} (0) \subset \bbQ [\Lambda_0]$,
    so that $e_*$ extends to a $\bbQ$-algebra homomorphism
    from the localization $\bbQ [\Lambda_0]_{\mathfrak{p}_e}$
    to $F$. We write this as a map
    \begin{equation}
        \label{eq-h-star-extended}
        e_* \colon \bbQ (\Lambda_0) \to F \cup \{ \infty \},
    \end{equation}
    which takes the value $\infty$
    on elements outside $\bbQ [\Lambda_0]_{\mathfrak{p}_e}$.
    
    We define the \emph{regularized sum} of $f$ to be
    \begin{equation}
        \sumbar_{\bix \in \Lambda} f (\bix) =
        e_* (S_\Lambda (c) \cdot \upe^{-\bix_0}) \cdot g (\bix_0),
    \end{equation}
    where $S_\Lambda$ is the map in~\eqref{eq-def-s-lambda},
    base-changed from $\bbQ$ to $F$.

    To see that this does not depend on the choice of $c, g$,
    let $f = c_1 \cdot g_1 = c_2 \cdot g_2$ be two of such choices,
    and let $e_1, e_2 \colon \Lambda_0 \to F^\times$ be the maps corresponding to $g_1, g_2$\,.
    Let $\Sigma = \Sigma_1 \cup \cdots \cup \Sigma_k$
    be a decomposition into disjoint simple sectors
    such that $c_1$ and $c_2$ are constant on each~$\Sigma_i$\,,
    with values $a_{1, i}, a_{2, i}$\,, respectively.
    The two choices give results
    \begin{equation}
        \sum_{i=1}^k {}
        (e_j)_* (S_{\Sigma_i - \bix_0}) \cdot
        a_{j, i} \, g_j (\bix_0),
    \end{equation}
    for $j = 1, 2$, and this is independent of $j$,
    since for each $i$, we have
    $a_{1, i} \, g_1 (\bix_0) \cdot e_1 (\bix) =
    a_{2, i} \, g_2 (\bix_0) \cdot e_2 (\bix)$ for all $\bix \in \Sigma_i - \bix_0$\,.
\end{definition}

\begin{remark}
    \label{rmk-reg-sum-sector}
    Using the notations in Definition~\ref{def-regularized-sum},
    for a sector $\Sigma = C \cap \Lambda$,
    where $C \subset V$ is a polytope, we often write
    \begin{equation}
        \sumbar_{\bix \in \Sigma} f (\bix)
    \end{equation}
    for the regularized sum
    $\sumbar_{\bix \in \Lambda} \chi_C (\bix) \cdot f (\bix)$,
    where $\chi_C$ is the characteristic function of~$C$.
    This notation was already used in~\eqref{eq-univ-reg-sum-simple}
    and~\eqref{eq-univ-reg-sum-sector}.

    In particular, if $\Sigma$ satisfies $\Sigma = \Sigma + \bix$
    for some non-zero $\bix \in V$, then by Lemma~\ref{lem-degen-cone},
    we always have
    \begin{equation}
        \label{eq-reg-sum-degen}
        \sumbar_{\bix \in \Sigma} f (\bix) = 0.
    \end{equation}
\end{remark}

We prove a version of Fubini's theorem for regularized sums.

\begin{lemma}
    \label{lem-reg-sum-fubini}
    Let $m \geq 1$ be an integer,
    and $F$ a field of characteristic zero.
    For $i = 1, \dotsc, m$,
    let $V_i$ be a finite-dimensional $\bbQ$-vector space,
    $\Lambda_i \subset V_i$ a lattice,
    and $c_i \in \PF (V_i; F)$ a polytope function.
    Write $V = \prod_{i=1}^m V_i$,
    $\Lambda = \prod_{i=1}^m \Lambda_i$,
    and write
    $c (\bix) = \prod_{i=1}^m c_i (\pr_i (\bix))$ for $\bix \in V$.
    Let $g \colon \Lambda \to F$ be a geometric series.
    Then
    \[
        \sumbar_{\bix_1 \in \Lambda_1} {} c_1 (\bix_1) \cdots
        \sumbar_{\bix_m \in \Lambda_m} {}
        c_m (\bix_m) \cdot g (\bix_1, \dotsc, \bix_m) =
        \sumbar_{\bix \in \Lambda}
        c (\bix) \cdot g (\bix) .
    \]
\end{lemma}

\begin{proof}
    It suffices to prove the lemma when each
    $c_i$ is the characteristic function of a polytope $C_i \subset V_i$,
    and it suffices to consider the universal regularized sum with
    $F = \bbQ (V)$ and $g (\bix) = \upe^{\bix}$.
    In this case, the lemma is equivalent to
    \begin{equation}
        \label{eq-univ-reg-sum-fubini}
        \prod_{i=1}^m S_{\Sigma_i} =
        S_{\Sigma} \, ,
    \end{equation}
    where $\Sigma_i = C_i \cap \Lambda_i$ and
    $\Sigma = \prod_{i=1}^m \Sigma_i \subset V$.
    Decomposing each $\Sigma_i$ as in Lemma~\ref{lem-lattice-simple-sector},
    we see that it is enough to prove~\eqref{eq-univ-reg-sum-fubini}
    when each $\Sigma_i$ is a simple sector.
    But this follows from~\eqref{eq-def-s-sigma-simple}.
\end{proof}

Finally, we prove a simple property of regularized sums
that was used in Remark~\ref{rmk-reg-sum-lie-brack}.

\begin{lemma}
    \label{lem-reg-sum-field-aut}
    Using the notations in Definition~\textnormal{\ref{def-regularized-sum},}
    suppose that a regularized sum
    \begin{equation}
        \sumbar_{\bix \in \Lambda} f (\bix)
    \end{equation}
    is defined.
    Let $\sigma \colon F \to F$ be a field automorphism.
    Then we have
    \begin{equation}
        \sumbar_{\bix \in \Lambda} {} \sigma ( f (\bix) ) =
        \sigma \biggl( \sumbar_{\bix \in \Lambda} f (\bix) \biggr).
    \end{equation}
\end{lemma}

\begin{proof}
    By linearity, we may assume that $f = c \cdot g$,
    with $g$ a geometric series and $c = \chi_C$,
    with $C$ a polytope and $C \cap \Lambda$ a simple sector,
    and we use the notations in~\eqref{eq-def-simple-sector}.
    Let $h = g (\bix_0)$, and let $a_i = e (\bie_i)$.
    It remains to prove that
    \begin{equation}
        \sumbar_{n_1, \dotsc, n_m \geq 0}
        \sigma (h a_1^{n_1} \cdots a_m^{n_m})
        = \sigma \biggl( \frac{h}{\prod_{i=1}^m {} (1 - a_i)} \biggr).
    \end{equation}
    But this is true since both sides are equal to
    $\sigma (h) \big/ {\prod_{i=1}^m (1 - \sigma (a_i))}$.
\end{proof}

\subsection{Relation to previous work}
\label{sect-reg-sum-and-szenes}

We explain a relationship between the regularized sum and
a type of convergent infinite sums studied by Szenes \cite{Szenes1998}.
The latter type of infinite sums appear in Witten's formula \cite[(5.21)]{Witten1992}
for intersection pairings on $\Mssfd_{(r,d)}$.
Szenes showed that such an infinite sum can be written as an iterated residue,
and this result was used by Jeffrey--Kirwan \cite{jeffrey-kirwan-2}
to prove that their formulae are equivalent to Witten's.

We will see that this type of iterated residue expressions may be 
rewritten as a regularized sum,
and we believe that via this process,
for intersection pairings on $\Mssfd_{(r,d)}$\,,
Witten's formula transforms to our regularized sum formulae
in \S\ref{sect-higher-rank}.
For the simpler case of the symplectic volume,
this is shown in \S\ref{sect-comparison}.

Let us explain this relationship.
For simplicity, we only consider the case that is involved in Witten's formulae,
as in \cite[\S5.3]{Szenes1998}.

Let $r > 0$ be an integer. Consider the $\bbC$-vector space
\begin{equation}
    V = \{
        \biz = (z_0, \dotsc, z_{r-1}) \in \bbC^r \mid
        z_0 + \cdots + z_{r-1} = 0
    \}.
\end{equation}
Let $\Lambda = V \cap \bbZ^{r}$, which is a lattice in $V$.
For $0 \leq i < j \leq r-1$, let
$H_{i,j} = \{ z_i = z_j \} \subset V$\,, and let
\begin{equation}
    H = \bigcup_{i,j} H_{i,j} \subset V.
\end{equation}

The Weyl group $W_r$ of $\mathrm{SU} (r)$
is isomorphic to the symmetric group on $r$ elements, $\mathfrak{S}_r$\,,
and acts on $V$ by permuting the coordinates $z_i$\,.

Define a new set of coordinates $\biy = (y_1, \dotsc, y_{r-1})$ on $V$
by $y_i = z_{i-1} - z_i$\,.

Let $f \colon V \dashrightarrow \bbC$
be a rational function,
whose poles can only lie in $H$.
Assume that $f$ is invariant under the action of the Weyl group
$W_r \simeq \mathfrak{S}_r$\,.
As in \cite[(5.3)]{Szenes1998} and \cite[Theorem~2.4]{jeffrey-kirwan-2},
for all $\bit \in V^\vee$
that are real (i.e.\ real linear combinations of the $y_i$), we have
\begin{equation}
    \label{eq-szenes-sum-iter-res}
    \sum_{ \biy \in \Lambda \setminus H }
    f (2 \uppi \upi \biy) \exp \bigl( 2 \uppi \upi \langle \bit, \biy \rangle \bigr)
    =
    \sum_{w \in W_{r-1}} {}
    \res_{y_1} \circ \cdots \circ \res_{y_{r-1}}
    \frac{
        f (\biy)
        \exp \langle \llbr w \bit \rrbr, \biy \rangle
    }
    {\prod_{i=1}^{r-1} {} (1 - \exp y_i)}
\end{equation}
whenever the left-hand side
is absolutely convergent,
where $W_{r-1} \subset W_r$ is the subgroup
of permutations of the first $(r-1)$ coordinates
$z_0, \dotsc, z_{r-2}$\,, isomorphic to $\mathfrak{S}_{r-1}$\,.
The notation $\llbr w \bit \rrbr$ denotes the unique element of $w \bit + \Lambda$
that can be expressed as $\sum_i a_i y_i$ with $0 \leq a_i < 1$ for all $i$.

Our applications will only involve some special choices of $\bit$.
For these $\bit$, the element $\llbr w \bit \rrbr$
does not depend on $w$ for any $w \in W_r$\,,
and hence $\exp \bigl( 2 \uppi \upi \langle \bit, \biy \rangle \bigr)$
is Weyl invariant for $\biy \in \Lambda$.
In this case, \eqref{eq-szenes-sum-iter-res} becomes
\begin{equation}
    \label{eq-szenes-sum-iter-res-2}
    \sum_{ \biy \in \Lambda^+ }
    f (2 \uppi \upi \biy) \exp \bigl( 2 \uppi \upi \langle \bit, \biy \rangle \bigr)
    =
    \frac{1}{r}
    \res_{y_1} \circ \cdots \circ \res_{y_{r-1}}
    \frac{
        f (\biy)
        \exp \langle \llbr \bit \rrbr, \biy \rangle
    }
    {\prod_{i=1}^{r-1} {} (1 - \exp y_i)} \, ,
\end{equation}
where $\Lambda^+ = \{ \biy \in \Lambda \mid y_i > 0 \text{ for all } i \}$.
This is a result of averaging over the Weyl group action.

Let $\Lambda' \subset V^\vee$ be the dual lattice of $\Lambda$,
i.e., $\Lambda' = \bigoplus_i \bbZ y_i$\,.
Define
\begin{equation}
    \textstyle \Lambda'_{\bit} = (\Lambda' + \bit) \cap
    \bigl\{ \sum_{i=1}^{r-1} a_i y_i \bigm| a_i \in \bbQ_{\geq0} \bigr\}
    \subset V^\vee.
\end{equation}
Then \eqref{eq-szenes-sum-iter-res-2} can be rewritten as a regularized sum
\begin{equation}
    \label{eq-szenes-equals-regularized}
    \sum_{ \biy \in \Lambda^+ }
    f (2 \uppi \upi \biy) \exp \bigl( 2 \uppi \upi \langle \bit, \biy \rangle \bigr)
    =
    \frac{1}{r}
    \res_{y_1} \circ \cdots \circ \res_{y_{r-1}}
    \sumbar_{\bia \in \Lambda'_{\bit}} f (\biy) \exp \langle \bia, \biy \rangle,
\end{equation}
by \eqref{eq-def-s-sigma-simple} and the definition of the regularized sum.

As Witten's formula for intersection pairings
is expressed as a sum over $\Lambda^+$ of this type,
and our result in \S\ref{sect-higher-rank} is a regularized sum
which can be transformed to look like the right-hand side of
\eqref{eq-szenes-equals-regularized},
as carried out in \S\ref{sect-comparison},
we expect \eqref{eq-szenes-equals-regularized}
to be the bridge that connects Witten's results with ours.

\section{Proof of Theorems~\ref{thm-inv-mss-main-as-reg-sum} and~\ref{thm-inv-mss-main}}
\label{sect-proof-of-main-theorem}

\subsection{An inductive formula for higher rank invariants}

We follow Joyce's~\cite[\S8.6]{Joyce2021} inductive algorithm
and express the invariants
$\smash{\inv{\Mss_{(r,d)}}}$ and
$\smash{\inv{\Mbss_{(r,d),1}}}$
in terms of the invariants of lower ranks,
so one can compute these invariants by induction on rank.
The main result will be Theorem~\ref{thm-reduction-main}.

\begin{definition}
    For $(r, d) \in K (X)$ with $r > 0$ and $r \nu + d > 0$, 
    where $\nu$ is as in Notation~\ref{ntn-n},
    define
    \begin{align}
        \label{eq-def-tau}
        \tau \colon \check{H}_{\bullet} (\Mbpl_{(r,d),1})
        & \longrightarrow \check{H}_{\bullet} (\Mpl_{(r,d)}) , \\
        \gamma & \longmapsto
        \frac{\upe^{(r,d)}}{\upe^{((r,d),1)}} \cdot \biggl(
            \Bigl( -\frac{\partial}{\partial s_{+,0,1}} \Bigr)^{r\nu+d-1} \,
            \acute{\xi} (\gamma)
        \biggr) \, \bigg|_{s_{+,0,l} = 0 \text{ for all } l > 0} \, . \notag
    \end{align}
    One of its important properties is that
    \begin{equation}
        \tau ( -[ \upe^{((0,0),1)}, \ \gamma ] ) = \gamma
    \end{equation}
    for any $\gamma \in \check{H}_{\bullet} (\Mpl_{(r,d)})$.
\end{definition}

The next lemma gives some combinatorial identities
that will be used in the proof of the main theorem.

\begin{lemma}
    \label{lem-comb-of-z-acute-is-zero}
    Let $m > 0$ be an integer,
    and let $d_0, \dotsc, d_m, r_0, \dotsc, r_m$
    be integers such that $r_i > 0$ for all $i$.
    \begin{enumerate}
        \item
        One has
        \upshape
        \begin{equation}
            \label{eq-sum-of-vee-coeff-equals-zero}
            \sum_{ \leftsubstack[8em]{
                \\[-2ex]
                & 0 \leq m' \leq m \\[-1.6ex]
                & \text{such that } \textstyle
                d_1 / r_1 \geq \cdots \geq d_{m'} / r_{m'} \geq
                (\sum_{i=0}^{m'} d_i) / (\sum_{i=0}^{m'} r_i),
                \\[-1.6ex]
                & \text{and } \textstyle
                (\sum_{i=0}^{m'} d_i) / (\sum_{i=0}^{m'} r_i)
                < d_{m'+1} / r_{m'+1} \leq \cdots \leq d_m / r_m \, .
                \\[-1.6ex]
                & \text{Define } 
                0 = a_0 < \cdots < a_l = m' 
                = b_0 < \cdots < b_{l'} = m \\[-1.6ex]
                & \text{such that for any } 0 < i < m', \ 
                d_i / r_i > d_{i+1} / r_{i+1} \text{ if and only if }
                i = a_j \text{ for some } 0 < j < l ,
                \\[-1.6ex]
                & \text{and for any } m' < i < m, \ 
                d_i / r_i < d_{i+1} / r_{i+1} \text{ if and only if }
                i = b_j \text{ for some } 0 < j < l'
            } }
            \frac{(-1)^{m'}}{\prod_{i=1}^l {} (a_i - a_{i-1})! \cdot \prod_{i=1}^{l'} {} (b_i - b_{i-1})!}
            = 0.
        \end{equation}
        \itshape
        
        \item
        One has
        \upshape
        \begin{equation}
            \label{eq-sum-of-semi-vee-coeff-equals-zero}
            \sum_{ \leftsubstack[8em]{
                \\[-2ex]
                & 0 \leq m' \leq m \\[-1.6ex]
                & \text{such that } \textstyle
                d_1 / r_1 \geq \cdots \geq d_{m'} / r_{m'} \geq
                (\sum_{i=0}^{m'} d_i) / (\sum_{i=0}^{m'} r_i),
                \\[-1.6ex]
                & \text{and } \textstyle
                (\sum_{i=0}^{m'} d_i) /
                (\sum_{i=0}^{m'} r_i)
                = d_{m'+1} / r_{m'+1} = \cdots = d_m / r_m \, .
                \\[-1.6ex]
                & \text{Define } 
                0 = a_0 < \cdots < a_l = m'
                \text{ such that for any } 0 < i < m'
                \\[-1.6ex]
                & d_i / r_i > d_{i+1} / r_{i+1} \text{ if and only if }
                i = a_j \text{ for some } 0 < j < l.
                \\[-1.6ex]
                & \text{Let $e$ be the number of $0 < i < m'$ with } \textstyle
                d_i / r_i = (\sum_{i=0}^{m} d_i) / (\sum_{i=0}^{m} r_i)
            } }
            \frac{(-1)^{m'+e} \, B_e}{\prod_{i=1}^l {} (a_i - a_{i-1})! \cdot (m-m'+1)!}
            = 0,
        \end{equation}
        \itshape
        where $B_e$ denotes the $e$-th Bernoulli number,
        with the convention that $B_1 = -1/2$.
    \end{enumerate}
\end{lemma}

\begin{proof}
    For (i),
    since the sum requires that $d_1 / r_1 \geq \cdots \geq d_{m'} / r_{m'}$
    and $d_{m'+1} / r_{m'+1} \leq \cdots \leq d_m / r_m$\,,
    if the sum is not empty, one must have
    \begin{equation}
        d_1 / r_1 \geq \cdots \geq d_{m_1} / r_{m_1}
        > d_{m_1+1} / r_{m_1+1} = \cdots = d_{m_2} / r_{m_2} <
        d_{m_2+1} / r_{m_2+1} \leq \cdots \leq d_m / r_m
    \end{equation}
    for some integers $0 \leq m_1 < m_2 \leq m$.
    In this case, we only need to prove that
    \begin{equation}
        \sum_{m'=m_1}^{m_2}
        \frac{(-1)^{m'}}{(m'-m_1)! \ (m_2-m')!} = 0,
    \end{equation}
    which is an elementary fact.
    
    For (ii), similarly,
    if the sum is not empty, one must have
    \begin{equation}
        d_1 / r_1 \geq \cdots \geq d_{m_1} / r_{m_1}
        > d_{m_1+1} / r_{m_1+1} = \cdots = d_m / r_m
    \end{equation}
    for some $0 \leq m_1 \leq m$.
    Thus, noticing that in each term of the sum,
    one has $e = m' - m_1$, we only need to prove that
    \begin{equation}
        \sum_{m'=m_1}^{m}
        \frac{B_{m'-m_1}}{(m'-m_1)! \ (m-m'+1)!} = 0,
    \end{equation}
    which is an identity of the Bernoulli numbers that holds whenever $m > 0$.
\end{proof}

\begin{theorem}
    \label{thm-reduction-main}
    Let $(r,d) \in K(X)$ with $r \geq 1$. Then
    \upshape
    \begin{multline}
        \label{eq-inv-mbss-equals-lower-rank}
        \inv{ \Mbss_{(r,d),1} } = \hspace{-2em}
        \sum_{ \leftsubstack{
            \\
            & (r,d) = (1,d_0) + (r_1,d_1) + \cdots + (r_m,d_m), \\[-.6ex]
            & m \geq 1, \ r_i > 0 \text{ for all } i, \\[-.6ex]
            & \text{such that }
            d_1 / r_1 \geq \cdots \geq d_m / r_m > d/r .
            \\[-.6ex]
            & \text{Define } 
            0 = a_0 < \cdots < a_l = m
            \text{ such that for any } 0 < i < m,
            \\[-.6ex]
            & d_i / r_i > d_{i+1} / r_{i+1} \text{ if and only if }
            i = a_j \text{ for some } 0 < j < l
        } }
        \frac{1}{\prod_{i=1}^l (a_i - a_{i-1})!} \cdot
        \bigl[ \bigl[ \dotsc \bigl[ \bigl[
            \inv{\Mbss_{(1,d_0),1}} \, , \ 
            \inv{\Mss_{(r_1,d_1)}} \bigr] , \\[-11.5ex]
            \shoveright{\inv{\Mss_{(r_2,d_2)}} \bigr] , \dotsc \bigr] , \ } \\
            \mathmakebox[\dwidth][r]{ \inv{\Mss_{(r_m,d_m)}}
        \bigr], }
        \xmathstrut[1]{0}
    \end{multline}
    \itshape
    and
    \upshape
    \begin{multline}
        \label{eq-inv-mss-equals-lower-rank}
        \inv{ \Mss_{(r,d)} } = \hspace{-2em}
        \sum_{ \leftsubstack[6em]{
                \\
                & (r,d) = (1,d_0) + (r_1,d_1) + \cdots + (r_m,d_m), \\[-.6ex]
                & m \geq 1, \ r_i > 0 \text{ for all } i, \\[-.6ex]
                & \text{such that }
                d_1 / r_1 \geq \cdots \geq d_m / r_m \geq d/r .
                \\[-.6ex]
                & \text{Define } 
                0 = a_0 < \cdots < a_l = m
                \text{ such that for any } 0 < i < m,
                \\[-.6ex]
                & d_i / r_i > d_{i+1} / r_{i+1} \text{ if and only if }
                i = a_j \text{ for some } 0 < j < l.
                \\[-.6ex]
                & \text{Let $e$ be the number of $i > 0$ with }
                d_i / r_i = d/r
        } }
        \mspace{-24mu}
        \frac{(-1)^e \, B_e}{\prod_{i=1}^l (a_i - a_{i-1})!} \cdot
        \tau \bigl( \bigl[ \bigl[ \dotsc \bigl[ \bigl[
            \inv{\Mbss_{(1,d_0),1}} \, , \ 
            \inv{\Mss_{(r_1,d_1)}} \bigr] , \\[-13.5ex] 
            \shoveright{\inv{\Mss_{(r_2,d_2)}} \bigr] , \dotsc \bigr] , \ } \\
            \mathmakebox[\dwidth][r]{ \inv{\Mss_{(r_m,d_m)}}
        \bigr] \bigr), }
        \xmathstrut[2]{0}
    \end{multline}
    \itshape
    where $B_e$ denotes the $e$-th Bernoulli number,
    with the convention that $B_1 = -1/2$,
    and $\tau$ is the map defined by \eqref{eq-def-tau}.
    In both expressions,
    only finitely many terms of the sum are non-zero.
\end{theorem}

\begin{proof}
    First, let us prove \eqref{eq-inv-mbss-equals-lower-rank}
    by induction on $r$.
    The case $r = 1$ is trivially true.
    If $r > 1$, then by Theorem~\ref{thm-wcf-1}, we have
    \begin{multline}
        \label{eq-inv-mbss-equals-brackets-of-inv}
        \inv{\Mbss_{(r,d),1}} = \hspace{-1.5em}
        \sum_{ \leftsubstack[5em]{
                \\
                & (r,d)=(r_0,d_0)+\cdots+(r_m,d_m), \\[-.6ex]
                & m \geq 1, \ r_i > 0 \text{ for all } i, \\[-.6ex]
                & \text{such that }
                d_0 / r_0 < d_1 / r_1 \leq \cdots \leq d_m / r_m \, .
                \\[-.6ex]
                & \text{Define } 
                0 = a_0 < \cdots < a_l = m
                \text{ such that for any } 0 < i < m,
                \\[-.6ex]
                & d_i / r_i < d_{i+1} / r_{i+1} \text{ if and only if }
                i = a_j \text{ for some } 0 < j < l,
        } }
        \frac{(-1)^{m-1}}{\prod_{i=1}^l (a_i - a_{i-1})!} \cdot
        \bigl[ \bigl[ \dotsc \bigl[ \bigl[
            \inv{\Mbss_{(r_0,d_0),1}} \, , \ 
            \inv{\Mss_{(r_1,d_1)}} \bigr] , \\[-11.5ex] 
            \shoveright{\inv{\Mss_{(r_2,d_2)}} \bigr] , \dotsc \bigr] , \ } \\
            \mathmakebox[\dwidth][r]{ \inv{\Mss_{(r_m,d_m)}}
        \bigr]. }
        \xmathstrut[1]{0}
    \end{multline}
    For each term in the sum
    \eqref{eq-inv-mbss-equals-brackets-of-inv},
    by the inductive hypothesis,
    we expand $\smash{\inv{\Mbss_{(r_0,d_0),1}}}$
    using \eqref{eq-inv-mbss-equals-lower-rank}
    and obtain
    \begin{multline}
        \label{eq-inv-mbss-equals-vee-sum}
        \inv{\Mbss_{(r,d),1}} = \hspace{-3.5em}
        \sum_{ \leftsubstack[10em]{
                \\[1.5ex]
                & (r,d)=(1,d_0)+(r_1,d_1)+\cdots+(r_m,d_m), \\[-.6ex]
                & 0 \leq m' < m, \ r_i > 0 \text{ for all } i, \\[-.6ex]
                & \text{such that } \textstyle
                d_1 / r_1 \geq \cdots \geq d_{m'} / r_{m'}
                \geq (\sum_{i=0}^{m'} d_i) / (1 + \sum_{i=1}^{m'} r_i),
                \\[-.6ex]
                & \text{and } \textstyle
                (\sum_{i=0}^{m'} d_i) / (1 + \sum_{i=1}^{m'} r_i)
                < d_{m'+1} / r_{m'+1} \leq \cdots \leq d_m / r_m \, .
                \\[-.6ex]
                & \text{Define } 
                0 = a_0 < \cdots < a_l = m' 
                = b_0 < \cdots < b_{l'} = m
                \\[-.6ex]
                & \text{such that for any } 0 < i < m', \ 
                d_i / r_i > d_{i+1} / r_{i+1} \text{ if and only if }
                i = a_j \text{ for some } 0 < j < l ,
                \\[-.6ex]
                & \text{and for any } m' < i < m, \ 
                d_i / r_i < d_{i+1} / r_{i+1} \text{ if and only if }
                i = b_j \text{ for some } 0 < j < l'
        } } \mspace{-30mu}
        \frac{(-1)^{m-m'-1}}{
            \prod_{i=1}^l (a_i - a_{i-1})! \cdot \prod_{i=1}^{l'} (b_i - b_{i-1})!
        } \cdot {} \\[-16ex] 
        \shoveright{ \bigl[ \bigl[ \dotsc \bigl[ \bigl[
            \inv{\Mbss_{(1,d_0),1}} \, , \ 
            \inv{\Mss_{(r_1,d_1)}} \bigr] , \ } \\[-.5ex]
            \shoveright{\inv{\Mss_{(r_2,d_2)}} \bigr] , \dotsc \bigr] , \ } \\[-.5ex]
            \mathmakebox[\dwidth][r]{ \inv{\Mss_{(r_m,d_m)}}
        \bigr]. }
        \xmathstrut[1.5]{0}
    \end{multline}
    By Lemma~\ref{lem-comb-of-z-acute-is-zero}~(i),
    \eqref{eq-inv-mbss-equals-vee-sum} is equivalent to
    \eqref{eq-inv-mbss-equals-lower-rank} for $(r, d)$,
    since the summation condition of
    \eqref{eq-inv-mbss-equals-vee-sum}
    is the same as that of
    \eqref{eq-sum-of-vee-coeff-equals-zero}
    with the case $m' = m$ excluded,
    where we set $r_0 = 1$.
    This completes the proof of \eqref{eq-inv-mbss-equals-lower-rank}.
    
    Next, we aim to prove \eqref{eq-inv-mss-equals-lower-rank}
    by showing the stronger result that
    \begin{multline}
        \label{eq-inv-mss-equals-lower-rank-strong}
        -\bigl[ \upe^{((0,0),1)}, \ \inv{\Mss_{(r,d)}} \bigr] = \\
        \sum_{ \leftsubstack[6em]{
                \\
                & (r,d) = (1,d_0) + (r_1,d_1) + \cdots + (r_m,d_m), \\[-.6ex]
                & m \geq 1, \ r_i > 0 \text{ for all } i, \\[-.6ex]
                & \text{such that }
                d_1 / r_1 \geq \cdots \geq d_m / r_m \geq d/r .
                \\[-.6ex]
                & \text{Define } 
                0 = a_0 < \cdots < a_l = m
                \text{ such that for any } 0 < i < m,
                \\[-.6ex]
                & d_i / r_i > d_{i+1} / r_{i+1} \text{ if and only if }
                i = a_j \text{ for some } 0 < j < l.
                \\[-.6ex]
                & \text{Let $e$ be the number of $i > 0$ with }
                d_i / r_i = d/r
        } }
        \mspace{-9mu}
        \frac{(-1)^e \, B_e}{\prod_{i=1}^l (a_i - a_{i-1})!} \cdot
        \bigl[ \bigl[ \dotsc \bigl[ \bigl[
            \inv{\Mbss_{(1,d_0),1}} \, , \ 
            \inv{\Mss_{(r_1,d_1)}} \bigr] , \\[-13.5ex] 
            \shoveright{\inv{\Mss_{(r_2,d_2)}} \bigr] , \dotsc \bigr] , \ } \\
            \mathmakebox[\dwidth][r]{ \inv{\Mss_{(r_m,d_m)}}
        \bigr]. }
        \xmathstrut[2]{0}
    \end{multline}
    By Theorem~\ref{thm-wcf-2}, we have
    \begin{multline}
        \label{eq-bracket-e01-inv-mss-expansion}
        \bigl[ \upe^{((0,0),1)}, \ \inv{\Mss_{(r,d)}} \bigr] =
        -\inv{\Mbss_{(r,d),1}} + {} \\
        \sum_{ \substack{
            (r,d) = (r_1,d_1) + \cdots + (r_m,d_m), \\
            m \geq 2, \ r_i > 0, \ 
            d_i / r_i = d / r, \ \forall i
        } }
        \frac {(-1)^m} {m!} 
        \bigl[ \bigl[ \dotsc \bigl[ \bigl[
            \upe^{((0,0),1)}, \ 
            \inv{\Mss_{(r_1,d_1)}} \bigr], \ 
            \inv{\Mss_{(r_2,d_2)}} \bigr], \\[-3.5ex]
            \dotsc \bigr],
            \inv{\Mss_{(r_m,d_m)}}
        \bigr] .
    \end{multline}
    For each term in the sum
    \eqref{eq-bracket-e01-inv-mss-expansion},
    by the inductive hypothesis,
    we may expand $\bigl[ \upe^{((0,0),1)}, \ \allowbreak
    \smash{\inv{\Mbss_{(r_0,d_0),1}}} \bigr]$
    using \eqref{eq-inv-mss-equals-lower-rank-strong},
    so we obtain
    \begin{multline}
        \label{eq-inv-mss-equals-semi-vee-sum}
        -\bigl[ \upe^{((0,0),1)}, \ \inv{\Mss_{(r,d)}} \bigr] =
        \inv{\Mbss_{(r,d),1}} + {} \\
        \shoveleft{ \sum_{ \leftsubstack[6em]{
                \\
                & (r,d) = (1,d_0) + (r_1,d_1) + \cdots + (r_m,d_m), \\[-.6ex]
                & 0 \leq m' < m, \ r_i > 0 \text{ for all } i, \\[-.6ex]
                & \text{such that }
                d_1 / r_1 \geq \cdots \geq d_{m'} / r_{m'}
                = \cdots = d_m / r_m = d/r.
                \\[-.6ex]
                & \text{Define } 
                0 = a_0 < \cdots < a_l = m'
                \text{ such that for any } 0 < i < m',
                \\[-.6ex]
                & d_i / r_i > d_{i+1} / r_{i+1} \text{ if and only if }
                i = a_j \text{ for some } 0 < j < l.
                \\[-.6ex]
                & \text{Let $e$ be the number of $0 < i \leq m'$ with }
                d_i / r_i = d/r
        } }
        \frac{(-1)^{m-m'+e} \, B_e}{\prod_{i=1}^l (a_i - a_{i-1})! \cdot (m-m'+1)!} \cdot {} } \\[-12.5ex]
        \shoveright{ \bigl[ \bigl[ \dotsc \bigl[ \bigl[
            \inv{\Mbss_{(1,d_0),1}} \, , \ 
            \inv{\Mss_{(r_1,d_1)}} \bigr] , \ } \\
            \shoveright{\inv{\Mss_{(r_2,d_2)}} \bigr] , \dotsc \bigr] , \ } \\
            \mathmakebox[\dwidth][r]{ \inv{\Mss_{(r_m,d_m)}}
        \bigr]. }
        \raisetag{-1ex}
    \end{multline}
    By \eqref{eq-inv-mbss-equals-lower-rank},
    the term $\smash{\inv{\Mbss_{(r,d),1}}}$
    on the right-hand side of \eqref{eq-inv-mss-equals-semi-vee-sum}
    gives all the terms in the sum \eqref{eq-inv-mss-equals-lower-rank-strong}
    with $e = 0$.
    By Lemma~\ref{lem-comb-of-z-acute-is-zero}~(ii),
    the sum in \eqref{eq-inv-mss-equals-semi-vee-sum} equals
    the sum of all terms in \eqref{eq-inv-mss-equals-lower-rank-strong}
    with $e > 0$,
    since the summation condition of
    \eqref{eq-inv-mss-equals-semi-vee-sum}
    is the same as that of
    \eqref{eq-sum-of-semi-vee-coeff-equals-zero}
    with the case $m' = m$ excluded,
    where we set $r_0 = 1$.
    This completes the proof of \eqref{eq-inv-mss-equals-lower-rank-strong}.
\end{proof}

\subsection{Explicit formulae for Lie brackets}

\begin{definition}
    For $\bid \in \bbZ^r$, write
    \begin{align}
        \label{eq-def-z-d}
        Z (\bid) & =
        \bigl[ \bigl[ \dotsc \bigl[
            \inv{\Mss_{(1,d_0)}} \, , \ 
            \inv{\Mss_{(1,d_1)}} \bigr] , \dotsc \bigr] , \ 
            \inv{\Mss_{(1,d_{r-1})}}
        \bigr], \\
        \label{eq-def-z-acute-d}
        \acute{Z} (\bid) & =
        \bigl[ \bigl[ \dotsc \bigl[
            \inv{\Mbss_{(1,d_0),1}} \, , \ 
            \inv{\Mss_{(1,d_1)}} \bigr] , \dotsc \bigr] , \ 
            \inv{\Mss_{(1,d_{r-1})}}
        \bigr].
    \end{align}
\end{definition}

\begin{lemma}
    \label{lem-exp-dds}
    For any $l > 0$, we have
    \begin{equation}
        \exp \Bigl( a \, \frac{\partial}{\partial s_{1,0,l}} \Bigr) \,
        \bigl[
            \exp ( z \, D_{(r, d)} ) \, \sigma (-s_{1,2,2})
        \bigr] =
        \exp \Bigl(
            a r \, \frac{z^l}{l!}
        \Bigr) \cdot
        \bigl[
            \exp ( z \, D_{(r, d)} ) \, \sigma (-s_{1,2,2})
        \bigr].
    \end{equation}
\end{lemma}

\begin{proof}
    This follows from applying Lemma~\ref{lem-exp-z-d-expansion}
    to convert $D_{(r,d)}$ into $D_{(0,d)}$,
    then applying $\exp (a \, \partial / \partial s_{1,0,l})$,
    and converting back.
\end{proof}

\begin{theorem}
    \label{thm-z-d-explicit}
    We have
    \begin{align}
        \label{eq-z-d-explicit}
        Z (\bid) & =
        \upe^{(r,d)} \cdot
        \res_{z_{r-1}} \circ \cdots \circ \res_{z_1}
        \bigl( Z (\bid; \biz) |_{z_0=0} \bigr) + \im D_{(r,d)} \, , \\
        \label{eq-z-acute-d-explicit}
        \acute{Z} (\bid) & =
        \upe^{((r,d),1)} \cdot
        \res_{z_{r-1}} \circ \cdots \circ \res_{z_1} \circ \res_w 
        \bigl( \acute{Z} (\bid; w, \biz) |_{z_0=0} \bigr) + \im \acute{D}_{(r,d),1} \, ,
    \end{align}
    where
    \begin{align}
        \label{eq-def-z-d-z}
        Z (\bid; \biz) & =
        \frac{(-1)^{(g-1)r(r-1)/2+(r-1)(d-1)}}{\displaystyle
            \prod_{0 \leq i < j \leq r-1} \mspace{-9mu}
            (z_i - z_j)^{2g-2}
        } \cdot 
        \prod_{i=0}^{r-1} {}
        \bigl[
            \exp ( z_i \, D_{(1, d_i)} ) \, \sigma (-s_{1,2,2})
        \bigr], \\
        \acute{Z} (\bid; w, \biz) & =
        \frac{(-1)^{(g-1)r(r-1)/2+(r-1)(d-1)}}
        {\displaystyle
            \prod_{0 \leq i < j \leq r-1} \mspace{-9mu}
            (z_i - z_j)^{2g-2} \cdot
            \prod_{i=0}^{r-1} {} (w + z_i)^{\nu+d_i}
        } \cdot {} \notag \\
        \label{eq-def-z-acute-d-z}
        & \hspace{6em}
        \prod_{i=0}^{r-1} {}
        \biggl[
            \exp ( z_i \, D_{(1, d_i)} ) \,
            \sigma \Bigl( \frac{1}{w + z_i} - s_{1,2,2} \Bigr)
        \biggr] \cdot
        \rho (w),
    \end{align}
    where $\sigma$ and $\rho$ are defined by
    \eqref{eq-def-sigma} and \eqref{eq-def-rho}, respectively.
\end{theorem}

\begin{proof}
    We prove \eqref{eq-z-d-explicit}
    by induction on $r$.
    The case $r = 1$ is true by Theorem~\ref{thm-fund-mss-1-d}.
    
    For $r > 1$, write $\bid' = (d_0, \dotsc, d_{r-2})$
    and $d' = d_0 + \cdots + d_{r-2}$, so that 
    \begin{equation}
        Z (\bid) = \bigl[ Z (\bid'), \, \inv{\Mss_{(1,d_{r-1})}} \bigr].
    \end{equation}
    Write
    \begin{equation}
        A (z) =
        \sum_{l=1}^\infty {}
        (2g-2) \, (l-1)! \, z^{-l}
        \frac{\partial}{\partial s'_{1,0,l}} \, ,
    \end{equation}
    so that the vertex operation on sheaves
    given by \eqref{eq-vertex-algebra-sheaves}
    can be written as
    \begin{multline}
        \label{eq-y-explicit}
        Y \bigl( (\upe^{(1,d_{r-1})} \cdot \sigma (-s_{1,2,2})), \ z_{r-1} \bigr) \,
        Z (\bid') =
        \upe^{(r,d)} \cdot
        \frac{(-1)^{(g-1)(r-1) + d' + (r-1) d_{r-1}}}
        {z^{(2g-2)(r-1)}} \cdot {} \\
        \exp (z_{r-1} \, D_{(1, d_{r-1})})
        \exp (A (z_{r-1}))
        \Bigl( \sigma (-s_{1,2,2}) \cdot \frac{Z' (\bid')}{\upe^{(r-1,d')}} \Bigr) \,
        \Big|_{s'_{\smash{j,k,l}} = s_{j,k,l}} \ ,
    \end{multline}
    where $Z' (\bid')$ is the expression obtained from $Z (\bid')$
    by replacing all variables $s_{j,k,l}$ by primed variables $s'_{j,k,l}$\,.
    
    By Lemma~\ref{lem-exp-dds},
    the operator $\exp (A (z_{r-1}))$ in \eqref{eq-y-explicit}
    produces a factor of
    \begin{align*}
        & \phantom{{} = {}}
        {\exp} \biggl(
            \sum_{l=1}^{\infty} {}
            (2g-2) \,
            \frac{1}{l}
            \sum_{i=1}^{r-2} {}
            \frac{
                z_{i}^{l}
            } {z_{r-1}^{l}}
        \biggr) \\
        & = \exp \biggl(
            (2g-2)
            \sum_{i=1}^{r-2} {}
            \log
            \frac {z_{r-1}} {z_{r-1} - z_{i}}
        \biggr) \\
        & =
        \prod_{i=1}^{r-2} {}
        \Bigl(
            \frac{z_{r-1}}{z_{r-1} - z_i}
        \Bigr)^{2g-2} \\
        & =
        \frac{z_{r-1}^{(2g-2)(r-2)}}
        {\prod_{i=1}^{r-2} {} (z_{r-1} - z_i)^{2g-2}} \, .
        \numberthis
    \end{align*}
    Noticing that
    \begin{multline}
        \frac{(-1)^{(g-1)r(r-1)/2 + (r-1)(d-1)}}
        {\prod_{0 \leq i < j \leq r-1} {} (z_i - z_j)^{2g-2}} \Big|_{z_0 = 0} =
        -
        \frac{(-1)^{(g-1)(r-1)(r-2)/2 + (r-2)(d'-1)}}
        {\prod_{0 \leq i < j \leq r-2} {} (z_i - z_j)^{2g-2}} \Big|_{z_0 = 0} \cdot {} \\
        \frac{(-1)^{(g-1)(r-1) + d' + (r-1) d_{r-1}}}
        {z_{r-1}^{(2g-2)(r-1)}} \cdot 
        \frac{z_{r-1}^{(2g-2)(r-2)}}
        {\prod_{i=1}^{r-2} {} (z_{r-1} - z_i)^{2g-2}} \, ,
    \end{multline}
    we have proved \eqref{eq-z-d-explicit}.
    
    For \eqref{eq-z-acute-d-explicit},
    We also use induction on $r$.
    The case $r = 1$ is true by Theorem~\ref{thm-fund-mbss-1-d-1}.
    
    Write $\bid' = (d_0, \dotsc, d_{r-2})$
    and $d' = d_0 + \cdots + d_{r-2}$. Let
    \begin{multline}
        \label{eq-def-a-acute}
        \acute{A} (z) =
        \sum_{l=1}^\infty {}
        (2g-2) \, (l-1)! \, z^{-l}
        \frac{\partial}{\partial s'_{1,0,l}}
        -\sum_{l=0}^\infty {}
        l! \, z^{-l-1}
        \frac{\partial}{\partial s_{1,2,2}}
        \frac{\partial}{\partial r'_{l}} \\ {} +
        \sum_{l=1}^\infty {}
        (\nu+d_{r-1}) \, (l-1)! \, z^{-l}
        \frac{\partial}{\partial r'_{l}} \, ,
    \end{multline}
    with the understanding that
    $\partial / \partial r'_0 = 1$,
    so that \eqref{eq-vertex-algebra-pairs} becomes, in our case,
    \begin{multline}
        \label{eq-y-acute-explicit}
        \acute{Y} (\upe^{((1,d_{r-1}),0)} \cdot \sigma(-s_{1,2,2}), \ z_{r-1}) \,
        \acute{Z} (\bid') =
        \upe^{((r,d), 1)} \cdot
        \frac{(-1)^{(g-1)(r-1) + (r-1)d_{r-1} + d'}}
        {z_{r-1}^{(2g-2)(r-1)+\nu+d_{r-1}}} \cdot {} \\
        \exp (z_{r-1} \, D_{(1, d_{r-1})})
        \exp (\acute{A} (z_{r-1}))
        \Bigl( \sigma (-s_{1,2,2}) \cdot \frac{\acute{Z}' (\bid')}{\upe^{((r-1,d'),1)}} \Bigr) \,
        \Big|_{\acute{s}_{j,k,l} = \acute{s}'_{j,k,l}} \ .
    \end{multline}
    
    Similarly to the first part of the proof,
    let us compute the effects of
    the operator $\exp (\acute{A} (z_{r-1}))$.
    The first term of \eqref{eq-def-a-acute} is the same as $A (z)$,
    and produces a factor of
    \begin{align}
        \frac{z_{r-1}^{(2g-2)(r-2)}}
        {\prod_{i=1}^{r-2} {} (z_{r-1} - z_i)^{2g-2}} \, .
    \end{align}
    The second term turns $\sigma (-s_{1,2,2})$ into
    \begin{equation}
        \sigma \biggl(-s_{1,2,2}
            + \sum_{l=0}^{\infty} {}
            \frac{ (-w)^{l} } {z_{r-1}^{l+1}}
        \biggr) =
        \sigma \Bigl(
            \frac {1} { w + z_{r-1} } - s_{1,2,2}
        \Bigr).
    \end{equation}
    The third term produces a factor of
    \begin{align*}
        & \phantom{{} = {}}
        {\exp} \biggl(
            \sum_{l=1}^{\infty} {}
            (\nu+d_{r-1}) \,
            \frac{1}{l}
            \frac{(-w)^{l}} {z_{r-1}^{l}}
        \biggr) \\
        & = \exp \biggl(
            (\nu+d_{r-1})
            \log \frac {z_{r-1}} {w + z_{r-1}}
        \biggr) \\
        & =
        \Bigl(
            \frac {z_{r-1}} {w + z_{r-1}}
        \Bigr)^{\nu+d_{r-1}}.
        \numberthis
    \end{align*}
    One can then deduce \eqref{eq-z-acute-d-explicit}
    from this computation.
\end{proof}

Finally, we establish that
the elements $Z (\bid; \biz)$ and
$\acute{Z} (\bid; w, \biz)$ in \eqref{eq-def-z-d-z}--\eqref{eq-def-z-acute-d-z}
are geometric series with respect to $\bid$.

\begin{lemma}
    \label{lem-exp-z-d-expansion}
    For any integers $r, r', d, d'$, we have
    \begin{equation}
        \label{eq-lem-exp-zd-expansion}
        \exp \bigl( z D_{(r+r',d+d')} \bigr) =
        \exp \biggl(
            \sum_{l=1}^{\infty}
            \frac{z^l}{l!} 
            \, ( r' s_{1,0,l} + d' s_{1,2,l+1} )
        \biggr)
        \exp \bigl( z D_{(r,d)} \bigr),
    \end{equation}
    where $z$ is a formal variable of degree $-2,$
    and $D_{(r, d)}$ is as in~\eqref{eq-def-d-r-d}.
\end{lemma}

\begin{proof}
    Write
    \begin{equation}
        \Delta = D_{(0,0)} = \sum_{ \substack {
            (j, k) \in J \\
            l > k/2
        } }
        s_{j,k,l+1} \, \frac{\partial}{\partial s_{j,k,l}} \, ,
    \end{equation}
    and write $c_{(r',d')} = r' s_{1,0,1} + d' s_{1,2,2}$\,.
    
    We use the Zassenhaus formula,
    as in~\cite[Theorem~1.1]{zassenhaus},
    to expand the operator
    $\exp \bigl( z D_{(r+r',d+d')} \bigr) =
    \exp \bigl( z D_{(r,d)} + z c_{(r',d')} \bigr)$
    as an infinite product
    \begin{equation}
        \exp \bigl( z D_{(r+r',d+d')} \bigr) \circ
        \exp \bigl( z c_{(r',d')} \bigr) \circ
        \exp (C_2) \circ \exp (C_3) \circ \cdots ,
    \end{equation}
    where each $C_m$ is a linear combination of
    $(m-1)$-fold Lie brackets of
    $z D_{(r+r',d+d')}$ and $z c_{(r',d')}$.
    However, the only non-zero terms in these linear combinations are
    \begin{equation}
        \frac{(-1)^{m-1}}{m!}
        \operatorname{ad}(z D_{(1,d)})^{m-1} \, (z c_{(r',d')}) =
        \frac{(-1)^{m-1}}{m!} \,
        z^m \, \Delta^{m-1} \, c_{(r',d')} \, ,
    \end{equation}
    so that
    \begin{align*}
        \exp \bigl( z D_{(r+r',d+d')} \bigr)
        & =
        \exp \bigl( z D_{(r,d)} \bigr) \circ
        \prod_{m=1}^{\infty}
        \exp \biggl(
            \frac{(-1)^{m-1}}{m!} \,
            z^m \, \Delta^{m-1} \, c_{(r',d')}
        \biggr) \\
        & =
        \exp \bigl( z D_{(r,d)} \bigr) \circ
        \exp \biggl(
            \frac{1 - \upe^{-z \Delta}}{\Delta} \, c_{(r',d')}
        \biggr) \\
        & =
        \exp \biggl(
            \exp \bigl( \operatorname{ad} (z D_{(r,d)}) \bigr)
            \Bigl(
                \frac{1 - \upe^{-z \Delta}}{\Delta} \, c_{(r',d')}
            \Bigr)
        \biggr) \circ
        \exp \bigl( z D_{(r,d)} \bigr) \\
        & =
        \exp \Bigl( \frac{\upe^{z \Delta} - 1}{\Delta} \, c_{(r',d')} \Bigr)
        \exp \bigl( z D_{(r,d)} \bigr),
        \numberthis
    \end{align*}
    which is equivalent to~\eqref{eq-lem-exp-zd-expansion}.
\end{proof}

\begin{lemma}
    \label{lem-z-multiplicative}
    The maps
    \begin{align}
        \label{eq-z-multiplicative}
        Z (-; -) \colon \bbZ^r & \longrightarrow F_r^+, \\
        \label{eq-z-acute-multiplicative}
        \acute{Z} (-; -) \colon \bbZ^r & \longrightarrow \acute{F}_r^+,
    \end{align}
    defined by \eqref{eq-def-z-d-z}--\eqref{eq-def-z-acute-d-z},
    where $F_r^+$ and $\acute{F}_r^+$ are as in
    Definition~\textnormal{\ref{def-field-inf-many-variables}},
    are geometric series in the following sense:

    The image of~\eqref{eq-z-multiplicative} lies in
    a one-dimensional $F_r$-subspace of $F_r^+$,
    and~\eqref{eq-z-multiplicative} is a geometric series
    in this subspace,
    in the sense of Definition~\textnormal{\ref{def-regularized-sum}}.
    Similar statements are true for~\eqref{eq-z-acute-multiplicative},
    with $\acute{F}_r, \acute{F}_r^+$ in place of $F_r, F_r^+$.
\end{lemma}

\begin{proof}
    This follows from Lemma~\ref{lem-exp-z-d-expansion}
    and the explicit expressions~\eqref{eq-def-z-d-z}--\eqref{eq-def-z-acute-d-z}.
\end{proof}

\subsection{The coefficient function}

We prove a complicated combinatorial result, 
Theorem~\ref{thm-c-delta-main},
which is a key ingredient in the proof of our main theorems.
One of its special cases
was already demonstrated as the calculation~\eqref{eq-pict-comb-main}
in the sketch proof of the main result.

\begin{definition}
    \label{def-vector-space-vr}
    Let $r > 0$ be an integer. We define the $\bbQ$-vector space
    \begin{equation}
        V_r = \{
            \bix = (x_0, \dotsc, x_{r-1}) \in \bbQ^r \mid
            x_0 + \cdots + x_{r-1} = 0
        \}.
    \end{equation}
\end{definition}

\begin{definition}
    \label{def-c-delta}
    Let $I = \{ I (i) \colon V_r \to \bbQ \mid i = 1, \dotsc, k \}$
    be a linearly independent set of linear functions on $V_r$.
    We define an element $c_I \in \PF (V_r; \bbQ)$ as follows.
    For $\bix \in V_r$, if the system of inequalities
    \begin{equation}
        \label{eq-def-positive-region}
        I (i) (\bix) \geq 0, \quad i = 1, \dotsc, k
    \end{equation}
    holds, then we define 
    \begin{equation}
        c_I (\bix) = \frac{1}{m+1},
    \end{equation}
    where $m$ is the number of $i = 1, \dotsc, k$ with
    $I (i) (\bix) = 0$.
    If \eqref{eq-def-positive-region} does not hold,
    then define $c_I (\bix) = 0$.
    
    In other words, $c_I$ takes the value $1$
    in the bulk of the region defined by
    \eqref{eq-def-positive-region},
    and takes $1/(m+1)$ on its codimension $m$ boundary.
    
    The element $c_{\Delta} \in \PF (V_r; \bbQ)$ will be
    important, where
    \begin{equation}
        \Delta = \{ x_i + \cdots + x_{r-1} \mid i = 1, \dotsc, r-1 \}.
    \end{equation}
\end{definition}

\begin{definition}
    For an integer $r > 0$,
    let $P_r$ be the set of permutations
    $\alpha$ of the set $\{ 0, \dotsc, r-1 \}$,
    such that if $\alpha (i) = 0$, then
    \begin{equation}
        \alpha (0) > \cdots > \alpha (i) < \cdots < \alpha (r-1).
    \end{equation}
    Note that $|P_r| = 2^{r-1}$,
    since there is a bijection $P_r \to 2^{\{ 1, \dotsc, r-1 \}}$
    given by $\alpha \mapsto \{ i \mid \alpha^{-1} (i) < \alpha^{-1} (0) \}$.
\end{definition}

\begin{lemma}
    \label{lem-iter-brackets-sum}
    Let $\calL$ be a Lie algebra.
    Then for any elements $a_0, \dotsc, a_{r-1}, b \in \calL$, where $r > 0$,
    we have the identity
    \begin{multline}
        [[[ \dotsc [[a_0, a_1], a_2], \dotsc], a_{r-1}], b] = \\
        \sum_{\alpha \in P_r} {} (-1)^{\alpha^{-1} (0)+1} \,
        [[ \dotsc [[ b, a_{\alpha(0)} ], a_{\alpha(1)} ], \dotsc], a_{\alpha(r-1)}].
    \end{multline}
\end{lemma}

\begin{proof}
    We prove the identity by induction on $r$.
    If $r = 1$, the identity says that $[a_0, b] = -[b, a_0]$.
    If $r > 1$, applying the Jacobi identity, one has
    \begin{align*}
        & [[[[ \dotsc [[a_0, a_1], \dotsc], a_{r-2}], a_{r-1}], b] \\[1ex]
        = {} &
        [[[[ \dotsc [a_0, a_1], \dotsc], a_{r-2}], b], a_{r-1}] -
        [[[ \dotsc [a_0, a_1], \dotsc], a_{r-2}], [b, a_{r-1}]] \\[1ex]
        = {} &
        \sum_{\alpha \in P_{r-1}} {} (-1)^{\alpha^{-1} (0) + 1} \,
        [[[ \dotsc [[ b, a_{\alpha(0)} ], a_{\alpha(1)} ], \dotsc], a_{\alpha(r-2)}], a_{r-1}] + {} \\*[-2ex]
        & \hspace{6em} \sum_{\alpha \in P_{r-1}} {} (-1)^{\alpha^{-1} (0)} \,
        [[ \dotsc [[[ b, a_{r-1} ], a_{\alpha(0)} ], a_{\alpha(1)} ], \dotsc], a_{\alpha(r-2)}] \\
        = {} &
        \sum_{\alpha \in P_r} {} (-1)^{\alpha^{-1} (0) + 1} \,
        [[ \dotsc [[ b, a_{\alpha(0)} ], a_{\alpha(1)} ], \dotsc], a_{\alpha(r-1)}].
        \numberthis
    \end{align*}
\end{proof}

\begin{definition}
    Let $\PF^0 (V_r; \bbQ) \subset \PF (V_r; \bbQ)$
    be the subspace spanned by characteristic functions
    $\chi_{C}$ of polytopes $C$
    that contains a $1$-dimensional affine linear subspace of $V_r$.
    Define a quotient space
    \begin{equation}
        \label{eq-def-pftilde}
        \PFtilde (V_r; \bbQ) = \frac{\PF (V_r; \bbQ)}{\PF^0 (V_r; \bbQ)}.
    \end{equation}
    By Lemma~\ref{lem-degen-cone},
    it is valid to use elements of $\PFtilde (V_r; \bbQ)$
    as coefficients in a regularized sum,
    since using elements of $\PF^0 (V_r; \bbQ)$
    always gives zero.
\end{definition}

\begin{theorem}
    \label{thm-c-delta-main}
    In $\PFtilde (V_r; \bbQ)$, one has
    \begin{equation}
        \sum_{\alpha \in P_r} {} (-1)^{\alpha^{-1} (0)} \,
        (\alpha^{-1})^* (c_{\Delta}) = r \cdot c_\Delta \, ,
    \end{equation}
    where $(\alpha^{-1})^* (c_\Delta)$ is defined by
    \begin{equation}
        (\alpha^{-1})^* (c_\Delta) (x_0, \dotsc, x_{r-1}) =
        c_\Delta (x_{\alpha^{-1} (0)}, \dotsc, x_{\alpha^{-1} (r-1)}).
    \end{equation}
\end{theorem}

\begin{proof}
    For each $\alpha \in P_r$\,,
    let $I_\alpha = (\alpha^{-1})^* (\Delta) =
    \{ I_\alpha (i) \mid i = 1, \dotsc, r - 1 \}$, where
    \begin{equation}
        I_\alpha (i) = \sum_{j \geq i} x_{\alpha^{-1} (j)} \, ,
    \end{equation}
    be the set of inequalities that defines $(\alpha^{-1})^* (c_\Delta)$.
    
    Let
    \begin{equation}
        J_\alpha = \{ i = 1, \dotsc, r-1 \mid \alpha^{-1} (i) < \alpha^{-1} (0) \}.
    \end{equation}
    Let $\calI_\alpha$ be the set of subsets $I \subset I_\alpha$
    that contain all $I_\alpha (i)$ with
    $i \notin J_\alpha$. Define
    \begin{equation}
        c_{\alpha} = \sum_{I \in \calI_\alpha} {}
        (-1)^{|I_\alpha \setminus I|} \, c_I \,.
    \end{equation}
    Since $c_I = 0$ in $\PFtilde (V_r; \bbQ)$ whenever $|I| < r-1$,
    we see that 
    \begin{equation}
        \label{eq-alpha-star-c-delta-equals-c-alpha}
        (\alpha^{-1})^* (c_\Delta) = c_\alpha \quad \text{in} \quad \PFtilde (V_r; \bbQ) \, .
    \end{equation}
    We regard $c_\alpha$ as a modification of $(\alpha^{-1})^* (c_\Delta)$
    with its support moved to a different region.
    Indeed, by the inclusion--exclusion principle,
    we see that $c_\alpha$ is supported in the region
    defined by the set of inequalities
    $I'_\alpha = \{ I'_\alpha (i) \mid i = 1, \dotsc, r-1 \}$,
    where $I'_\alpha (i)$ is defined by
    \begin{equation}
        I'_\alpha (i) =
        (-1)^{[i \in J_\alpha]} \, I_\alpha (i) \, ,
    \end{equation}
    where $[i \in J_\alpha]$ is $1$ if $i \in J_\alpha$, or $0$ otherwise.
    
    However, $c_\alpha$ is different from $c_{I'_\alpha}$\,,
    since the values on the boundary of the region can be different.
    Indeed, for integers $m \geq m' \geq 0$, write
    \begin{equation}
        \upL^{m+1}_{m'+1} = \frac{1}{(m+1) \, \binom{m}{m'}}
    \end{equation}
    for the entries of the Leibniz harmonic triangle.
    By an elementary argument,
    one shows that if $\bix \in V_r$ is in the support of $c_\alpha$, then
    \begin{equation}
        \label{eq-c-alpha-x-equals-leibniz}
        c_\alpha (\bix) = \upL^{m+1}_{m'+1},
    \end{equation}
    where $m$ is the number of $i$ such that $I'_\alpha (i) (\bix) = 0$,
    and $m'$ is the number of $i \in J_\alpha$
    such that $I'_\alpha (i) (\bix) = 0$.
    
    We claim that $\supp c_\alpha \subset \supp c_\Delta$.
    Indeed, suppose that $\bix \in \supp c_\alpha$.
    For $i = 1, \dotsc, r-1$, we wish to show that
    the $i$-th inequality of $\Delta$ is satisfied,
    i.e., $\smash{\sum_{j=i}^{r-1} x_j} \geq 0$.
    Let 
    \begin{align}
        A_i & = \{ j = 1, \dotsc, r-1 \mid \alpha^{-1} (j-1) < i \text{ and } \alpha^{-1} (j) \geq i \} , \\
        B_i & = \{ j = 1, \dotsc, r-1 \mid \alpha^{-1} (j-1) \geq i \text{ and } \alpha^{-1} (j) < i \} .
    \end{align}
    Then $A_i \cap J_\alpha = \varnothing$ and $B_i \subset J_\alpha$.
    It follows that
    \begin{equation}
        \sum_{j=i}^{r-1} x_j =
        \sum_{j \in A_i} I_\alpha (j) (\bix) -
        \sum_{j \in B_i} I_\alpha (j) (\bix) =
        \sum_{j \in A_i \cup B_i} I'_\alpha (j) (\bix) \geq 0,
    \end{equation}
    and this proves the claim.
    
    Now, we start the main part of the proof.
    We prove the stronger result that for each $i = 0, \dotsc, r-1$, one has
    \begin{equation}
        \label{eq-main-comb-lemma-strong}
        \sum_{ \substack{ \alpha \in P_r \\ \alpha^{-1} (0) = i } } c_\alpha =
        c_\Delta.
    \end{equation}
    By \eqref{eq-alpha-star-c-delta-equals-c-alpha},
    this will imply the desired result.
    
    Fix $i \in \{ 0, \dotsc, r-1 \}$, and fix $\bix \in \supp c_\Delta$.
    Let $Q_i \subset P_r$
    be the set of $\alpha \in P_r$ with
    $\alpha^{-1} (0) = i$ and $\bix \in \supp c_\alpha$.
    Let $0 = a_0 < a_1 < \cdots < a_p = r$ such that
    $\{ a_0, \dotsc, a_p \}$ is the set of all $j = 0, \dotsc, r$
    satisfying $\sum_{j'=j}^{r-1} x_{j'} = 0$.
    
    For each $\beta \in P_p$, let $Q_{i,\beta} \subset Q_i$
    be the set of all $\alpha \in Q_i$ satisfying the following property:
    Write $a'_j = a_{\beta^{-1} (j)+1} - a_{\beta^{-1} (j)}$
    for $j=0,\dotsc,p-1$. Then
    \begin{equation}
        \label{eq-def-q-i-beta}
        \textstyle
        \alpha \bigl( \{ a_{j}, \dotsc, a_{j+1}-1 \} \bigr) =
        \bigl\{
            \sum_{j'=0}^{\beta(j)-1} a'_{\beta^{-1} (j')}, \dotsc,
            \sum_{j'=0}^{\beta(j)} a'_{\beta^{-1} (j')} - 1
        \bigr\}
    \end{equation}
    for all $j=0,\dotsc,p-1$.
    
    Let $j_0 \in \{ 0, \dotsc, p-1 \}$
    be the element such that $a_{j_0} \leq i < a_{j_0+1}$\,.
    We claim that
    \begin{equation}
        Q_i = \bigcup_{ \substack{ \beta \in P_p \\ \beta^{-1} (0) = j_0 } } Q_{i,\beta} \, .
    \end{equation}
    To prove the claim, let us first try to see that
    for all $\alpha \in Q_i$, one has
    $\alpha (\{ a_{j_0}, \dotsc, a_{j_0+1} - 1 \}) =
    \{ 0, \dotsc, a'_{\smash{j_0}}-1 \}$.
    Suppose the contrary, and let
    $h \in \{ 0, \dotsc, a'_{\smash{j_0}}-1 \}$ be the least element
    that is not in $\alpha (\{ a_{j_0}, \dotsc, a_{j_0+1} - 1 \})$.
    Then either (i)~$\alpha^{-1} (h) = a_{j_0} - 1$\,,
    or (ii)~$\alpha^{-1} (h) = a_{j_0+1}$.
    For case (i),
    one must have $\alpha^{-1} (\{ 0, \dotsc, h-1 \}) =
    \{ a_{j_0}, \dotsc, a_{j_0} + h - 1 \}$,
    so that the inequality $I'_\alpha (h)$ says that
    $\sum_{h' \notin \{ a_{j_0}, \dotsc, a_{j_0} + h - 1 \}} x_{h'} \leq 0$.
    Note that $h \in J_\alpha$ in this case.
    But by the definition of the numbers $a_j$, we have 
    $\sum_{h' < a_{j_0}} x_{h'} = 0$.
    This means that
    $\sum_{h' \geq a_{j_0}+h} x_{h'} \leq 0$.
    The equality cannot be attained, since
    $a_{j_0}+h$ is not one of the numbers $a_j$\,.
    Therefore, $\sum_{h' \geq a_{j_0}+h} x_{h'} < 0$,
    contradicting with the assumption that $\bix \in \supp c_\Delta$.
    For case (ii), similarly,
    one must have $\alpha^{-1} (\{ 0, \dotsc, h-1 \}) =
    \{ a_{j_0+1} - h, \dotsc, a_{j_0} - 1 \}$,
    so that the inequality $I'_\alpha (h)$ says that
    $\sum_{h' \notin \{ a_{j_0+1} - h, \dotsc, a_{j_0+1} - 1 \}} x_{h'} \geq 0$,
    where we note that $h \notin J_\alpha$ in this case.
    But $\sum_{h' \geq a_{j_0+1}} x_{h'} = 0$, so
    $\sum_{h' < a_{j_0+1} - h} x_{h'} \geq 0$, and hence
    $\sum_{h' \geq a_{j_0+1} - h} x_{h'} \leq 0$.
    Similarly, the equality cannot be attained, and we obtain a contradiction.
    Thus, we have shown that
    $\alpha (\{ a_{j_0}, \dotsc, a_{j_0+1} - 1 \}) =
    \{ 0, \dotsc, a'_{\smash{j_0}}-1 \}$.
    A similar argument shows that
    $\alpha$ must belong to $Q_{i,\beta}$ for some $\beta$,
    and this proves the claim.
    
    Next, we wish to show that for each $\beta \in P_p$ with
    $\beta^{-1} (0) = j_0$, one has
    \begin{equation}
        \label{eq-sum-over-q-i-beta-equals-leibniz}
        \sum_{\alpha \in Q_{i,\beta}}
        c_\alpha (\bix) =
        \upL^p_{j_0+1}.
    \end{equation}
    For each $\alpha \in Q_{i,\beta}$ and $j = 1, 2, \dotsc, r$,
    let $\alpha_j$ denote the sequence
    $\langle \alpha^{-1} (0), \dotsc, \alpha^{-1} (j-1) \rangle$.
    For each sequence $\gamma = \langle \gamma(0), \dotsc, \gamma(j-1) \rangle$
    of $j$ integers,
    let $Q_{i,\beta} (\gamma)$ be the set of $\alpha \in Q_{i,\beta}$
    with $\alpha_j = \gamma$, and denote
    \begin{equation}
        l (\gamma) = \sum_{\alpha \in Q_{i,\beta} (\gamma)} c_\alpha (\bix).
    \end{equation}
    For a sequence $\gamma$ of $j$ integers such that
    $Q_{i,\beta} (\gamma) \neq \varnothing$,
    and for $j' = 1,\dotsc,j$, let $\gamma_{(j')}$
    denote the sequence $\langle \gamma(0), \dotsc, \gamma(j'-1) \rangle$.
    Let $R_{i,\beta} (\gamma)$ be the set of sequences $\gamma'$ of $(j+1)$ integers
    such that $\gamma'_{(j)} = \gamma$ and
    $Q_{i,\beta} (\gamma') \neq \varnothing$.
    We see that if $j < r$,
    then $R_{i,\beta} (\gamma)$ consists of $2$ elements if
    $\smash{\sum_{j'=0}^{j-1}} \, x_{\gamma (j')} = 0$
    and $j \neq \smash{\sum_{j''=0}^{j'}} \, a'_{j''}$ for any $j'$,
    or $1$ element otherwise.
    Indeed, the two possible elements are
    $\langle \gamma, \min \gamma - 1 \rangle$ and
    $\langle \gamma, \max \gamma + 1 \rangle$,
    and one of them might not be allowed since
    $\bix$ must satisfy the set of inequalities $I'_\alpha \geq 0$
    for any $\alpha \in Q_{i,\beta}$, and \eqref{eq-def-q-i-beta}
    must be satisfied as well.
    Let $h (\gamma)$ be the number of $j' = 1, \dotsc, j-1$ such that
    $|R_{i,\beta} (\gamma_{(j')})| = 2$,
    and let $h' (\gamma)$ be the number of $j' = 1, \dotsc, j-1$ such that
    $|R_{i,\beta} (\gamma_{(j')})| = 2$ and
    $\gamma (j') = \min \gamma_{(j'-1)} - 1$.
    We claim that
    \begin{equation}
        \label{eq-l-gamma-equals-leibniz}
        l (\gamma) = \upL^{h (\gamma) + p}_{h' (\gamma) + j_0 + 1}\,.
    \end{equation}
    Indeed, if the length of $\gamma$ is $j = r$,
    then \eqref{eq-l-gamma-equals-leibniz}
    is a result of \eqref{eq-c-alpha-x-equals-leibniz},
    since in this case,
    $Q_{i,\beta} (\gamma)$ contains a unique element $\alpha$
    (it is non-empty by assumption),
    and $I'_\alpha (j') (\bix) = 0$
    if and only if either $|R_{i,\beta} (\gamma_{j'})| = 2$
    or $j' = \smash{\sum_{j\third=0}^{j\second}} a'_{j\third}$ for some $j\second$.
    If the length of $\gamma$ is $j < r$,
    then $|R_{i,\beta} (\gamma)|$ is $1$ or $2$.
    If it is $1$, then \eqref{eq-l-gamma-equals-leibniz}
    follows from the same equation for the unique element of $R_{i,\beta} (\gamma)$.
    If it is $2$, then \eqref{eq-l-gamma-equals-leibniz}
    follows from the identity
    \begin{equation}
        \upL^{h (\gamma) + p}_{h' (\gamma) + j_0 + 1} =
        \upL^{h (\gamma) + p + 1}_{h' (\gamma) + j_0 + 1} +
        \upL^{h (\gamma) + p + 1}_{h' (\gamma) + j_0 + 2}.
    \end{equation}
    This proves \eqref{eq-l-gamma-equals-leibniz},
    and taking $\gamma = \langle i \rangle$ proves
    \eqref{eq-sum-over-q-i-beta-equals-leibniz}.
    Finally, since the number of $\beta \in P_p$ with $\beta^{-1} (0) = j_0$
    is $\binom{p-1}{j_0}$, one has
    \begin{equation}
        \sum_{ \substack{ \beta \in P_p \\ \beta^{-1} (0) = j_0 } }
        \sum_{\alpha \in Q_{i,\beta}}
        c_\alpha (\bix) =
        \dbinom{p-1}{j_0} \cdot \upL^p_{j_0+1} =
        \frac{1}{p} = c_\Delta (\bix).
    \end{equation}
    This completes the proof of \eqref{eq-main-comb-lemma-strong}.
\end{proof}

\begin{corollary}
    \label{cor-comb-main}
    For any lattice $\Lambda \subset V_r$ and any
    $\gamma \in \check{H}_{\bullet} (\Mbpl_{>0}; \bbQ)$, one has
    \begin{multline}
        \sumbar_{\bid \in \Lambda} c_\Delta (\bid) \cdot
        [ \gamma, \, Z (\bid) ] = \\
        r \cdot \sumbar_{\bid \in \Lambda} c_\Delta (\bid) \cdot
        \bigl[ \bigl[ \dotsc \bigl[ \bigl[
            \gamma , \ 
            \inv{\Mss_{(1,d_0)}} \bigr] , \ 
            \inv{\Mss_{(1,d_1)}} \bigr] , \dotsc \bigr] , \ 
            \inv{\Mss_{(1,d_{r-1})}}
        \bigr],
    \end{multline}
    where $Z (\bid)$ is defined by \eqref{eq-def-z-d}.
\end{corollary}

\begin{proof}
    By Theorem~\ref{thm-c-delta-main}, the right-hand side equals
    \begin{align*}
        & \sum_{\alpha \in P_r} {} (-1)^{\alpha^{-1} (0)} \cdot
        \sumbar_{\bid \in \Lambda} (\alpha^{-1})^* (c_{\Delta}) (\bid) \cdot {} \\*
        & \hspace{4em} \bigl[ \bigl[ \dotsc \bigl[ \bigl[
            \gamma , \ 
            \inv{\Mss_{(1,d_0)}} \bigr] , \ 
            \inv{\Mss_{(1,d_1)}} \bigr] , \dotsc \bigr] , \ 
            \inv{\Mss_{(1,d_{r-1})}}
        \bigr] \\[1ex]
        = {} &
        \sum_{\alpha \in P_r} {} (-1)^{\alpha^{-1} (0)} \cdot
        \sumbar_{\bid \in \Lambda} c_{\Delta} (\bid) \cdot {} \\
        & \hspace{4em} \bigl[ \bigl[ \dotsc \bigl[ \bigl[
            \gamma , \ 
            \inv{\Mss_{(1,d_{\alpha (0)})}} \bigr] , \ 
            \inv{\Mss_{(1,d_{\alpha (1)})}} \bigr] , \dotsc \bigr] , \ 
            \inv{\Mss_{(1,d_{\alpha (r-1)})}}
        \bigr] \\[1ex]
        = {} &
        \sumbar_{\bid \in \Lambda} c_{\Delta} (\bid) \cdot 
        [ \gamma, \, Z (\bid) ],
        \numberthis
    \end{align*}
    where we applied Lemma~\ref{lem-iter-brackets-sum} in the last step.
    Note that when we have a finite sum of regularized sums
    with the same coefficient polytope function (which is $c_\Delta$ in this case),
    we may swap the order of the two summation signs.
\end{proof}

\subsection{Invariants as regularized sums}

We combine the results from the previous subsections,
and complete the proof of our main results,
expressing the invariants $\inv{\Mss_{(r,d)}}$ and
$\inv{\Mbss_{(r,d),1}}$ as regularized sums.

\begin{definition}
    \label{def-lambda-r-d}
    Let $r > 0$ and $d$ be integers.
    Let $\Lambda_{r,d} \subset \bbZ^r$ be the set of integer sequences
    $\bid = (d_0, \dotsc, d_{r-1})$ such that $d_0 + \cdots + d_{r-1} = d$.
    We also regard $\Lambda_{r,d}$ as a lattice in $V_r$,
    in the sense of Definition~\ref{def-cone-lattice}, via the projection
    \begin{align*}
        \numberthis
        \pi_{r,d} \colon \Lambda_{r,d} & \longrightarrow V_r \, , \\
        \bid & \longmapsto \bid - \frac{d}{r} \cdot (1, \dotsc, 1),
    \end{align*}
    so that if $c \in \PF (V_r; \bbQ)$ and $\bid \in \Lambda_{r,d}$\,,
    we may write $c (\bid)$ for $c (\pi_{r,d} (\bid))$.
\end{definition}

\begin{lemma}
    \label{lem-jigsaw}
    Let $\bid = (d_0, \dotsc, d_{r-1}) \in \Lambda_{r,d}$. Then
    \upshape
    \begin{equation}
        \sum_{ \leftsubstack[6em]{
            \\[-1ex]
            & 1 = k_0 < \cdots < k_m = r, \\[-2ex]
            & \text{such that }
            \tilde{d}_1 / \tilde{r}_1 \geq \cdots \geq \tilde{d}_m / \tilde{r}_m \geq d/r ,
            \\[-2ex]
            & \text{where } \textstyle
            (\tilde{r}_i, \tilde{d}_i) =
            \sum_{k_{i-1} \leq j < k_i} {} (1, d_j).
            \\[-1.5ex]
            & \text{Write } 
            \bid_i = (d_{k_{i-1}}, \dotsc, d_{k_i \mspace{2mu} -1})
            \in \Lambda_{\smash{\tilde{r}_i, \tilde{d}_i}}\,.
            \\[-1.5ex]
            & \text{Define } 
            0 = a_0 < \cdots < a_l = m
            \text{ such that for any } 0 < i < m,
            \\[-2ex]
            & \tilde{d}_i / \tilde{r}_i > \tilde{d}_{i+1} / \tilde{r}_{i+1} \text{ if and only if }
            i = a_j \text{ for some } 0 < j < l.
            \\[-2ex]
            & \text{Let $e$ be the number of $i > 0$ with }
            \tilde{d}_i / \tilde{r}_i = d/r
        } }
        \frac{(-1)^e \, B_e}{\prod_{i-1}^l {} (a_i - a_{i-1})!} \cdot
        \prod_{i=1}^m c_{\Delta_{\tilde{r}_i}} (\bid_i) =
        c_\Delta (\bid),
    \end{equation}
    \itshape
    where $B_e$ denotes the $e$-th Bernoulli number,
    and $c_{\Delta_{\tilde{r}_i}} \in \PF (V_{\tilde{r}_i}; \bbQ)$
    denotes the version of
    $c_\Delta$ with $V_{\tilde{r}_i}$ in place of $V_r$
    in Definition~\textnormal{\ref{def-c-delta}}.
\end{lemma}

\begin{proof}
    For $i = 1, \dotsc, r-1$,
    let $y_i = d_1 + \cdots + d_i$, and let
    \begin{equation}
        S_{\bid} = \{ (i, y_i) \mid i = 1, \dotsc, r-1 \} \subset \bbR^2.
    \end{equation}
    Let $\bar{S}_{\bid} \subset \bbR^2$ be the convex hull of
    $\{ (i, y_i - y) \mid i = 1, \dotsc, r-1, \ y \geq 0 \}$,
    and let $S'_{\smash{\bid}} = S_{\bid} \cap \partial \bar{S}_{\bid}$.
    Let
    \begin{align}
        A_{\bid} & = \{ i \mid (i, y_i) \in S'_{\smash{\bid}} \}, \\
        B_{\bid} & = \Bigl\{
            i \in A_{\bid} \Bigm|
            \text{there exists } i' < i < i\second 
            \text{ in } A_{\bid}
            \text{ with } \frac{y_i-y_{i'}}{i-i'} = \frac{y_{i\second}-y_i}{i\second-i}
        \Bigr\}.
    \end{align}
    For each non-zero term in the sum,
    we necessarily have
    \begin{equation}
        B_{\bid} \subset \{ k_0, \dotsc, k_m \} \subset A_{\bid}.
    \end{equation}
    Moreover, if there are non-zero terms,
    then $S_{\bid}$ must be contained in the region
    $(y_{r-1}-y)/((r-1)-x) \geq d/r$.
    Equivalently, the inequalities
    $(d_i + \cdots + d_{r-1})/(r-i) \geq d/r$ must hold for all
    $i=1, \dotsc, r-1$.
    This proves the lemma in the case $c_{\Delta} (\bid) = 0$.
    
    Write $B_{\bid} = \{ b_0, \dotsc, b_s \}$,
    where $b_0 < \cdots < b_s$\,,
    and for each $i = 1, \dotsc, s$,
    let $b'_{\smash{i,1}} < \cdots < b'_{\smash{i,t_i}}$ be the elements of $A_{\bid}$
    that lie between $b_{i-1}$ and $b_i$\,, where $t_i \geq 0$.
    
    Suppose that $c_\Delta (\bid) \neq 0$.
    If $e = 0$, then we have strict inequalities
    $(d_i + \cdots + d_{r-1})/(r-i) > d/r$ for all $i=1, \dotsc, r-1$.
    In this case, the lemma reduces to the combinatorial identity
    \begin{equation}
        \prod_{i=1}^s
        \sum_{ \substack{
            0 < m_i \leq t_i \\
            0 = c_0 < \cdots < c_{m_i} = t_i
        } }
        \frac{1}{m_i! \cdot \prod_{j=1}^{m_i} {} (c_j - c_{j-1})}
        = 1,
    \end{equation}
    which can be seen by looking at the coefficient of $x^{t_i}$
    in the identity
    \begin{equation}
        \exp \bigl( - \log (1-x) \bigr) = \frac{1}{1 - x},
    \end{equation}
    where $x$ is a formal variable.
    
    If $e > 0$, then we have to prove that
    \begin{align}
        \Biggl(
            \prod_{i=1}^{s-1}
            \sum_{ \substack{
                0 < m_i \leq t_i \\
                0 = c_0 < \cdots < c_{m_i} = t_i
            } } \mspace{-12mu}
            \frac{1}{m_i! \cdot \prod_{j=1}^{m_i} {} (c_j - c_{j-1})}
        \Biggr) \cdot \mspace{-12mu}
        \sum_{ \substack{
            0 < m_s \leq t_s \\
            0 = c_0 < \cdots < c_{m_s} = t_s
        } } \mspace{-12mu}
        \frac{(-1)^{m_s} B_{m_s}}{m_s! \cdot \prod_{j=1}^{m_s} {} (c_j - c_{j-1})}
        = \frac{1}{t_s + 1}.
        \raisetag{2ex}
    \end{align}
    The term $\prod_{i=1}^{s-1} {} ({\cdots})$ is $1$ as before,
    and the equation follows from examining the
    coefficient of $x^{t_s}$ in the identity
    \begin{equation}
        \frac{\log (1-x)}{\exp \bigl( \log (1-x) \bigr) - 1} =
        - \frac{\log (1-x)}{x}.
    \end{equation}
\end{proof}

\begin{theorem}
    \label{thm-main-restate}
    We have
    \begin{equation}
        \label{eq-inv-mss-equals-sum-over-lattice}
        \inv{\Mss_{(r,d)}} =
        \sumbar_{\bid \in \Lambda_{r,d}}
        \frac{1}{r} \, c_\Delta (\bid) \cdot Z (\bid),
    \end{equation}
    where $c_\Delta$ is as in Definition~\textnormal{\ref{def-c-delta}},
    and $Z (\bid)$ is defined by \eqref{eq-def-z-d}.
\end{theorem}

\begin{proof}
    The regularized sum is well-defined
    by Lemma~\ref{lem-z-multiplicative},
    where we choose the elements $Z (\bid; \biz) |_{z_0 = 0}$
    as representatives in cosets of $\im D$,
    as in Definition~\ref{def-reg-sum-of-residues}.

    We prove the result by induction on $r$.
    For $r = 1$, \eqref{eq-inv-mss-equals-sum-over-lattice}
    says that $\smash{\inv{\Mss_{(1,d)}}} = Z (d)$,
    which holds by definition.
    
    For $r > 1$, by Theorem~\ref{thm-reduction-main},
    and by the induction hypothesis, we have
    \begin{multline}
        \label{eq-inv-mss-equals-sum-of-brackets-of-z}
        \inv{ \Mss_{(r,d)} } = \hspace{-1em}
        \sum_{ \leftsubstack[6em]{
            \\
            & (r,d) = (1,d_0) + (r_1,d_1) + \cdots + (r_m,d_m), \\[-.6ex]
            & m \geq 1, \ r_i > 0 \text{ for all } i, \\[-.6ex]
            & \text{such that }
            d_1 / r_1 \geq \cdots \geq d_m / r_m \geq d/r .
            \\[-.6ex]
            & \text{Define } 
            0 = a_0 < \cdots < a_l = m
            \text{ such that for any } 0 < i < m,
            \\[-.6ex]
            & d_i / r_i > d_{i+1} / r_{i+1} \text{ if and only if }
            i = a_j \text{ for some } 0 < j < l.
            \\[-.6ex]
            & \text{Let $e$ be the number of $i > 0$ with }
            d_i / r_i = d/r
        } }
        \frac{(-1)^e \, B_e}{\prod_{i=1}^l (a_i - a_{i-1})!} \cdot {} \\
        \tau \Biggl( \ 
            \sumbar_{ \leftsubstack[5em]{
                \\[-.6ex]
                & d_i = d_{i,1} + \cdots + d_{i,r_i}\,, \\[-.6ex]
                & i = 1, \dotsc, m
            } } {}
            \biggl(
                \prod_{i=1}^m
                \frac{1}{r_i} \,
                c_{\Delta_{r_i}} (\bid_i)
            \biggr) \cdot
            [[ \dotsc [[
                \acute{Z} (d_0) , Z (\bid_1) ], Z (\bid_2) ],
                \dotsc], Z (\bid_m)
            ]
        \Biggr),
        \raisetag{2ex}
    \end{multline}
    where we write $\bid_i = (d_{i,1}, \dotsc, d_{i,r_i-1})$
    for $i = 1, \dotsc, m$,
    and $c_{\Delta_{r_i}}$ denotes the function
    $c_\Delta$ in $\PF (V_{r_i}; \bbQ)$.
    
    We then apply Corollary~\ref{cor-comb-main} $m$ times
    to each term of \eqref{eq-inv-mss-equals-sum-of-brackets-of-z},
    each time to the Lie bracket with one of the $Z (\bid_i)$.
    This shows that
    \begin{multline}
        \label{eq-inv-mss-equals-sum-of-z-acute}
        \inv{ \Mss_{(r,d)} } = \hspace{-1em}
        \sum_{ \leftsubstack{
                \\
                & (r,d) = (1,d_0) + (r_1,d_1) + \cdots + (r_m,d_m), \\[-.6ex]
                & m \geq 1, \ r_i > 0 \text{ for all } i, \\[-.6ex]
                & \text{such that }
                d_1 / r_1 \geq \cdots \geq d_m / r_m \geq d/r .
                \\[-.6ex]
                & \text{Define } 
                0 = a_0 < \cdots < a_l = m
                \text{ such that for any } 0 < i < m,
                \\[-.6ex]
                & d_i / r_i > d_{i+1} / r_{i+1} \text{ if and only if }
                i = a_j \text{ for some } 0 < j < l.
                \\[-.6ex]
                & \text{Let $e$ be the number of $i > 0$ with }
                d_i / r_i = d/r
        } }
        \frac{(-1)^e \, B_e}{\prod_{i=1}^l (a_i - a_{i-1})!} \cdot {} \\
        \tau \Biggl( \ 
            \sumbar_{ \leftsubstack[6.4em]{
                & d_i = d_{i,1} + \cdots + d_{i,r_i}\,, \\[-.6ex]
                & i = 1, \dotsc, m
            } } {}
            \biggl( \prod_{i=1}^m c_{\Delta_{r_i}} (\bid_i) \biggr) \cdot
            \acute{Z} (d_0, \bid_1, \dotsc, \bid_m)
        \Biggr).
    \end{multline}
    
    For a fixed $\bid \in \Lambda_{r,d}$\,,
    by Lemma~\ref{lem-jigsaw},
    the sum of all terms in
    \eqref{eq-inv-mss-equals-sum-of-z-acute} satisfying
    $(d_0, \bid_1, \dotsc, \bid_m) = \bid$
    equals $\tau (c_\Delta (\bid) \cdot \acute{Z} (\bid))$.
    This shows that 
    \begin{equation}
        \label{eq-inv-mss-equals-sum-of-z-acute-2}
        \inv{ \Mss_{(r,d)} } = 
        \tau \biggl( \,
            \sumbar_{ \bid \in \Lambda_{r,d} }
            c_\Delta (\bid) \cdot \acute{Z} (\bid)
        \biggr).
    \end{equation}
    
    Finally, applying Corollary~\ref{cor-comb-main} once again,
    we see that
    \begin{equation}
        \sumbar_{ \bid \in \Lambda_{r,d} }
        c_\Delta (\bid) \cdot
        \acute{Z} (\bid) =
        \frac{1}{r}
        \sumbar_{ \bid \in \Lambda_{r,d} }
        c_\Delta (\bid) \cdot
        \bigl[ -\upe^{((0,0),1)}, \, Z (\bid) \bigr],
    \end{equation}
    and the result follows from \eqref{eq-inv-mss-equals-sum-of-z-acute-2}.
\end{proof}

\begin{bproof}[Proof of Theorem~\ref{thm-inv-mss-main-as-reg-sum}]
    Part \ref{item-inv-mss-main-as-reg-sum}
    was restated and proved as
    Theorem~\ref{thm-main-restate}.
    
    For \ref{item-inv-mbss-main-as-reg-sum},
    by Theorem~\ref{thm-wcf-2}, one has
    \begin{multline}
        \inv{\Mbss_{(r,d),1}} =
        \sum_{ \leftsubstack{
            \\[-1ex]
            & (r,d) = (r_1,d_1) + \cdots + (r_m,d_m), \\[-.6ex]
            & m \geq 1, \ r_i > 0, \ 
            d_i / r_i = d / r \text{ for all } i
        } }
        \frac {(-1)^m} {m!} \cdot
        \bigl[ \bigl[ \dotsc \bigl[ \bigl[
            \upe^{((0,0),1)}, \ 
            \inv{\Mss_{(r_1,d_1)}} \bigr], \\[-5.5ex]
            \shoveright{\inv{\Mss_{(r_2,d_2)}} \bigr] , \dotsc \bigr] , \ } \\
            \inv{\Mss_{(r_m,d_m)}}
        \bigr] .
    \end{multline}
    By Theorem~\ref{thm-main-restate},
    and by an argument similar to how we obtained
    \eqref{eq-inv-mss-equals-sum-of-z-acute}, we see that
    \begin{multline}
        \label{eq-inv-mbss-equals-sum-of-z-acute}
        \inv{\Mbss_{(r,d),1}} =
        \sum_{ \leftsubstack{
            \\[-1ex]
            & (r,d) = (r_1,d_1) + \cdots + (r_m,d_m), \\[-.6ex]
            & m \geq 1, \ r_i > 0, \ 
            d_i / r_i = d / r \text{ for all } i
        } }
        \frac {(-1)^{m-1}} {m!} \cdot {} \\
        \sumbar_{ \leftsubstack[6.4em]{
            & d_i = d_{i,1} + \cdots + d_{i,r_i}\,, \\[-.6ex]
            & i = 1, \dotsc, m
        } } {}
        \biggl( \prod_{i=1}^m c_{\Delta_{r_i}} (\bid_i) \biggr) \cdot
        \acute{Z} (\bid_1, \dotsc, \bid_m),
    \end{multline}
    where we write $\bid_i = (d_{i,1}, \dotsc, d_{i,r_i-1})$
    for $i = 1, \dotsc, m$,
    and $c_{\Delta_{r_i}}$ denotes the function
    $c_\Delta$ in $\PF (V_{r_i}; \bbQ)$.
    
    For fixed $\bid \in \Lambda_{r,d} \cap \supp c_{\Delta}$\,,
    collect all terms in \eqref{eq-inv-mbss-equals-sum-of-z-acute}
    that involve $\acute{Z} (\bid)$.
    Let 
    \begin{equation}
        I = \Bigl\{ i = 1, \dotsc, r-1 \Bigm|
        \frac{d_i + \cdots + d_{r-1}}{r-i} = \frac{d}{r} \Bigr\},
    \end{equation}
    where the $d_i$ are the components of $\bid$.
    We see that the coefficient of $\acute{Z} (\bid)$
    in \eqref{eq-inv-mbss-equals-sum-of-z-acute} is
    \begin{align*}
        &
        \sum_{ \leftsubstack{
            \\[-2ex]
            & r = r_1 + \cdots + r_m, \\[-.6ex]
            & r_i > 0, \ 
            \tilde{r}_i \in I \text{ for all } 0<i<m, \\[-.6ex]
            & \text{where } \tilde{r}_i = r_1 + \cdots + r_i
        } }
        \frac {(-1)^{m-1}} {m!} \cdot
        \prod_{i=1}^m \frac{1}{|I \cap \{ \tilde{r}_{i-1} + 1, \dotsc, \tilde{r}_i - 1\}| + 1} \\
        = {} &
        \sum_{ \leftsubstack{
            \\[-2ex]
            & |I| + 1 = a_1 + \cdots + a_m, \\[-.6ex]
            & a_i > 0
        } }
        \frac {(-1)^{m-1}} {m!} \cdot \prod_{i=1}^m \frac{1}{a_i} \\
        = {} &
        \begin{cases}
            1, & |I| = 0, \\
            0, & |I| > 0.
        \end{cases}
        \numberthis
    \end{align*}
    Therefore, we see that $\acute{Z} (\bid)$ contributes to
    \eqref{eq-inv-mbss-equals-sum-of-z-acute}
    if and only if $\bid$ is in the interior of $\supp c_\Delta$
    (and not on the boundary),
    so that by an argument similar to the proof of Theorem~\ref{thm-main-restate},
    we obtain
    \begin{equation}
        \inv{\Mbss_{(r,d),1}} =
        \sumbar_{ \leftsubstack{
            \\[-3ex]
            & d = d_0 + \cdots + d_{r-1} \\[-1.6ex]
            & (d_i + \cdots + d_{r-1})/(r-i) > d/r, \ i = 1, \dotsc, r-1
        } } 
        \acute{Z} (\bid).
        \hspace{4em}
    \end{equation}
\end{bproof}

\begin{bproof}[Proof of Theorem~\ref{thm-inv-mss-main}]
    To prove \eqref{eq-main-explicit}, 
    by Theorem~\ref{thm-main-restate},
    it suffices to show that
    \begin{equation}
        \label{eq-xi-of-final-result}
        \xi \biggl( \upe^{(r,d)} \cdot
            \prod_{i=0}^{r-1} {}
            \bigl[
                \exp ( z_i \, D_{(1, d_i)} ) \, \sigma (-s_{1,2,2})
            \bigr]
        \biggr) = \upe^{(r,d)} \cdot
        \prod_{i=0}^{r-1} {}
        \bigl[
            \exp ( \tilde{z}_i \, D_{(1, d_i)} ) \, \sigma (-s_{1,2,2})
        \bigr].
    \end{equation}
    This is because the $\partial / \partial s_{1,0,1}$ of the right-hand side is $0$,
    and
    \begin{multline}
        \prod_{i=0}^{r-1} {}
        \bigl[
            \exp ( \tilde{z}_i \, D_{(1, d_i)} ) \, \sigma (-s_{1,2,2})
        \bigr] = \\
        \exp \Bigl( -\frac{z_0 + \cdots + z_{r-1}}{r} \, D_{(r,d)} \Bigr)
        \prod_{i=0}^{r-1} {}
        \bigl[
            \exp ( z_i \, D_{(1, d_i)} ) \, \sigma (-s_{1,2,2})
        \bigr],
    \end{multline}
    so the difference of the two sides lies in $\im D_{(r,d)}$.
    Therefore, \eqref{eq-xi-of-final-result} must hold by Lemma~\ref{lem-xi-well-def}.
    
    Next, let us prove \eqref{eq-main-explicit-summed}
    from \eqref{eq-main-explicit}.
    For each subset $I = \{ i_1, \dotsc, i_m \} \subset \{ 1, \dotsc, r-1 \}$,
    with $i_1 < \cdots < i_m$, we define a subset
    $\Delta_I \subset V_r$ (see Definition~\ref{def-vector-space-vr})
    by
    \begin{equation}
        \Delta_I = \biggl\{
            \bix \in V_r \biggm|
            \begin{aligned}
                & x_i + \cdots + x_{r-1} \geq 0, \quad i = 1, \dotsc, r-1, \\
                & x_i + \cdots + x_{r-1} = 0, \quad i \in I
            \end{aligned}
        \biggr\}.
    \end{equation}
    Then we observe that
    \begin{equation}
        \label{eq-c-delta-equals-chi-delta-i}
        c_{\Delta} = \sum_{|I| = m} \frac{(-1)^{m-1}}{m} \, \chi_{\Delta_I}\,,
    \end{equation}
    where $\chi_{\Delta_I} \in \PF (V_r; \bbQ)$ is
    the characteristic function of the cone $\Delta_I$.
    
    Now, consider the lattice $\Lambda_{r,d}$ of Definition~\ref{def-lambda-r-d},
    which we are summing over.
    Then $\Delta_I \cap \pi_{r,d} (\Lambda_{r,d})$ is non-empty if and only if
    for each $i \in I$, one has $i d / r \in \bbZ$.
    If this is the case, then $\Delta_I \cap \pi_{r,d} (\Lambda_{r,d})$
    is a simple sector of rank $r-1-m$,
    whose apex is the point
    \begin{equation}
        \Bigl(
            \Bigl\lfloor \frac{d}{r} \Bigr\rfloor, \ 
            \Bigl\lfloor \frac{2d}{r} \Bigr\rfloor -
            \Bigl\lfloor \frac{d}{r} \Bigr\rfloor, \ \dotsc, \ 
            d -
            \Bigl\lfloor \frac{(r-1)d}{r} \Bigr\rfloor
        \Bigr) \in \Lambda_{r,d}\,,
    \end{equation}
    and the basis vectors are
    $\bie_i - \bie_{i-1}$ for $i \in \{ 1, \dotsc, r-1 \} \setminus I$.
    Finally, using Lemma~\ref{lem-exp-z-d-expansion} to write
    \begin{multline}
        \prod_{i=0}^{r-1} {}
        \bigl[
            \exp ( \tilde{z}_i \, D_{(1, d_i)} ) \, \sigma (-s_{1,2,2})
        \bigr] = \\
        \exp \biggl(
            \sum_{l=1}^{\infty}
            \sum_{i=0}^{r-1}
            \frac{\tilde{z}_i^l}{l!} 
            \, \Bigl( d_i + 
            \Bigl\lfloor \frac{id}{r} \Bigr\rfloor - 
            \Bigl\lfloor \frac{(i+1)d}{r} \Bigr\rfloor \Bigr) \,
            s_{1,2,l+1}
        \biggr) \cdot {} \\
        \prod_{i=0}^{r-1} {}
        \biggl[
            \exp \Bigl( \tilde{z}_i \,
                D_{\left(1, \,
                    \lfloor\mspace{-3mu} \frac{(i+1)d}{r} \mspace{-3mu}\rfloor -
                    \lfloor\mspace{-3mu} \frac{id}{r} \mspace{-3mu}\rfloor
                \right)}
            \Bigr) \, \sigma (-s_{1,2,2})
        \biggr],
    \end{multline}
    we can evaluate the regularized sum in \eqref{eq-main-explicit}
    using \eqref{eq-c-delta-equals-chi-delta-i}.
    This proves \eqref{eq-main-explicit-summed}.
    
    The proof of \eqref{eq-main-pairs-explicit}
    is analogous to that of \eqref{eq-main-explicit}.
    
    For \eqref{eq-main-pairs-explicit-summed},
    define $\Delta' \subset V_r$ by
    \begin{equation}
        \Delta' = \{ \bix \in V_r \mid
        x_i + \cdots + x_{r-1} > 0, \quad i = 1, \dotsc, r-1 \}.
    \end{equation}
    Then $\Delta' \cap \pi_{r,d} (\Lambda_{r,d})$
    is a simple sector of rank $r-1$, whose apex is the point
    \begin{equation}
        \Bigl(
            \Bigl\lceil \frac{d}{r} \Bigr\rceil - 1, \ 
            \Bigl\lceil \frac{2d}{r} \Bigr\rceil -
            \Bigl\lceil \frac{d}{r} \Bigr\rceil, \ \dotsc, \ 
            \Bigl\lceil \frac{(r-2)d}{r} \Bigr\rceil -
            \Bigl\lceil \frac{(r-1)d}{r} \Bigr\rceil, \ 
            d -
            \Bigl\lceil \frac{(r-1)d}{r} \Bigr\rceil + 1
        \Bigr) \in \Lambda_{r,d}\,,
    \end{equation}
    and the basis vectors are
    $\bie_i - \bie_{i-1}$ for $i = 1, \dotsc, r-1$.
    The rest of the proof is analogous to that of
    \eqref{eq-main-explicit-summed}.
\end{bproof}

\phantomsection
\addcontentsline{toc}{section}{References}
\sloppy
\printbibliography

\par\noindent\rule{0.38\textwidth}{0.4pt}
{\par\noindent\small
\hspace*{2em}Mathematical Institute, University of Oxford, Oxford OX2 6GG, United Kingdom.\\[-2pt]
\hspace*{2em}Email: \texttt{bu@maths.ox.ac.uk}
}

\end{document}